\documentclass[preprint,12pt]{elsarticle}

\usepackage[margin=0.7in]{geometry}
\usepackage{amssymb}
\usepackage{amsthm}
\usepackage{subfig}
\usepackage{color}
\usepackage{graphicx,graphbox}
\usepackage{natbib}
\usepackage{amsmath}
\usepackage{ulem}

\newcommand{\revMR}[1]{{\color{black}#1}}
\newcommand{\revMRR}[1]{{\color{black}#1}}
\newcommand{\revAB}[1]{{\color{black}#1}}
\newcommand{\revRP}[1]{{\color{black}#1}}
\newcommand{\revMC}[1]{{\color{black}#1}}

\begin{document}

	\begin{frontmatter}

	\newpage


	\title{Extrapolated DIscontinuity Tracking for complex  2D  shock interactions}

		
	\author[lab1]{Mirco Ciallella }
	\author[lab1]{Mario Ricchiuto }
	\author[lab2]{Renato Paciorri }
	\author[lab3]{Aldo Bonfiglioli}
	\address[lab1]{Team CARDAMOM, INRIA, Univ. Bordeaux, CNRS, Bordeaux INP, IMB, UMR 5251, 200 Avenue de la Vieille Tour, 33405 Talence cedex, France}
	\address[lab2]{Dip. di Ingegneria Meccanica e Aerospaziale, Universit\`a di Roma ``La Sapienza'', Via Eudossiana 18, 00184 Rome, Italy}
	\address[lab3]{Scuola di Ingeneria  - Universit\`a degli Studi della Basilicata, Viale dell'Ateneo Lucano 10, 85100 Potenza, Italy}
	\begin{abstract}
 	A new shock-tracking technique that avoids 
        re-meshing the computational grid around the moving shock-front 
        was recently proposed by the authors~\cite{EST2020}. 
        \revMC{The method combines the unstructured shock-fitting~\cite{Paciorri1} approach, developed in the last decade
        by some of the authors, with ideas coming from embedded boundary methods.
        In particular, second-order extrapolations based on Taylor series expansions are employed to trasfer the solution and retain high order of accuracy.}
	This paper describes the basic idea behind the new method and further algorithmic improvements 
        which make the extrapolated Discontinuity Tracking Technique ($eDIT$) capable of dealing with complex shock-topologies featuring
	shock-shock and shock-wall interactions occuring in steady problems.
        \revMC{This method paves the way to a new class of shock-tracking techniques truly independent on the mesh structure and flow solver.} 
	Various test-cases are included to  prove the potential of the method,  demonstrate the  key features of the methodology, \revMC{and thoroughly evaluate several technical aspects related to 
              the extrapolation from/onto the shock, and their impact on accuracy, and  conservation}.

        \end{abstract}

	\begin{keyword}
	Shock-tracking \sep shock-fitting 
        \sep  unstructured grids \sep embedded boundary  
        \sep shock interactions \sep Taylor expansions 
	\end{keyword}

	\end{frontmatter}

	\section{Introduction}\label{introduction}
\revRP{The numerical models nowadays used in the simulation of compressible high-speed flows within
the continuum framework can be cast into two main categories:
\revAB{shock-capturing (SC) and shock-fitting (SF).}
\revAB{SC can be traced back to the 1950s paper by VonNeumann and Richtmyer~\cite{doi:10.1063/1.1699639}
and lays its foundations in the mathematical theory of weak solutions~\cite{lax1973hyperbolic}.
The key ingredients of SC discretizations are very well known in the CFD community and
need not to be recalled here. 
SF techniques, on the contrary, rely on a very different standpoint},
which consists in identifying the shocks as lines within a 
two-dimensional (2D) flow-field (surfaces in 3D)
and computing their motion, upstream and downstream states according to the Rankine-Hugoniot 
\revAB{jump relations}.
\revAB{Fitted-shocks made their first appearance in two NACA technical
reports} written by Emmons~\cite{Emmons1,Emmons2} in the 1940s, but
it is between the late 1960s and the end of the 1980s that the technique has been extensively developed by
Gino Moretti and collaborators under the name ``Shock-Fitting'' 
and, starting in the early 1980s, by James Glimm and collaborators under the name ``Front-Tracking''.
Moretti's SF techniques (\revAB{both the ``boundary'' and ``floating'' variants)} were designed 
for structured-grid solvers, \revAB{the only kind of CFD codes in use at that time.
The development of structured-grid SF codes turned out to be
algorithmically} \revAB{fairly} complex, especially 
when \revAB{dealing with} flows \revAB{featuring shock-shock and shock-boundary} interactions~\cite{Moretti2}.
\revAB{Algorithmic complexity contributed to make CFD practitioners
``shy away from shock-fitting''~\cite{moretti2002thirty-six} 
to embrace the myth of ``one code for all flows'' that has so much
contributed to the popularity of SC discretizations}.}

\revAB{Starting in the late 1980s, however}, the CFD community has shown increasing interest towards unstructured
meshes \revAB{because of their} ability to easily \revAB{describe}
complex geometries and the possibility of locally changing the mesh size to adapt it to the local flow features.
Taking advantage of the greater flexibility offered by unstructured meshes, 
Paciorri and Bonfiglioli~\cite{Paciorri1} developed a new unstructured SF technique \revAB{suitable for
vertex-centered solvers operating with simplicial elements (triangles in 2D  and tetrahedra in 3D)}. 
Their approach alleviated many of the
\revAB{algorithmic difficulties incurred} by the SF techniques within the structured-grid setting.
In recent years, the unstructured SF technique originally proposed in~\cite{Paciorri1}
has been further improved making it capable of dealing with multiple interacting
discontinuities~\cite{Paciorri2}
and unsteady flows~\cite{Paciorri4,Campoli2017} in 2D, as well as 3D flows~\cite{Paciorri3}, thus opening 
a new route for simulating compressible ``shocked'' flows.
In particular, not only shocks and contact discontinuities are fitted, but also their interaction points,
\revAB{for example} the triple point that arises in Mach reflections~\cite{Paciorri5}.
\revAB{Most of the aforementioned developments have been included in an
open-source platform~\cite{campoli2021undifi}
with the goal of fostering even more the interest in such methods}.

\revRP{Thanks to these new developments the interest in shock-fitting/tracking methods has seen a 
resurgence with the proposal of several implicit and 
explicit algorithms~\cite{zou2017shock,corrigan2019moving,zahr2020implicit,zou2021,chang2019adaptive,chen2020convergence,huang2021robust,luke2021novel,zahr2021high}. 
A limitation of this technique is that it heavily relies on the flexibility of triangular and tetrahedral grids
to locally produce a fitted unstructured grid around the discontinuities.
}

The initial idea \revAB{that motivated the present work} comes from the similarity between the constraints 
arising from SF, as conceived by Paciorri and Bonfiglioli in~\cite{Paciorri1}, 
and those related to the construction of 
boundary-fitted grids for simulating flows around complex geometries. 
\revAB{More precisely, the SF technique~\cite{Paciorri1} requires that at each time-step the 
computational grid is locally re-meshed to follow the moving shock-front.
This approach turns out to be best suited to vertex-centered CFD solvers, 
but it can be problematic if a different type of data representation is used.
The $eDIT$ approach, which has been described by the authors in~\cite{EST2020} 
and is further developed in this paper, allows to
overcome this limitation, because it avoids re-meshing by using ideas borrowed from}
immersed and embedded boundary methods developed 
to allow a flexible management of complex geometries. These two approaches rely on different philosophies, however.\\
Immersed methods are based on an extension of the flow equations outside the physical domain 
(typically within solid bodies). This extension is formulated using some smooth
approximation of the Dirac delta function to localize the
boundary, as well as to impose the boundary condition. These methods are relatively old, and based on the original
ideas of Peskin~\cite{IB0}. Finite element and unstructured mesh extensions for elliptic PDEs 
as well as for incompressible, and compressible flows have been discussed in~\cite{IB1,IB2,IB3}.\\
Embedded methods, on the other hand, solve the PDEs only in the physical domain, while
replacing the exact boundary with some less accurate approximation, combined
with some weak enforcement of the boundary condition. There is a certain number of techniques
to perform this task, which go from the combination of XFEM-type methods with
penalization or Nitsche's type approaches~\cite{EM0},
to several types of cut finite element methods with improved stability~\cite{EM1,EM2}, to
approximate domain methods such as the well known ghost-fluid method~\cite{EM3,EM4},
and the more recent shifted boundary method (SBM)~\cite{Scovazzi1,Scovazzi2}.\\
As in the latter, we impose modified conditions on surrogate shock-manifolds, acting as boundaries 
between the shock-upstream and shock-downstream regions. The values of the flow variables imposed
on these surrogate boundaries are extrapolated from the tracked shock-front, accounting for
the non-linear jump and wave propagation conditions, as done in the unstructured SF approach of~\cite{Paciorri1}.
As in shock-fitting and front-tracking methods~\cite{Glimm2016}, the shock-front is discretized
by an independent lower-dimensional mesh, and its position, as well as
the position of the two surrogate boundaries, are themselves part of the computational result. \\

In this paper, we focus on the \revAB{new features coded} in the algorithm
to \revAB{enable} the computation of shock-shock and shock-wall interactions. In particular, when several
discontinuities are \revAB{mutually} interacting, or interact with a solid boundary,
the \revAB{following issues} must be addressed:
\begin{enumerate}
\item \revAB{identification of the various sub-domains obtained 
when the discontinuities cut through the computational domain};
\item \revAB{calculation of the velocity of those points where different discontinuities mutually interact
(for example the triple point in a Mach reflection)
or reach the boundary};
\item \revAB{the algorithm used to transfer the dependent variables between the shock-front
and the surrogate boundaries must be applicable
 also when the shape of the surrogate boundaries becomes very complex, which typically occurs
in the neighborhood of an interaction}.
\end{enumerate}

This paper is organized as follows: \S~\ref{generalities} \revAB{introduces the governing conservation equations}; 
\S~\ref{algorithm} \revAB{describes} the structure of the algorithm; \S~\ref{results} focuses on the various test-cases considered; 
finally in \S~\ref{conclusion} we draw some conclusions.

\section{Generalities}\label{generalities}

We consider the numerical approximation of solutions of the  steady limit of the Euler equations 
reading:
\begin{equation}\label{eq:euler0}
\partial_t \mathbf{U}+\nabla\cdot\mathbf{F} =0 \quad \text{in}\quad\Omega\subset\mathbb{R}^2
\end{equation}
with conserved variables and fluxes given by:
\begin{equation}\label{eq:euler0a}
\mathbf{U}=\left[
\begin{array}{c}
\rho \\ \rho \mathbf{u}  \\ \rho E
\end{array}
\right]\;,\;\;
\mathbf{F}= \left[
\begin{array}{c}
\rho \mathbf{u}\\ \rho \mathbf{u} \otimes \mathbf{u} + p \mathbb{I} \\ \rho H \mathbf{u}
\end{array}
\right] 
\end{equation}
having denoted by $\rho$ the mass density, by $\mathbf{u}$ the velocity, by $p$ the pressure, and with $E= e +  \mathbf{u}\cdot\mathbf{u}/2$ the specific total energy, $e$ being the specific internal energy.
Finally, the total specific enthalpy is  $H=h+ \mathbf{u}\cdot\mathbf{u}/2$, with  $h=e+p/\rho$ the specific enthalpy. For simplicity in this paper we  work with the classical perfect gas
equation of state:
\begin{equation}\label{eq:EOS}
p=(\gamma-1)\rho e
\end{equation}
with $\gamma$ the constant (for a perfect gas) ratio of specific heats. However,
note that the method discussed allows in principle to handle any other type of gas, see e.g.~\cite{pepe2015unstructured}.

In all applications involving high-speed flows, solutions of~(\ref{eq:euler0})  are  only piecewise continuous.
In $d$ space dimensions, discontinuities are represented by $d-1$ manifolds governed by the  well known Rankine-Hugoniot jump conditions reading:
\begin{equation}\label{eq:euler1}
	[\![\mathbf{F}]\!]\cdot\mathbf{n} =\omega[\![\mathbf{U}]\!]
\end{equation}
having denoted by $\mathbf{n}$ the local \revAB{unit vector normal} to the shock, by $[\![\cdot]\!]$ the corresponding jump of a quantity across the discontinuity,
and with $\revAB{\omega}$ the normal component of the shock speed.\\
As discussed in the introduction, the method proposed exploits ideas from two different approaches:
the unstructured shock-fitting method~\cite{Paciorri1} and subsequent works~\cite{Paciorri2,Paciorri3,Paciorri4}; the shifted boundary method by~\cite{Scovazzi1} and subsequent works~\cite{Scovazzi2,Scovazzi3}.\\
\revMR{Although the paper focuses on triangle based nodal solvers, the method proposed  is also \revAB{applicable} to \revRP{structured and  unstructured} cell-centered solvers~\cite{Assonitis2021}}.

%
%
\section{Extrapolated DIscontinuity Tracking ($eDIT$)}\label{algorithm}
%
\subsection{\revMR{Generalities}}\label{setting}

To illustrate the algorithmic features of the $eDIT$ method, let us consider a two-dimensional domain and a discontinuity front 
crossing the domain at a given time $t$.
\revMR{As shown \revAB{in Fig.~\ref{pic1}}, the front is described  using a collection of edges  
whose endpoints are marked by squares.
Throughout the paper  these two entities will be referred to as 
 \textit{front-edges} or \textit{discontinuity-edges}, and \textit{front-points} or
 \textit{discontinuity-points} respectively, and together they constitute what will be referred to as 
\textit{front-mesh} or \textit{discontinuity-mesh}.}
\revMR{The computational domain itself is discretized by means of
a background  mesh. 
As already said, we consider here solvers based on a nodal variable arrangement on triangular grids,
but the extrapolation proposed readily applies to cell-centered solvers}. \revAB{Figure~\ref{pic1} clearly shows that}
 the position of the discontinuity-points is completely 
independent of the location of the grid-points of the \textit{background mesh}. 
\revAB{Regardless of} \revMR{the method used to solve the Euler equations on the background mesh, the representation of the flow variables across the discontinuity
is discontinuous. In our case, two  values of the flow variables are stored in each front-\revAB{point}. Other arrangements (edge based for example) are of course possible.}
\revMR{The flow solvers considered here are based on a continuous nodal approximation}, 
\revAB{which amounts to store one solution value at each grid-point of the background mesh}. 
\revMR{As already said, this is not  essential for the method developed.}

 Assume that at time \emph{t} the solution is known at all grid- and discontinuity-points. 
The computation of the subsequent time level \emph{t+$\Delta$t} \revMR{using  the $eDIT$ method}  can be split into several steps that will 
be described in detail in the following sub-sections.
%
\subsection{Geometrical setting}\label{cell removal}
\begin{figure}
\begin{center}
\subfloat[Front-mesh laid on top of the background-mesh.]{%
\includegraphics[width=0.35\textwidth]{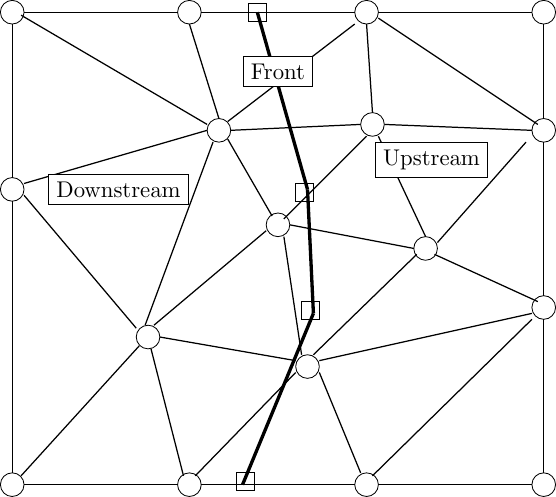}\label{pic1}}\qquad
\subfloat[Front-mesh, computational-mesh and surrogate boundaries.]{%
\includegraphics[width=0.35\textwidth]{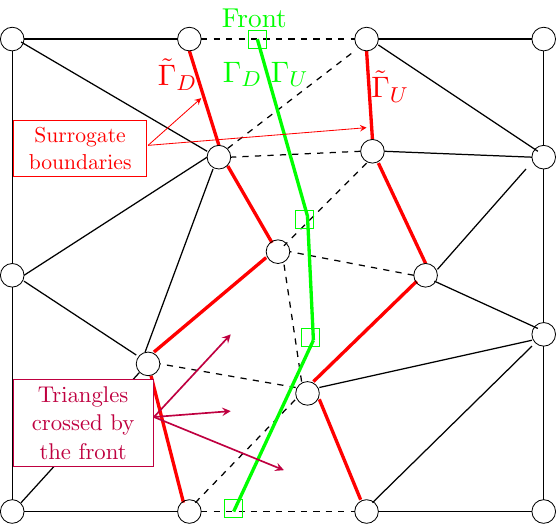}\label{pic2}}
\caption{The computational-mesh is obtained by removing those cells of the background mesh that are crossed by the front-mesh.}
\end{center}
\end{figure}

%
The first step consists in \revMR{flagging} the triangles crossed by the front.  
\revMR{As shown on  Fig.~\ref{pic2}, this 
leads to the creation of a cavity of flagged elements enclosing the discontinuity.}
\revMR{Here,} in contrast to the unstructured technique proposed in~\cite{Paciorri1},   
\revMR{or to more classical methods \revAB{such as the} boundary shock fitting~\cite{Salas}, the mesh within or in proximity of this cavity} is not \revMR{modified to be conformal
with the discontinuity. This cavity separates  two sub-domains, and allows to  define two surrogate boundaries, which are the intersections between the boundaries of the sub-domains and the boundaries of the front cavity.} 
\revMR{The mesh of the sub-domains separated by the front cavities will be referred to as the}
 \textit{computational mesh}. 
\revAB{The computational mesh is identical to the background mesh, except for the removal of the flagged elements enclosing the discontinuity.}
\revMR{In the simplest setting of a single discontinuity,} 
the upstream and downstream surrogate boundaries, drawn using red lines in Fig.~\ref{pic2}, 
will be called $\tilde{\Gamma}_U$ and $\tilde{\Gamma}_D$. \revMR{In the following we also refer to them as surrogate-discontinuities.}
The discontinuity-boundary, 
representing the actual front position, will be referred to as $\Gamma$ 
and its upstream and downstream sides as $\Gamma_U$ and $\Gamma_D$, respectively. 
Contrary to $\tilde{\Gamma}_U$ and $\tilde{\Gamma}_D$, $\Gamma_U$ and $\Gamma_D$ are superimposed
\revMR{and, of course, coincide with $\Gamma$}.

\revAB{Compared to its original version~\cite{EST2020},
the current algorithm is capable of dealing with several discontinuity-fronts which
implies that the computational domain is split into several, disjoint, sub-domains}.
%
%

\revMR{To} \revAB{apply the Rankine-Hugoniot jump relations~\eqref{eq:euler1} within
each pair of discontinuity-points}, \revMR{we need to define the tangent  and normal directions,}
$\boldsymbol{\tau}$ and $\mathbf{n}$, \revMR{ in correspondence of each of these points.} 
The computation of these vectors is carried out using finite-difference formulae
which involve the coordinates of the front-point itself and those of the neighboring front-points
that belong to its \revAB{range} of influence.
Specifically, centred finite difference
formulae will be used when the downstream state is subsonic whereas upwind formulae are employed whenever
the downstream state is supersonic. Further details are given \revRP{in~\cite{EST2020,assonitis2021numerical}}.
The unit vector $\boldsymbol{n}$ is conventionally assumed to be directed from
the high entropy (downstream) side towards the low entropy (upstream) side of the discontinuity,
\revAB{see Fig.~\ref{inter_alg}}.
It is also convenient to express the components of the velocity vector in the
$(\boldsymbol{n},\boldsymbol{\tau})$ reference frame which is locally attached to the discontinuity:
\begin{equation}
u_n = \mathbf{u} \cdot \mathbf{n} \;\;\; \text{and} \;\;\; u_\tau = \mathbf{u} \cdot \boldsymbol{\tau}
\end{equation}
\subsection{Definition of \revMR{the sub-domains} } 
\begin{figure}
\centering
\includegraphics[width=0.5\textwidth,trim={1cm 1cm 0.5cm 1cm},clip]{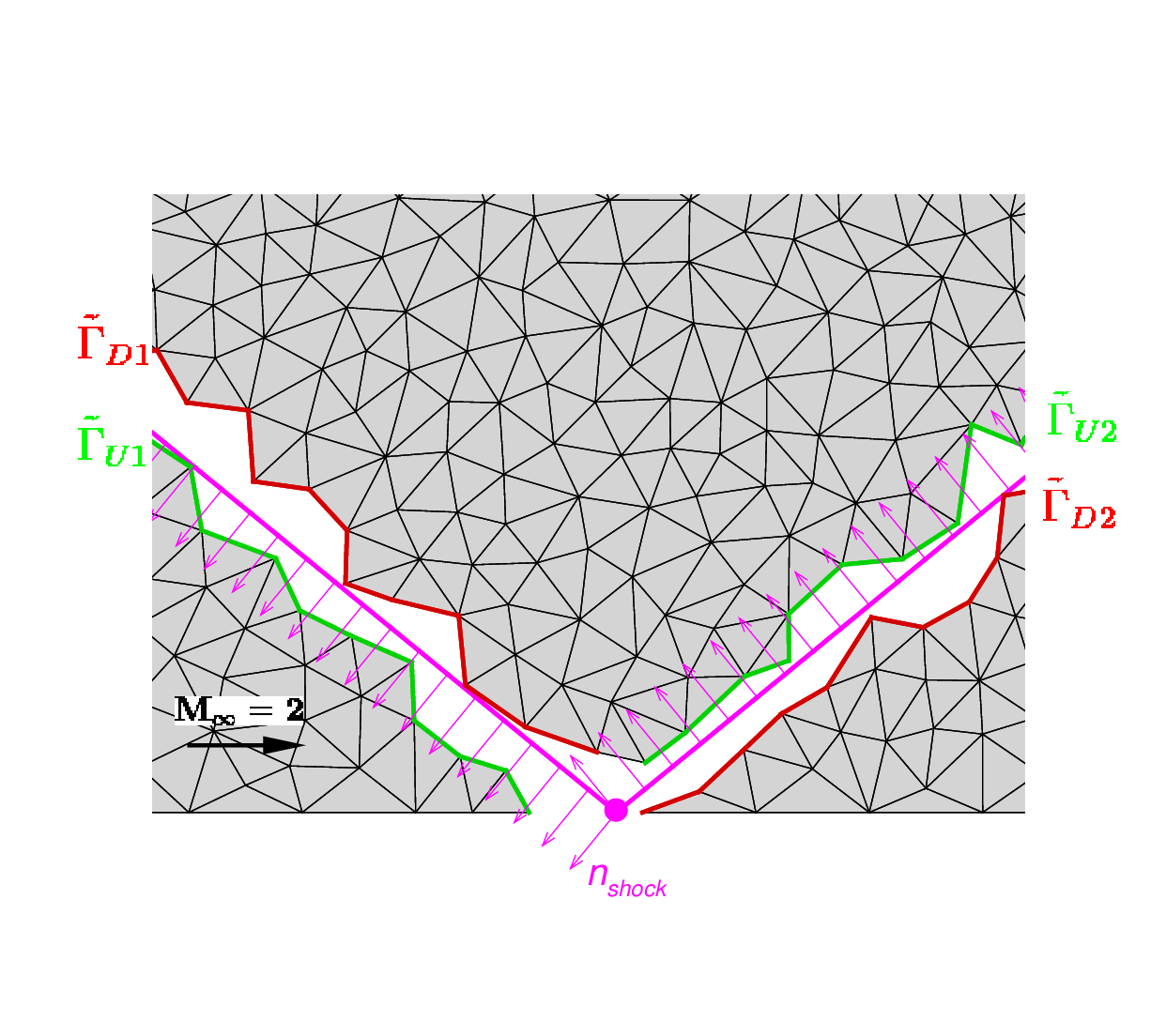}
\caption{Definition of surrogate boundaries when multiple fronts interact;
	\revAB{in the regular reflection shown here, the surrogate boundaries marked in green are located
	upstream w.r.t.\ the incident and reflected shocks, whereas those marked in red are located downstream}.}\label{inter_alg}
\end{figure}

%
\revMR{In principle the flow may evolve 
$n$ different  
fronts, which need to be handled numerically.  The identification of these discontinuities is in itself a challenging problem (see e.g. \cite{PA:20,B:20} and references therein), which is out of the scopes of this work.
Here we assume to be given in advance a set  $\{\mathcal{F}_j\}_{j=1}^n$ of $n$ discontinuities, as well as their topology and nature (shocks and/or contact discontinuities), 
and the set of numbered sub-domains $\{\Omega_l\}_{l=1}^m$ separated by the discontinuities. We also start from a brute set of the  ensemble of all \revAB{points} and edges of \revAB{the} surrogate discontinuities, which we denote by $\tilde\Gamma$.  We process this ensemble as follows: }
\revMR{
\begin{enumerate}
	\item associate to each point in $\tilde\Gamma$ the \revAB{index}
		of the closest discontinuity (distance measured by orthogonal projection);
	\item associate to each point in $\tilde\Gamma$ the {\textit U}pstream/{\textit D}ownstream flag based on its position w.r.t.\ the orientation of the front normals;
	\item for each discontinuity assemble the corresponding arrays of upstream and downstream surrogate boundaries (edge collection) as better shown in Fig.~\ref{inter_alg};
	\item build a pointer providing the explicit mapping between actual and surrogate discontinuities.
\end{enumerate}
}
\revMR{The  connectivity obtained allows to move easily from one surrogate discontinuity to the corresponding 
front-mesh, or to the other surrogate corresponding to the same discontinuity 
(e.g.\ move \revAB{between the} upstream and downstream surrogates $\tilde\Gamma_{Uj}$
and $\tilde\Gamma_{Dj}$ of the two shocks of Fig.~\ref{inter_alg}),
as well as between different surrogates bounding the same sub-domain (e.g.\ from $\tilde\Gamma_{D1}$ to $\tilde\Gamma_{U2}$ \revAB{in Fig.~\ref{inter_alg}}).}
%

\subsection{Solution update using the CFD solver}\label{cfd solver}
\begin{figure}
\begin{center}
\includegraphics[width=0.35\textwidth]{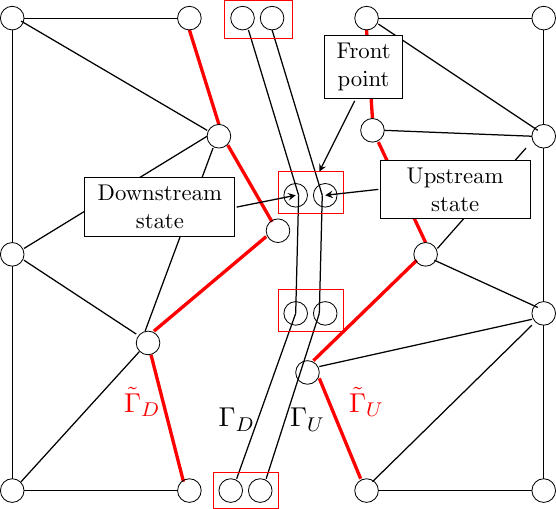}
\end{center}
\caption{The solution update is performed using the computational mesh.}\label{pic4}
\end{figure}

The solution is updated to time level \emph{t+$\Delta$t} using  a 
shock-capturing code. 
 The flow computations are performed on the 
 non-communicating \revAB{sub-}domains separated by the \revMR{front cavity,
 and including the surrogate discontinuities $\tilde{\Gamma}_U$ and $\tilde{\Gamma}_D$ as boundaries}
 (see Fig.~\ref{pic4}). 
 \revMR{As we will see in the next section, the boundary conditions imposed on the surrogate discontinuities are defined
 starting from the values of the flow variables on $\Gamma_U$ and $\Gamma_D$ which, as already said, are
in general different}, \revAB{because $\Gamma$ is a discontinuity.}\\
%
Concerning the solver used in this paper, it is based on a \emph{Residual Distribution  (RD)} method evolving in time 
approximations of the values of the flow variables in grid-points. The method has several appealing characteristics, 
including the possibility of defining genuine multidimensional upwind strategies for Euler flows, by means of a wave decoupling 
exploiting  appropriately preconditioned forms of the equations~\cite{Bonfiglioli1}. 
By combining ideas from both the stabilized finite element and finite volume methods, these schemes allow to achieve 
second order of accuracy and monotonicity preservation with a compact stencil of nearest neighbors.
The interested reader can refer to~\cite{dr2017,ar2017} and references therein for an in-depth review of this family of methods, 
as well as to~\cite{Bonfiglioli1,Bonfiglioli2} and references therein for some specific choices of the implementation used here.
%
%
\subsection{\revMR{Surrogate discontinuity  conditions update}} \label{transfer}
The flow solver provides updated nodal values at all grid-points of the computational mesh 
at time level \emph{t+$\Delta$t}. 
When dealing with \textit{shock waves}, the shock-upstream surrogate boundary, $\tilde{\Gamma}_U$ behaves like 
a supersonic outflow and, therefore, no boundary conditions should 
be applied. \revMR{However, along the}
shock-downstream surrogate
$\tilde{\Gamma}_D$ 
the 
flow is subsonic \revAB{in the shock-normal direction} \revMR{and in principle only the}
\revAB{characteristic variable\footnote{\revAB{Hereafter, characteristic variables shall also be referred to as Riemann variables}.} conveyed by the} \revMR{slow acoustic wave has been correctly updated, 
while boundary conditions for the remaining ones}
(the fast acoustic, entropy and vorticity waves) \revMR{should be imposed. The situation is similar on both sides of a \textit{contact discontinuity}}.\\


\revMR{Updating correctly the interface conditions across the embedded discontinuity is very delicate, and it is the key of the method proposed. It involves several steps, which are detailed hereafter.}
\subsubsection{Computational mesh to discontinuity mesh transfer}\label{transfer1}
\begin{figure}[htb]
\begin{center}
\subfloat[From the surrogate boundary $\tilde{\Gamma}_U$ to the shock-upstream side of the front $\Gamma_U$.]{\includegraphics[width=0.40\textwidth]{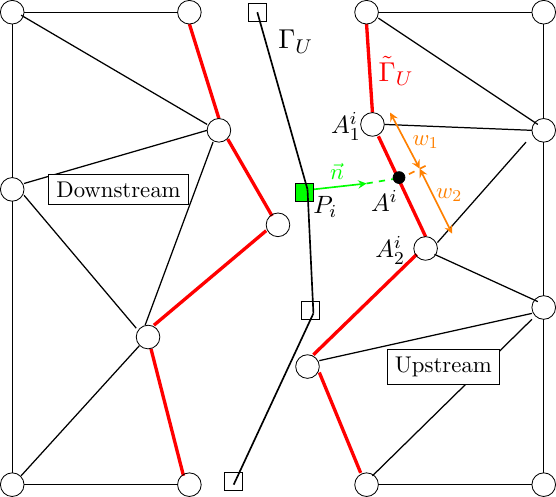}}\qquad 
\subfloat[From the surrogate boundary $\tilde{\Gamma}_D$ to the shock-downstream side of the front $\Gamma_D$.]{\includegraphics[width=0.40\textwidth]{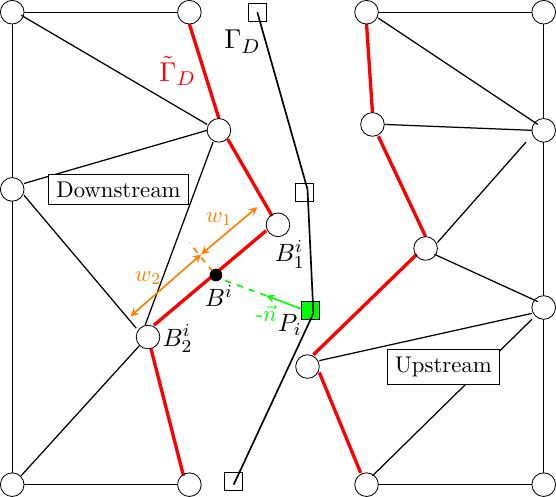}}
\end{center}
	\caption{First transfer \revAB{between the the surrogate boundaries and the front-mesh}.}\label{pic8}
\end{figure}

\revMR{The transfer of the flow data from the computational mesh to the front-mesh is necessary before imposing the interface conditions.
In this work we use a CFD solver strongly relying on the use of Roe's parameter vector $\mathbf{Z}=\sqrt{\rho}\left(1,H,u,v\right)^t$~\cite{Roe1,deconinck93} as the dependent variable,
so the transfer is also done in terms of $\mathbf{Z}$. This is evidently not a necessary choice, and any other choice of state vector is   acceptable.}
\revMR{This} first transfer is required to update the 
\revMR{flow variables along} 
$\Gamma_U$ and $\Gamma_D$. \revMR{Following \cite{Scovazzi3,EST2020}, this is achieved by first defining a map:
\begin{equation}\label{sbm1map}
\tilde{\mathbf{x}}=M(\mathbf{x})
\end{equation}
where, as shown in Fig.~\ref{pic8}, the point
$\tilde{\mathbf{x}} = \revAB{\mathbf{x}(A^i)}$ belonging to the 
computational mesh is the projection in the direction  of the front-normal  $\mathbf{n}$ of a \revAB{front-}point 
$\mathbf{x} = \revAB{\mathbf x(P_i)}$ of the front-mesh. 
Then we use  a forward Taylor series expansion to express the data along the discontinuity
in terms of the values of the flow variables and of their derivatives  on the computational mesh. 
A first-order transfer reads:}
%
\begin{equation}\label{sbm1order1}
\mathbf{Z}(\mathbf{x})\,=\,\mathbf{Z}(\tilde{\mathbf{x}})\,+\,\mathcal{O}(\|\mathbf{x}\,-\,\tilde{\mathbf{x}}\|)\,,
\end{equation}
\revMR{while a second order extrapolation is written as:} 
\begin{equation}\label{eq5}
\mathbf{Z}(\mathbf{x})\,=\,\mathbf{Z}(\tilde{\mathbf{x}})\,+\,\nabla  \mathbf{Z}(\tilde{\mathbf{x}})\cdot(\mathbf{x}\,-\,\tilde{\mathbf{x}})\,+\,\mathcal{O}(\|\mathbf{x}\,-\,\tilde{\mathbf{x}}\|^2).
\end{equation}
\revMR{With the exception of few methods based on a $C^1$ approximation (see e.g.~\cite{giorgiani}),
in general, the gradients of the flow variables are undefined at mesh-\revAB{points}. 
To avoid handling specific singular cases, we have
implemented~\eqref{eq5} using interpolated reconstructed nodal gradients.} 
%
%
%
Note that, in order to achieve an overall second order of accuracy in the calculation
of $\mathbf{Z}(\mathbf{x})$, the approximation of the gradient in Eq.~\eqref{eq5} only needs to be consistent, i.e.\ \revMR{at least} first-order-accurate.\\

%
%
\revMR{In practice we proceed as shown in Fig.~\ref{pic8}. 
When moving from the computational mesh to $\Gamma_U$ (see Fig.~\ref{pic8}a)} 
\revMR{a \revAB{front-}point $P_i$ is mapped to a} point $A^i$.
The value\revMR{s} of the dependent variables and \revMR{of} their gradients in 
 $A^i$ are \revMR{interpolated from the neighbors} 
$A_1^i$ and $A_2^i$: 
\begin{equation}\label{eq6}
\phi(A^i)  = w_2\,\phi(A_1^i) \,+ w_1\,\phi(A_2^i)
\end{equation}
where $\phi$ is either $\mathbf{Z}$ or $\nabla \mathbf{Z}$, and $w_1$ and $w_2$ are the weights,
equal to the normalized distances between $A^i$ and grid-points $A_1^i$ and $A_2^i$.
The evaluation of the gradient in the grid-points of the surrogate boundaries
may be performed using different approaches, as reported in~\S~\ref{gradient reconstruction}.
Once the value of $\mathbf{Z}$ \revAB{and $\nabla \mathbf{Z}$} in 
point $A^i$ has 
been computed using Eq.~(\ref{eq6}), 
$\mathbf{Z}$ in $P_i$ is computed by means of Eq.~\eqref{eq5}, 
having set $\mathbf{x} = \mathbf{x}(P_i)$ and $\tilde{\mathbf{x}} = \mathbf{x}(A^i)$.\\
%
%
\revMR{When moving from the computational mesh to \revAB{$\Gamma_D$} the}
procedure is identical:  \revMR{the \revAB{front-}point $P_i$ is} 
updated 
\revMR{using its mapped point} $B^i$ (see Fig.~\ref{pic8}b). \revMR{Values of the relevant quantities are again linearly interpolated, and} 
%
$\mathbf{Z}$ 
is computed by means of Eq.~\eqref{eq5},  having set $\mathbf{x} = \mathbf{x}(P_i)$ and $\tilde{\mathbf{x}} = \mathbf{x}(B^i)$.
%
\subsubsection{\revMR{Riemann \revAB{variables} and jump conditions enforcement}}\label{shocks calc}
\revMR{The solution updated by the flow solver and extrapolated to the discontinuity mesh does not take into account the coupling between the different sub-domains. 
\revRP{This coupling is done here based on the  Rankine-Hugoniot jump relations.}
On the upstream side $\Gamma_U$ \revAB{of a \textit{shock wave}} all the \revAB{Riemann variables} are transported into the shock, so the data extrapolated  from the \revAB{computational} mesh on this side is assumed to be correct.
On the shock-downstream side, however, 
this \revAB{holds true} only  for the Riemann 
variable associated with the slow acoustic wave that moves towards the shock (the subscript $_D$ denotes  values along $\Gamma_D$): 
\begin{equation}\label{eq4}
 R_D = 
  \;a^{t+\Delta t}_{D}\,+\,\frac{\gamma-1}{2}\,(u_n)^{t+\Delta t}_{D}
\end{equation}
We need to provide relations to compute along $\Gamma_D$ the  values of the Riemann \revAB{variables} corresponding to the  fast acoustic wave,
as well as to the entropy and   vorticity waves. These relations are \revAB{supplied} by the Ranking-Hugoniot \revAB{jump} conditions which, for a perfect gas, can be recast as:
\begin{equation}\label{RH}
        \begin{array}{ccl}
        \rho^{t+\Delta t}_D (u_n)^{t+\Delta t}_D -\,\omega_s^{t+\Delta t} \rho^{t+\Delta t}_D&=&         \rho^{t+\Delta t}_U(u_n)^{t+\Delta t}_U -\,\omega_s^{t+\Delta t} \rho^{t+\Delta t}_U  \\[5pt]
        \rho_D^{t+\Delta t}((u_n)_D^{t+\Delta t}-\,\omega_s^{t+\Delta t})^2+p_D^{t+\Delta t}&=&\rho_U^{t+\Delta t}((u_n)_U^{t+\Delta t}-\omega_s^{t+\Delta t})^2+p_U^{t+\Delta t}\\[5pt]
        \dfrac{\gamma}{\gamma-1}\dfrac{p_D^{t+\Delta t}}{\rho_D^{t+\Delta t}} + \dfrac{1}{2}((u_n)_D^{t+\Delta t}-\omega_s^{t+\Delta t})^2&=&\dfrac{\gamma}{\gamma-1}\dfrac{p_U^{t+\Delta t}}{\rho_U^{t+\Delta t}}+ \dfrac{1}{2}((u_n)_U^{t+\Delta t}-\omega_s^{t+\Delta t})^2 \\[15pt]
        (u_{\tau})_D^{t+\Delta t}&=&(u_{\tau})_U^{t+\Delta t}
        \end{array}
\end{equation}
with $\omega_s^{t+\Delta t}$ the shock speed. \revRP{The algebraic non-linear system made up of 
Eqs.~\eqref{eq4} and~\eqref{RH}} allows to compute the five unknowns $\left(\rho_D,\, \mathbf{u}_D,\, p_D,\, 
\omega_s\right)$ given
the known upstream state $\left(\rho_U,\, \mathbf{u}_U,\, p_U\right)$ 
and the \revAB{slow acoustic} Riemann \revAB{variable} $R_D$.}\\

%
\revMR{For a \textit{contact discontinuity}, we use as input the 
following set of \revAB{Riemann variables}, which we assume to be correctly updated in the computational \revAB{mesh}:} 
\begin{equation}\label{cd1}
\begin{split}
R_D  = 
 \;a^{t+\Delta t}_D \,+\,\frac{\gamma-1}{2}\,(u_n)^{t+\Delta t}_D\;, &\;\;  
  R_U  = 
	\;a^{t+\Delta t}_U   \,\revRP{-}\,\frac{\gamma-1}{2}\,(u_n)^{t+\Delta t}_U   \\[5pt]
 s_D  = 
\; \dfrac{p^{t+\Delta t}_D}{(\rho^{t+\Delta t}_D)^{\gamma}} \;,& \;\; 
s_U  = 
 \; \dfrac{p_U^{t+\Delta t}}{(\rho_U^{t+\Delta t})^{\gamma}}\\[10pt]
V_D  = 
\;(u_{\tau})_D^{t+\Delta t} \;,&\;\; 
 V_U  = 
\;(u_{\tau})^{t+\Delta t}_U 
\end{split}
\end{equation}
\revMR{Note that \revAB{when dealing with a contact discontinuity} the situation is essentially symmetric and
the jump \revAB{relations} simplify somewhat, and can be shown to reduce to:}
\begin{equation}\label{CD}
        \begin{array}{ccl}
         p_U^{t+\Delta t}\,&=&\,p_D^{t+\Delta t}  \\[5pt]
        (u_n)_U^{t+\Delta t} \,&=&\, \omega_{cd}^{t+\Delta t}  \\[5pt]
        (u_n)_D^{t+\Delta t} \,&=&\, \omega_{cd}^{t+\Delta t}  
        \end{array}
\end{equation}
\revAB{The algebraic non-linear system made up of Eqs.~(\ref{cd1}) and~(\ref{CD}) allows to compute the nine unknowns
 $(\rho_U,\, \mathbf{u}_U,\, p_U,\, \rho_D,\, \mathbf{u}_D,\, p_D,\, \omega_{cd})$, given the six
Riemann variables $(R, s, V)$
on the two sides of the discontinuity and the jump relations~(\ref{CD})}.\\

%
\revAB{For both shocks and contact discontinuities, a non-linear system of algebraic equations} 
\revMR{needs to be solved in each discontinuity-point. This is done 
using the Newton-Raphson root-finding algorithm,
thus providing the correct states and $\omega$ at time level \emph{t+$\Delta$t}}.\\
%
%
Whenever two or more discontinuities interact \revMR{special care is required for} the front-points belonging to \revMR{different discontinuities.}
The numerical treatment of these \revMR{points is done in a case by case manner and detailed in the sections related to the specific tests } 
 in~\S~\ref{results}. 
%
\subsubsection{Discontinuity mesh to computational mesh  transfer}\label{transfer2}
Once the \revMR{value of the unknowns have been corrected along the discontinuity mesh to account for the jump conditions, the corrected values need 
to be \revAB{transferred back to} the surrogate discontinuities. To do this, we proceed in a similar manner as before} \revAB{by writing:}
%
\begin{equation}\label{secondtransfer}
\mathbf{Z}(\tilde{\mathbf{x}})\,=\,\mathbf{Z}(\mathbf{x})\,-\,\nabla  \mathbf{Z}(\tilde{\mathbf{x}})\cdot(\mathbf{x}\,-\,\tilde{\mathbf{x}})\,+\,\mathcal{O}(\|\mathbf{x}\,-\,\tilde{\mathbf{x}}\|^2).
\end{equation} 
\revMR{Note that there is a substantial difference in the direction in which the Taylor expansion is performed which is now a backward expansion from the arrival point (the surrogate discontinuity) to the point where the data is \revAB{available} (the discontinuity mesh). As in the previous case, several choices are
possible to evaluate the gradient involved in this transfer.} 
We shall focus on \revMR{this} 
in~\S~\ref{gradient reconstruction}. 
\revAB{In} the first-order accurate case, Eq.~\eqref{secondtransfer} \revAB{simplifies} to:
\begin{equation}\label{sbm1order2}
\mathbf{Z}(\tilde{\mathbf{x}})\,=\,\mathbf{Z}(\mathbf{x}) \,+\,\mathcal{O}(\|\mathbf{x}\,-\,\tilde{\mathbf{x}}\|)
\end{equation}
\revMR{Finally note that for} 
shock-waves 
only 
grid-points on $\tilde{\Gamma}_D$ need to be updated. 
%
\subsection{Front displacement and nodal re-initialization}\label{shock displacement}
The new position of the front at time level \emph{t+$\Delta$t} is computed by displacing all discontinuity-points 
using the following first-order-accurate (in time) formula:
\begin{equation}\label{displ}
\mathbf{x}\left(P^{t+\Delta t}\right)\,=\,\mathbf{x}\left(P^t\right)\,+\,
\omega^{t+\Delta t}\,\mathbf{n}\,\Delta t
\end{equation}
where $\mathbf{x}\left(P\right)$ denotes the geometrical location of the shock-points, and 
$ \omega^{t+\Delta t}$  may represent either $\omega_s^{t+\Delta t}$ or $\omega_{cd}^{t+\Delta t}$.
The use of a first-order-accurate temporal integration formula \revMR{has no impact} 
as long as steady flows are of interest. 
For unsteady flows, second-order-accurate time integration formulae should be used, as done for example in~\cite{Paciorri4,Campoli2017}.\\

\begin{figure}[t]
\begin{center}
\includegraphics[width=0.6\textwidth]{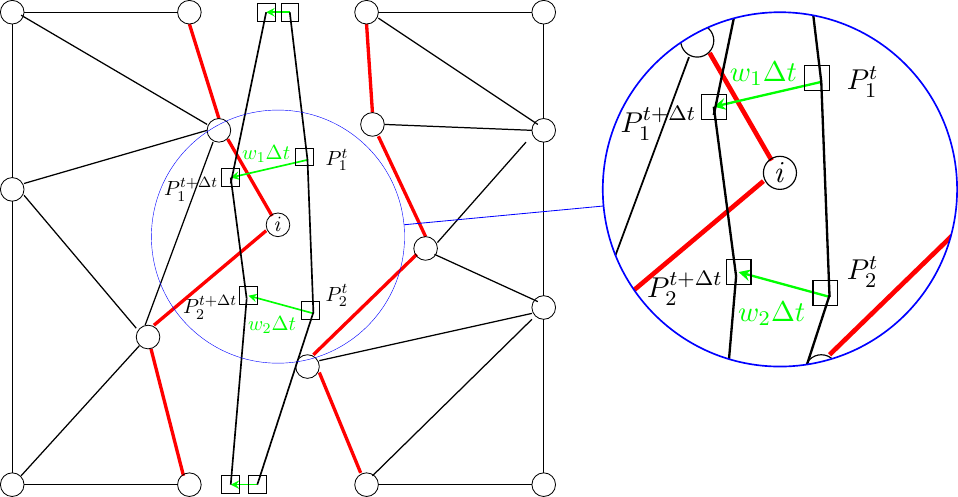}
\end{center}
\caption{The front overtakes a grid-point of the background mesh during its motion.}\label{pic6}
\end{figure}

Since the discontinuity can freely float over the background triangulation, \revAB{it may happen that it
crosses} the surrogate boundaries. \revAB{This situation has been sketched in Fig.~\ref{pic6},} 
where grid-point $i$ has been overtaken by the moving front. 
Whenever this happens, the flow state within grid-point $i$ has to be changed accordingly.  
%
%
\revMR{In particular, the flow variables in these points are} 
re-computed through an interpolation. 
To evaluate whether a grid-point $i$ has been overtaken or not by the discontinuity
a simple approach has been implemented. 
Since $\mathbf{n}$ points always upstream, the sign of the following scalar product allows us to 
distinguish the grid-points located upstream from those located downstream.
\begin{equation}  
\alpha_i\,=\,(\mathbf{x}_p - \mathbf{x}_i)\cdot \mathbf{n}\,=\,
	\begin{cases} <  0 \;,i\; \text{is upstream}\\ > 0 \;,i\; \text{is downstream} \end{cases}
\end{equation}
where $\mathbf{x}_i$ and $\mathbf{x}_p$ are, respectively, the coordinates of a grid-point $i$ and its projection over the closest  
front-edge and $\mathbf{n}$ is the normal vector computed in the closest front-point to $\mathbf{x}_i$. 
If $\alpha_i$ changes sign when the front moves, this means that grid-point $i$ has been overtaken and its state has to be updated.
%
\subsection{\revMR{Evaluation of the nodal gradients}}\label{gradient reconstruction}
\revMR{The technique used  to extrapolate the flow variables back and forth between the surrogate discontinuities and the discontinuity mesh involves 
the knowledge of nodal gradients which need to be reconstructed. The solvers used in this work are based on collocated nodal methods, so we will discuss
the recovery of the nodal gradients in  this specific case. The generalization to other cases is of course possible,
but left out of this work.\\

As for the enforcement of the jump conditions  we distinguish two main cases. For a shock wave, in the upstream domain we  already said that all characteristic information runs into the shock.
For this situation, it seems natural to evaluate the nodal gradients along the surrogate discontinuity based on the data pre-computed by the flow solver within the upstream domain.
To this end we have tested two gradient recovery strategies: }
%
\begin{itemize}
	\item  a Green-Gauss (GG) gradient reconstruction, \revMR{essentially boiling down to an area weighted average of the gradients in the cells surrounding a \revAB{grid-point};} 
\item an \revMR{area-}weighted version 
of Zienkiewicz-Zhu's \revMR{patch super-convergent method} (ZZ)~\cite{ZZ}.
\end{itemize}
\revMR{Explicit formulas for both methods are provided in~\ref{appendixA}}.

\begin{figure}
\begin{center}
\includegraphics[align=c,width=0.45\textwidth]{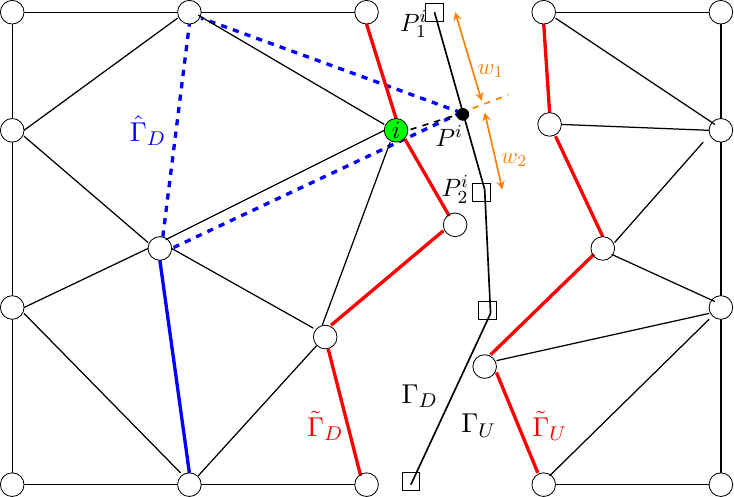}\includegraphics[align=c,width=0.450\textwidth]{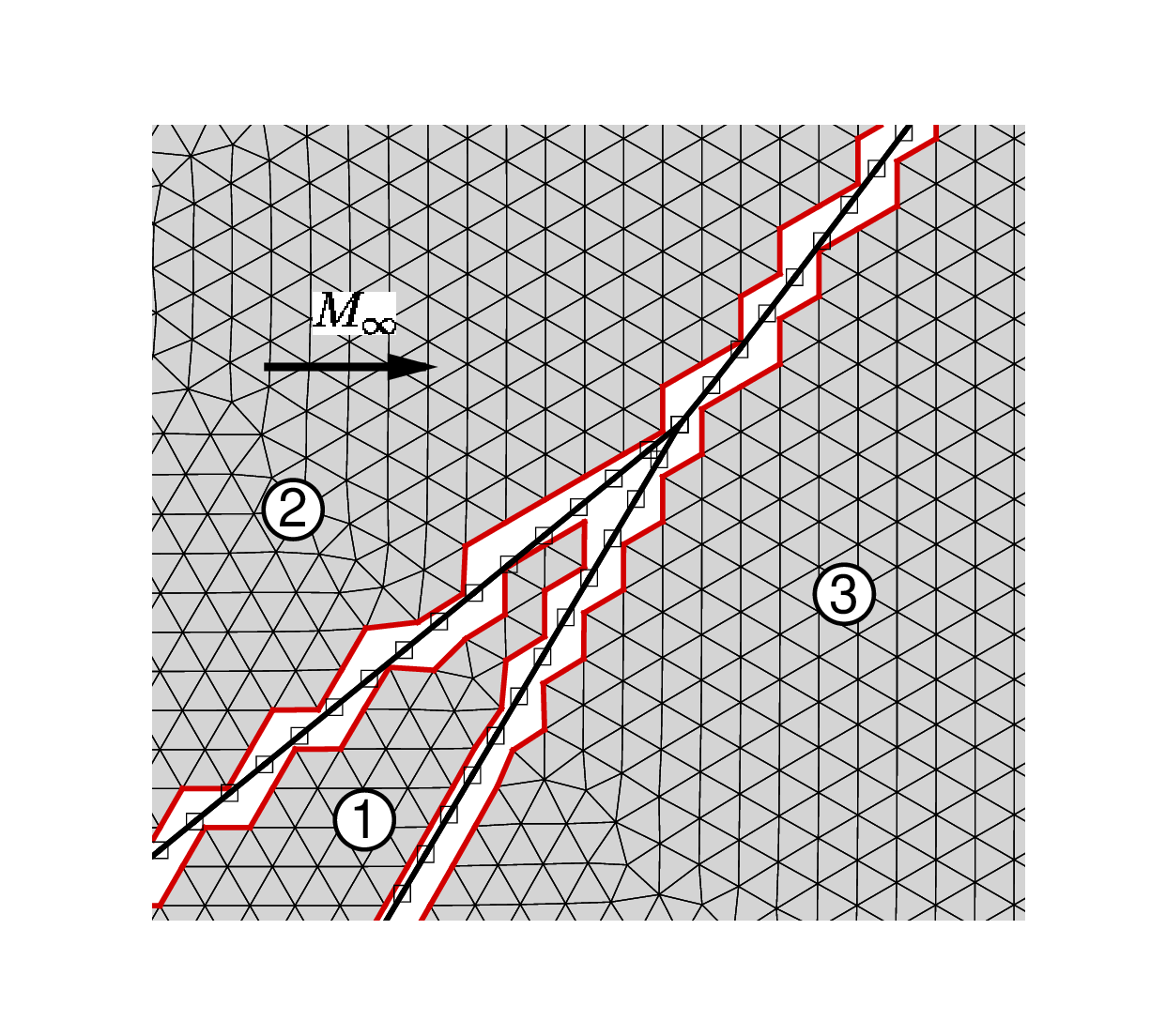}
\end{center}
\caption{Left: auxiliary triangle defining the nodal gradient with corrected data on the front. Right: a problematic example of an interaction for which no auxiliary triangle can be built close to the interaction point in domain 1.}\label{SSinterProblem}
\end{figure}

\revAB{Within the region downstream of a shock,} \revMR{as well as for contact discontinuities, 
as discussed in the previous sections, the data computed by the flow solver must be corrected to account 
for the jump conditions.
This fact has led us in the past to include the corrected data in the evaluation of the nodal gradients in these cases~\cite{EST2020}.}
%
%
%
In \revMR{the \revAB{aforementioned} reference, which did not account for discontinuity interactions,} 
 the gradient term of Eq.~\eqref{secondtransfer} \revMR{is replaced} with the cell-wise gradient of an
auxiliary triangle \revAB{(marked using a dashed line in Fig.~\ref{SSinterProblem})} \revMR{containing the surrogate 
\revAB{grid-point} to be updated (\revAB{grid-point $i$ in Fig.~\ref{SSinterProblem}}), and }
\revMR{defined} by a point on the discontinuity\revMR{-mesh} (\revAB{point $P^i$
in Fig.~\ref{SSinterProblem}}), and two \revAB{grid-}\revMR{points \revAB{of} the computational mesh, 
not belonging to the surrogate discontinuity}. 
However,  when \revMR{interactions are present, singular topological situations may arise in which this approach cannot be used. For example, for}
 region~$1$ in Fig.~\ref{SSinterProblem}, \revMR{a long strip of one row of elements is trapped within  two discontinuities and the interaction point.  
For the \revAB{grid-points} along the boundary of this strip there are no neighboring inner \revAB{grid-points} 
to construct the auxiliary triangle. Several \revAB{other gradient-reconstruction} formulations have \revAB{also} been tested. 
}
The simplest involves using a one-sided recovered gradient not accounting for the Rankine-Hugoniot correction. 
The latter
%
not only turned out to be  very useful for computing interacting discontinuities, but \revMR{much simpler, especially in view of possible extensions to three space dimensions.}
%
\revMR{In the following, we will refer to 
the method of~\cite{EST2020} as $eST$, while $eDIT_{GG}$ (GG reconstruction), and   $eDIT_{ZZ}$ (ZZ reconstruction) will denote the extrapolated discontinuity tracking obtained using one sided reconstructions.}
Finally, $eDIT_{FO}$ (first order) \revAB{refers to simulations run with a second-order-accurate
discretization of the governing PDEs,
but only first-order accurate
data transfer between the surrogate and actual discontinuities, i.e.\ by
using Eqs.~(\ref{sbm1order1}) and~(\ref{sbm1order2})}.

\revMRR{

\section{Remarks on conservation}\label{remarksconservation}
\begin{figure}
\begin{center}
\includegraphics[width=0.40\textwidth]{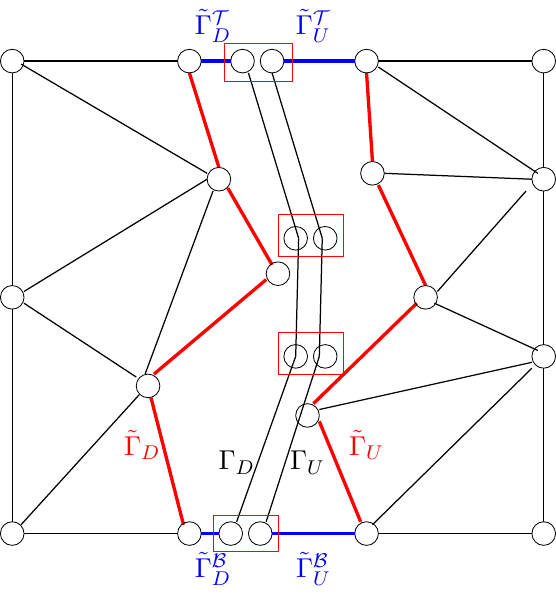}
\end{center}
\caption{Area of the cavity $\Omega_e$ when open shock geometries are considered.}\label{picOpenShock}
\end{figure}

We make a small detour to discuss the issue of conservation.  
A first important remark is that  the notion of conservation  is essential when considering the approximation of  shocks.
Failing to properly account  for conservation of mass, momentum, and energy would lead to a wrong approximation of 
these features, both in terms of position and strength. In smooth regions, as well as across contact discontinuities,
the use of non-conservative approaches is less critical, and in some cases even advantageous  (e.g.\ \cite{ABGRALL1996150,crd02,dr2017,ABGRALL201810}).

It is also important to realize that the definition of discrete conservation is associated to the identification of both
a   geometrical  cell over which conservation is expressed, and a numerical flux expressing local conservation.
The interested reader can refer to \cite{ar2017,goncalves22} for a discussion. 
Typically, global conservation over the domain is thus expressed as
\begin{equation}\label{global_cons} 
\int\limits_{\partial\Omega} \mathbf{\hat F}_n(\mathbf{Z}) \text{d}\Gamma =0
\end{equation}
with $\mathbf{\hat F}_n$ the numerical flux associated to the boundary conditions. 
Now note that, even for exact quadrature,  the numerical flux is not equal to the physical one $\mathbf{F}_n(\mathbf{Z})$, 
and the difference between the two is in general of the order of  the truncation error of the discretization. So conservation is still
verified within some numerical approximation.\\

The situation is similar for the embedded approach proposed here, for which, despite the exact imposition of the jump conditions across 
discontinuities,  conservation can  only be  measured within the truncation of the extrapolation method.  To be more precise,  
let us focus on the   configuration of Fig.~\ref{picOpenShock}. 
In the figure, the domains on the right of $\tilde\Gamma_U$ is the upstream domain, while the domain on the left
of  $\tilde\Gamma_D$ is  the downstream one. The flow is assumed to go from the right to the left. A discontinuity is placed in the middle
of the domain. We have denoted by $\Gamma_U$ and $\Gamma_D$ the upstream and downstream  sides of the  discontinuity.
The discretization of all the domains,
 including the upstream/downstream boundaries and the discontinuity,  is shown in the figure. In particular, as already discussed,
unlike in the unstructured shock-fitting method~\cite{Paciorri1}, for the $eDIT$ method the shock-edges are no longer part of the computational mesh.
In the region between the shock mesh and the computational mesh, the conservation law is then replaced by the
extrapolation procedure, which is the main source of loss of conservation. Note however, that this loss only occurs in the smooth parts of the flow. 
In particular,  denoting by $\Omega_e^U$ the region enclosed between $\Gamma^U$ and $\tilde \Gamma^U$,
given a smooth exact solution of the Euler equations we can write easily the upstream consistency estimate 
 \begin{equation*}
 \begin{split}
\!\!\mathcal{I}_h^U(\mathbf{Z}_h) := &\oint_{\partial\Omega_e^U} \!
\mathbf{F}_n(\mathbf{Z}_h) \text{d}\Gamma =
 \oint_{\partial\Omega_e^U} \!\left[\mathbf{F}_n(\mathbf{Z}_h)-\mathbf{F}_n(\mathbf{Z}_h^{\textsf{ex}})
 \right]\text{d}\Gamma \\ =&
 \oint_{\partial\Omega_e^U} \!\left[
 \mathbf{F}_n\!\left(\mathbf{Z}_h(\tilde{\mathbf{x}}) + \nabla \mathbf{Z}_h(\tilde{\mathbf{x}})\cdot( \mathbf{x} -\tilde{\mathbf{x}}) \right)-
  \mathbf{F}_n\!\left(\mathbf{Z}_h^{\textsf{ex}}(\tilde{\mathbf{x}}) + \nabla \mathbf{Z}_h^{\textsf{ex}}(\tilde{\mathbf{x}})\cdot( \mathbf{x} -\tilde{\mathbf{x}}) \right)
\right]\text{d}\Gamma +\mathcal{O}(\mathbf{d}^2)
\end{split}
\end{equation*}
Where $\mathbf{Z}^{\textsf{ex}}$ is the Roe's parameter vector evaluated for the exact solution.  
Formally replacing the nodal values of the solution with those of the exact one, we obtain the consistency estimate
$\mathcal{I}_h^U(\mathbf{Z}_h^{\textsf{ex}})= \mathcal{O}(\mathbf{d}^2)$. 
Using the exact same arguments, we can write the estimate $\mathcal{I}_h^D(\mathbf{Z}_h^{\textsf{ex}})= \mathcal{O}(\mathbf{d}^2)$
for the loss of conservation in the downstream region  $\Omega_e^D$ enclosed between $\Gamma^D$ and $\tilde \Gamma^D$.

Let us now set $\tilde \Gamma^T=\tilde \Gamma^T_U+\tilde \Gamma^T_D$ and 
 $\tilde \Gamma^B=\tilde \Gamma^B_U+\tilde \Gamma^B_D$.
 By construction of the method proposed here, the jump condition $[\![ \mathbf{F}_n(\mathbf{Z}_h)]\!]_{\Gamma}=0$ is exactly satisfied.  This
allows readily to write  the global consistency estimate on conservation (cf.\ \ref{picOpenShock} for the notation)
\begin{equation}\label{cons}
\int\limits_{\tilde \Gamma_U}\mathbf{F}_n(\mathbf{Z}_h^{\textsf{ex}}) \text{d}\Gamma + \int\limits_{\tilde \Gamma_D}\mathbf{F}_n(\mathbf{Z}_h^{\textsf{ex}}) \text{d}\Gamma  +
\int\limits_{\tilde \Gamma^T}\mathbf{F}_n(\mathbf{Z}_h^{\textsf{ex}}) \text{d}\Gamma 
+
\int\limits_{\tilde \Gamma^B}\mathbf{F}_n(\mathbf{Z}_h^{\textsf{ex}}) \text{d}\Gamma = 
\mathcal{I}_h(\mathbf{Z}_h^{\textsf{ex}})
= \mathcal{O}(\mathbf{d}^2)
\end{equation}
reducing to
\begin{equation}\label{cons1}
\int\limits_{\tilde \Gamma_U}\mathbf{F}_n(\mathbf{Z}_h^{\textsf{ex}}) \text{d}\Gamma + \int\limits_{\tilde \Gamma_D}\mathbf{F}_n(\mathbf{Z}_h^{\textsf{ex}}) \text{d}\Gamma  
=\mathcal{I}_h(\mathbf{Z}_h^{\textsf{ex}})
= \mathcal{O}(\mathbf{d}^2)
\end{equation}
whenever the out/inflow on the top/bottom boundaries is zero. So the conservation error is solely controlled from the accuracy of the extrapolation formula.
This means that first or second order is obtained in terms of  conservation, and thus in terms of position and magnitude of the discontinuities,
depending on whether the gradient correction is included or not. This will be verified in practice in the numerical results section.
Note that this behaviour is  better than what any fully conservative  capturing method can provide, and can of course be further improved by enhancing the extrapolation accuracy.
}

 	\section{Numerical Results}\label{results}
\revMR{To illustrate the capabilities of our method, we provide here examples representative of several}
types of interactions \revMR{which} can occur in \revMR{gas-dynamics}. 
 The solutions obtained \revMR{with the extrapolated tracking method will be compared to {\it i}) full-fledged computations, {\it ii}) hybrid simulations
 in which only some of the discontinuities are explicitly tracked while the others are captured, and {\it iii}) solutions obtained with the standard} 
 shock-capturing method. 
\revMR{To evaluate the gradient recovery strategies, a}
 grid-convergence analysis \revMR{is} 
 presented for two  test-cases (transonic source flow and blunt body problem) to \revMR{asses quantitative differences }
between the approach described in~\cite{EST2020} and the more flexible ones \revMR{proposed here to deal with interactions}.

        \subsection{Transonic source flow}\label{sourceflow}
%
This test-case \revMR{is}
 very useful 
  because of the availability of the analytical solution, which  \revMR{allows} 
\revAB{to} perform grid-convergence studies~\cite{EST2020,Paciorri6,Campobasso2011}. 
Assuming that the analytical velocity field has a purely radial velocity component, it may be easily verified that the
two-dimensional, compressible Euler equations, written in a polar coordinate system, become identical to those governing a
compressible quasi-one-dimensional flow \revMR{with a} 
nozzle area  \revMR{variation linear w.r.t.\ the} 
radial distance from the pole of the reference frame.
\begin{figure*}[h!]
\centering
\subfloat[Sketch of the computational domain.]{%
\includegraphics[width=0.45\textwidth,trim={3.5cm 2.9cm 3cm 2cm},clip]{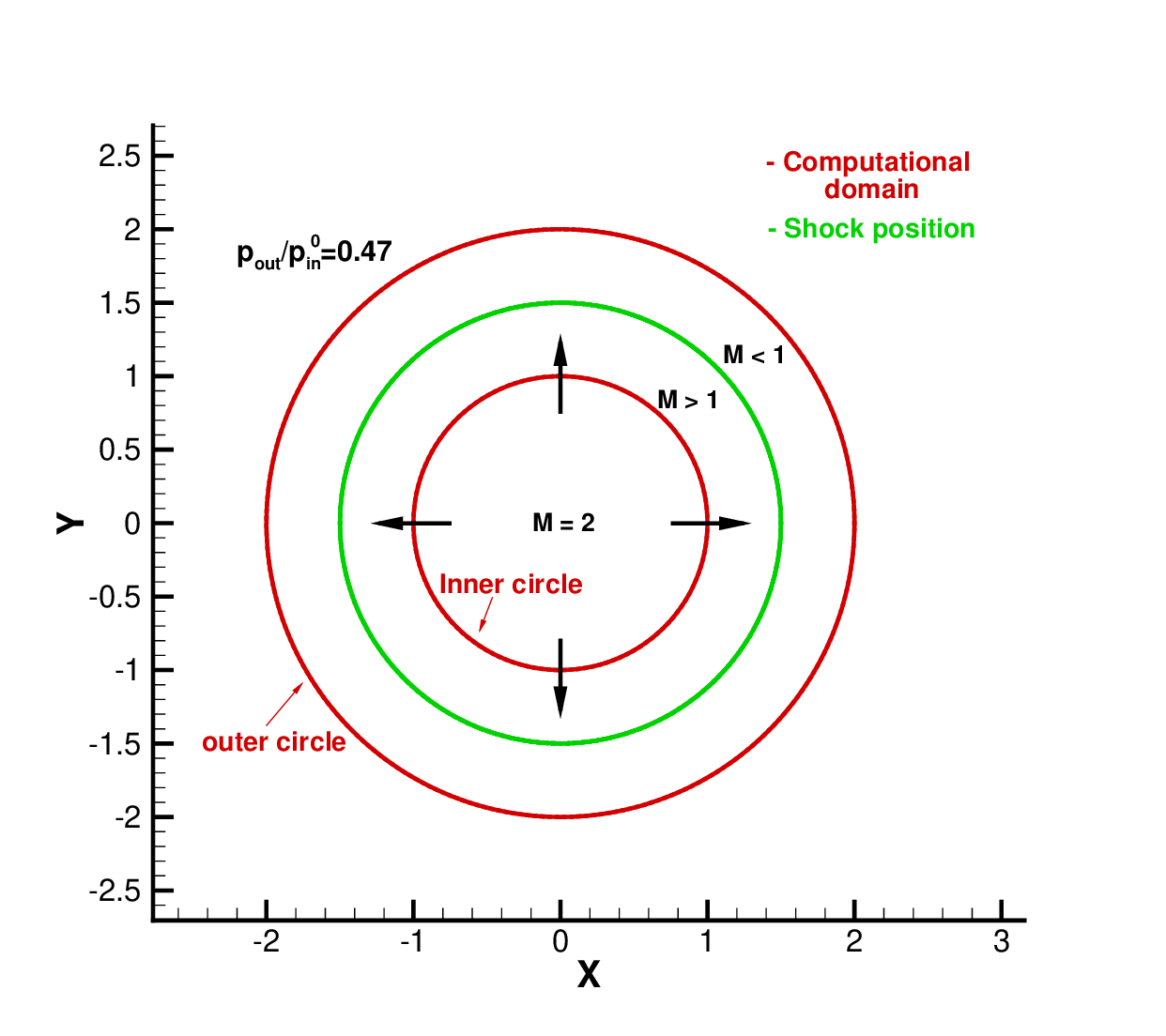}\label{pic20a}}
\subfloat[Detail of the level 0 unstructured mesh inside the first gradient.]{%
\includegraphics[width=0.45\textwidth]{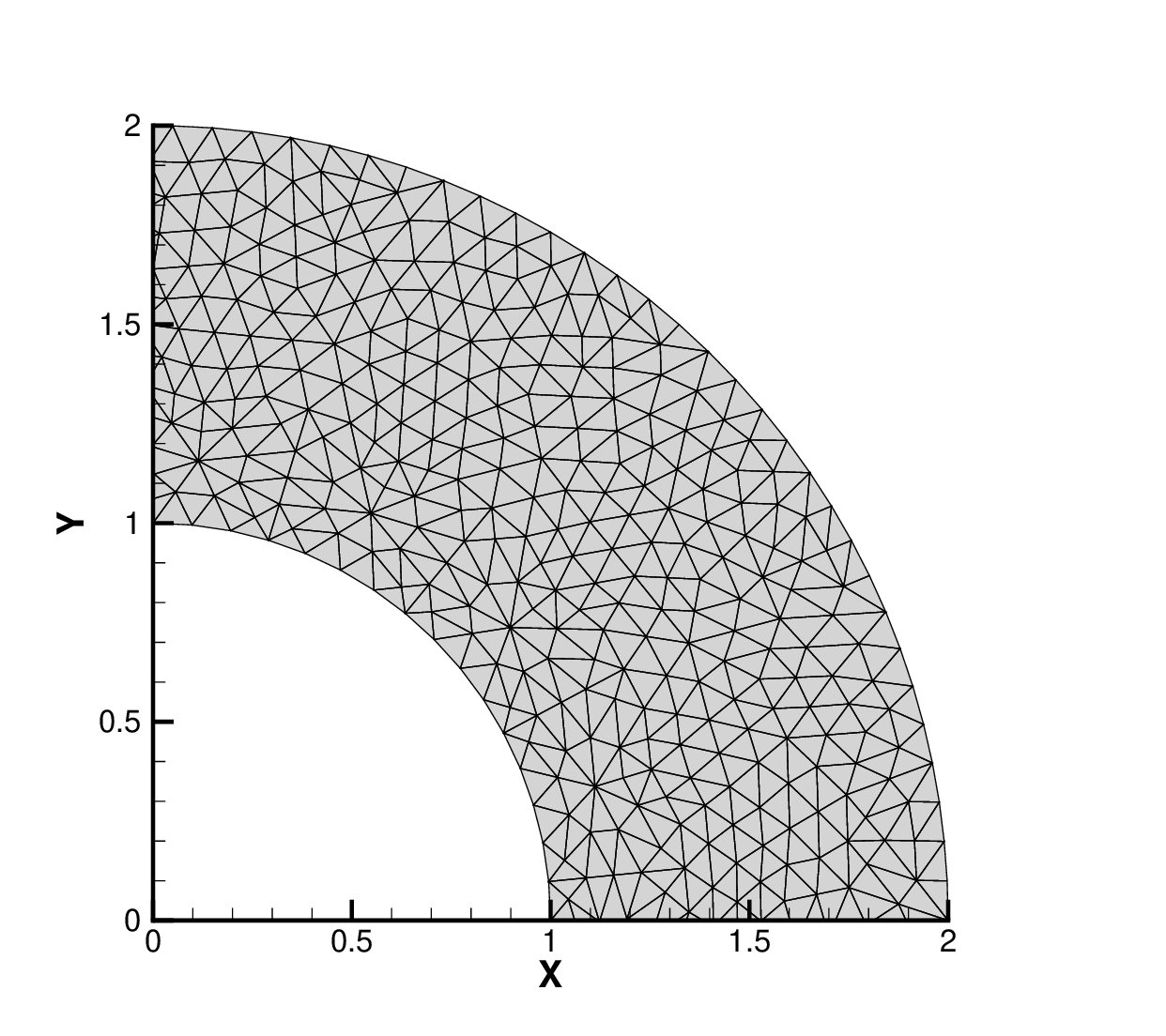}\label{pic20b}}
\caption{Transonic source flow.}\label{pic20}
\end{figure*}
The computational domain consists in the annulus sketched in Fig.~\ref{pic20a}:
the ratio between the radii of the outer and inner circles ($L = r_{in}$) has been set equal to $r_{out}$/$r_{in}$ = 2.
A transonic (shocked) flow has been simulated by imposing a supersonic inlet flow at $M = 2$ on the inner circle and
a ratio between the outlet static and inlet total pressures $p_{out}/p_{in}^0 = 0.47$ such that the shock forms at $r_{sh}$/$r_{in}$ = 1.5.
The Delaunay mesh shown in Fig.~\ref{pic20b}, which contains 6,916 grid-points and 13,456 triangles, has been
generated using the \texttt{\ttfamily gmsh} mesh generator~\cite{geuzaine2009gmsh} in such a way that no 
systematic alignment occurs between the edges of the triangulation and the circular iso-contour lines of the analytical solution.
By doing so, the discrete problem is made truly two-dimensional. The sequence of nested triangulations that have been employed 
for $SC$ and all $eDIT$ computations are summarized in Tab.~\ref{tab1} \revMR{for both the background and the computational grids.} 
\begin{table}[htb]
\scriptsize
\centering
\caption{Transonic source flow: characteristics of the background and computational meshes used to perform the grid-convergence tests.}\label{tab1}
\begin{tabular}{cccccc} \hline\hline
        &\multicolumn{2}{c}{\textit{Background grid}} &\textit{Both grids}  &\multicolumn{2}{c}{\textit{Computational grid}} \\[0.5mm] \cline{2-6}
        Grid level & Grid-points & Triangles & $h$        & Grid-points & Triangles \\[0.5mm] \hline
        0          &  1,369      &  2,548    & 0.5286E-01 &   1,369     &  2,328    \\
        1          &  5,286      &  10,192   & 0.2641E-01 &   5,286     &  9,788    \\
        2          &  20,764     &  40,768   & 0.1320E-01 &   20,764    &  39,968   \\
        3          &  82,296     &  163,072  & 0.6602E-02 &   82,296    &  161,454  \\
\hline\hline
\end{tabular}
\end{table}
 
\revAB{The discretization error for a mesh of size $h$, $\epsilon_h$, is defined
as the difference} between the numerical solution, $u_h$, and the analytical one, $u_0$:
\begin{equation}\label{error1}
\epsilon_h(\mathbf{x}) \,=\,u_h(\mathbf{x})\,-\,u_0(\mathbf{x}) 
\end{equation}
The availability of the exact solution allows to compute the discretization error \textit{locally} using 
Eq.~\eqref{error1}
or, \textit{globally}, by computing the $L_q$ norms of the discretization error over the entire computational 
domain $\Omega$ using Eq.~\eqref{error2}:
\begin{equation}\label{error2}
L_q(\epsilon_h) \,=\,\biggl( \frac{1}{|\Omega|} \int_{\Omega} |\epsilon_h(\mathbf{x})|^q\,d\Omega \biggr)^{1/q} 
,\;\;\;\;\;\;\;\; q\,=\,1,\,2,\,\infty.
\end{equation}
Herein, the $L_1$ norm has been employed and computed using Gaussian quadrature rules.\\
Different versions of the extrapolated shock/discontinuity tracking method will be compared to $SC$ 
computations by doing a thorough grid-convergence analysis.
Results obtained from the aforementioned approaches 
have been separately displayed in the grid-convergence plots for the shock-upstream (supersonic) and shock-downstream (subsonic)
regions, i.e.\ $r < r_{sh}$ and $r > r_{sh}$. In all convergence plots the 
errors are computed \revMR{in terms of the }
parameter vector $\mathbf{Z}$. \\

\begin{figure*}
\centering
\subfloat[$\sqrt{\rho}$]{\includegraphics[width=0.5\textwidth]{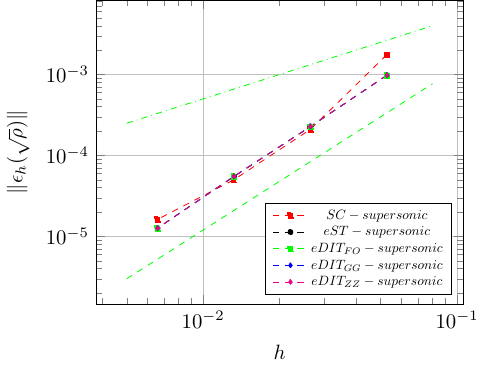}}
\subfloat[$\sqrt{\rho}H$]{\includegraphics[width=0.5\textwidth]{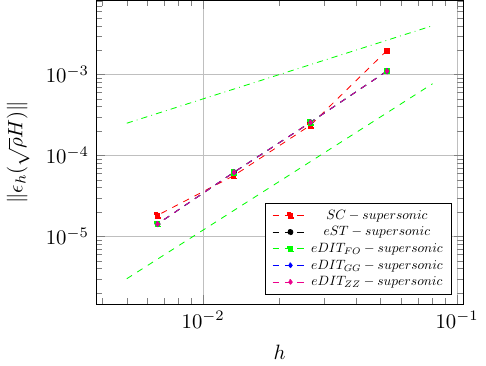}}
\qquad
\subfloat[$\sqrt{\rho}u$]{\includegraphics[width=0.5\textwidth]{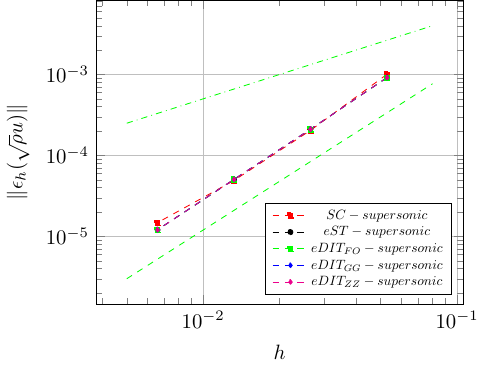}}
\subfloat[$\sqrt{\rho}v$]{\includegraphics[width=0.5\textwidth]{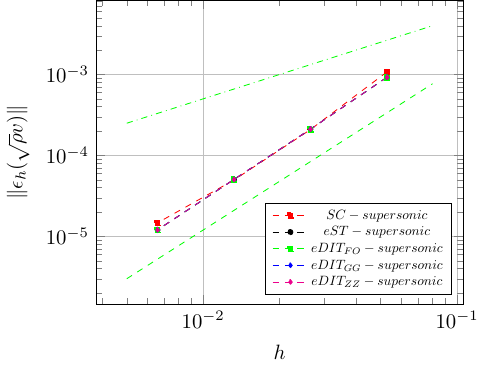}}
	\caption{Transonic source flow: grid-convergence analysis \revAB{within} the shock-upstream (supersonic) region.}\label{Q1D-supersonic}
\end{figure*}

\begin{figure*}
\centering
\subfloat[$\sqrt{\rho}$]{ \includegraphics[width=0.5\textwidth]{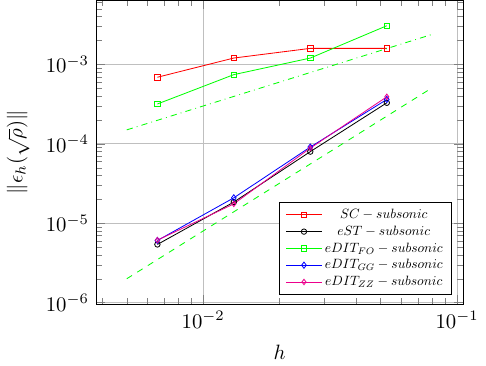}}
\subfloat[$\sqrt{\rho}H$]{\includegraphics[width=0.5\textwidth]{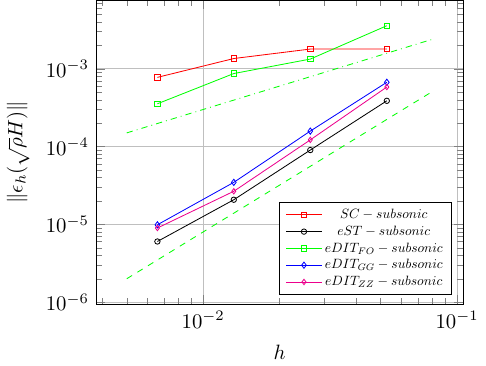}}
\qquad
\subfloat[$\sqrt{\rho}u$]{\includegraphics[width=0.5\textwidth]{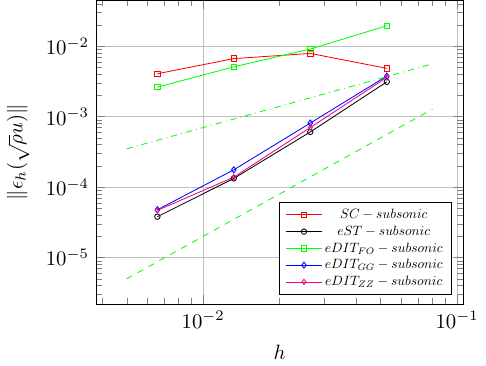}}
\subfloat[$\sqrt{\rho}v$]{\includegraphics[width=0.5\textwidth]{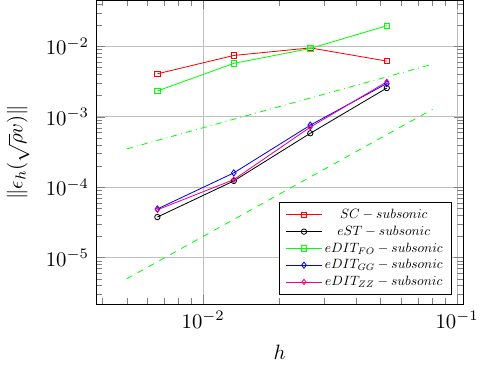}}
	\caption{Transonic source flow: grid-convergence analysis \revAB{within} the shock-downstream (subsonic) region.}\label{Q1D-subsonic}
\end{figure*}

\begin{figure*}
\centering
\subfloat[$eST$]{\includegraphics[width=0.33\textwidth]{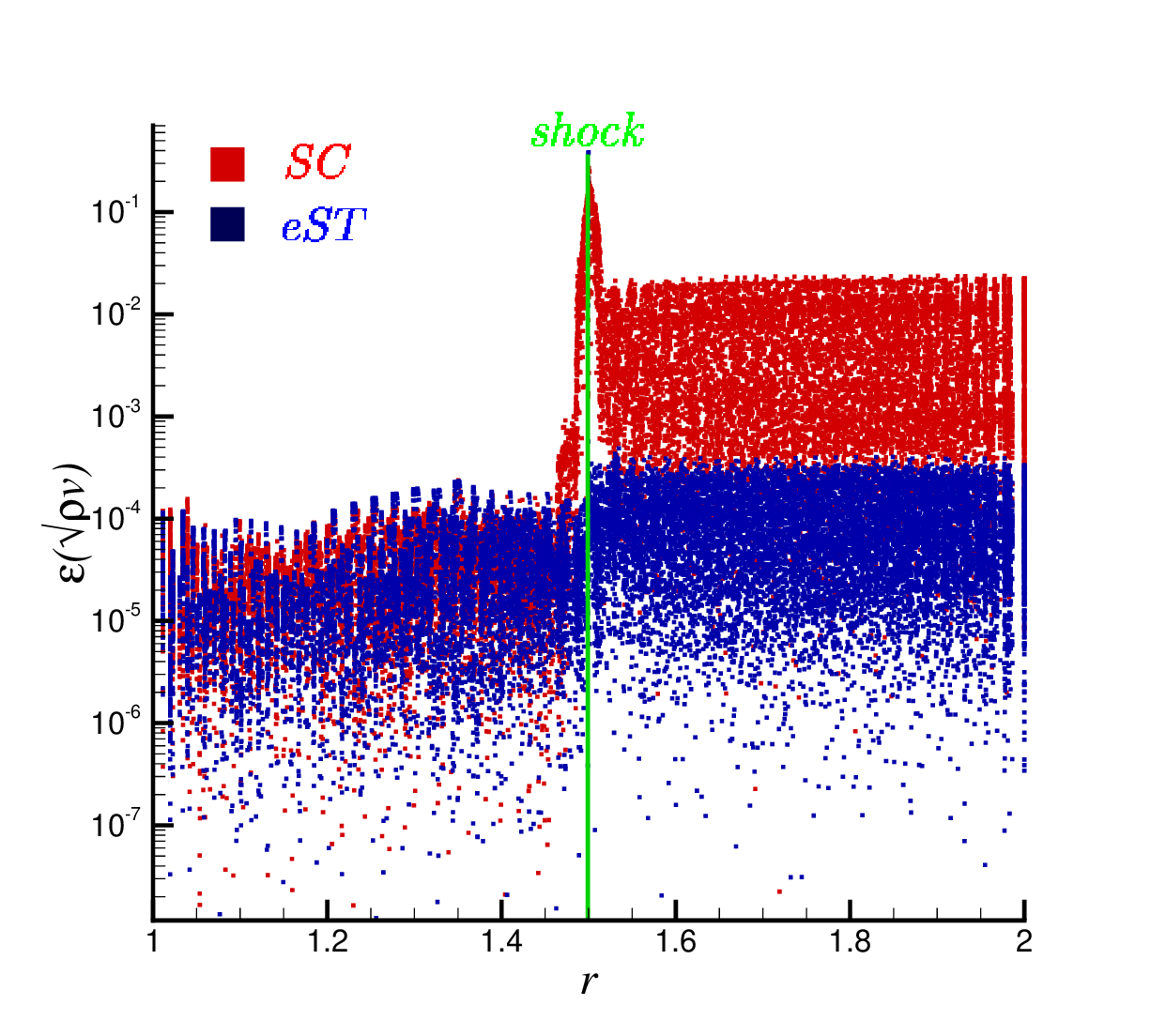}}
\subfloat[$eDIT_{GG}$] {\includegraphics[width=0.33\textwidth]{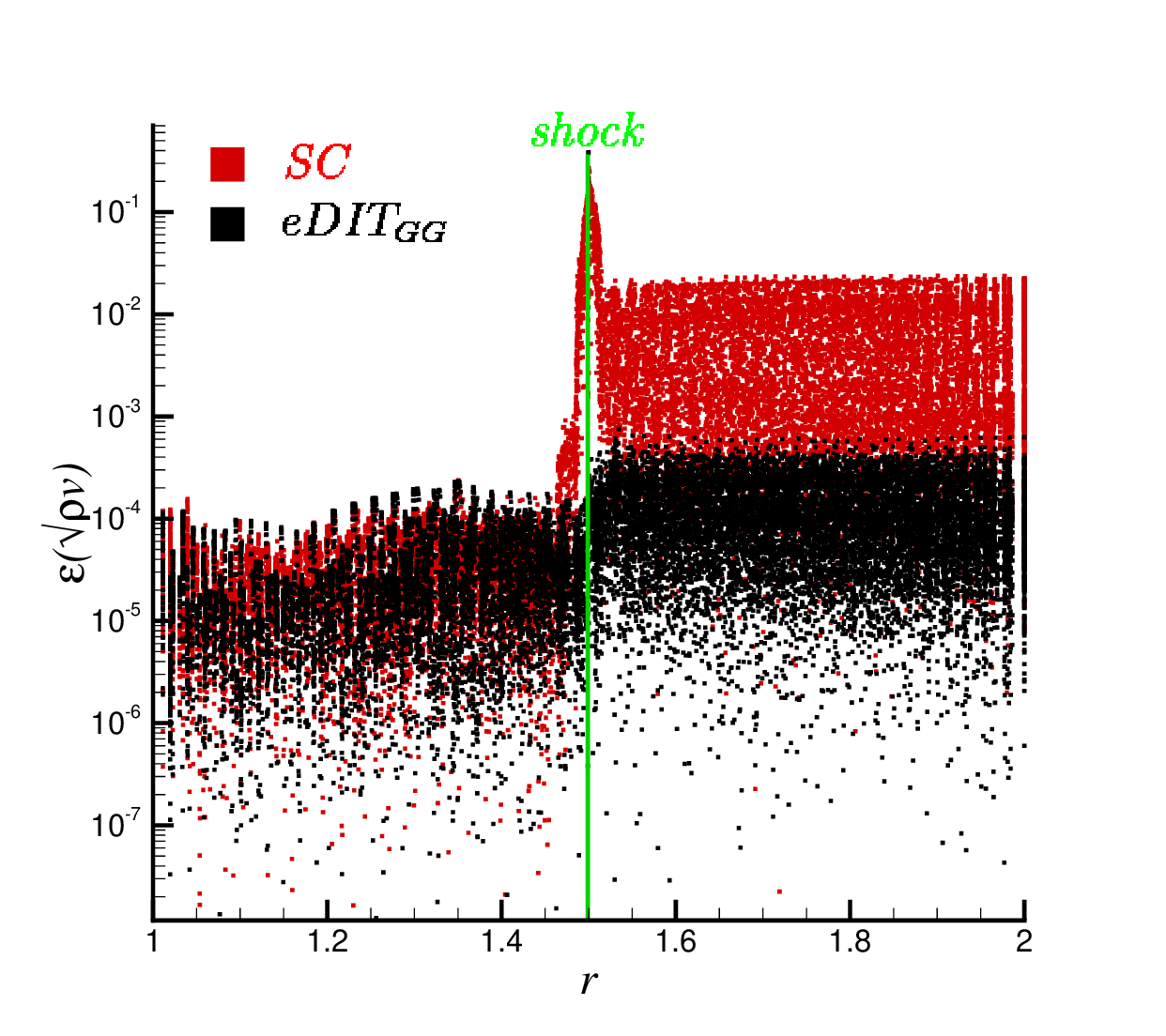}}
\subfloat[$eDIT_{ZZ}$] {\includegraphics[width=0.33\textwidth]{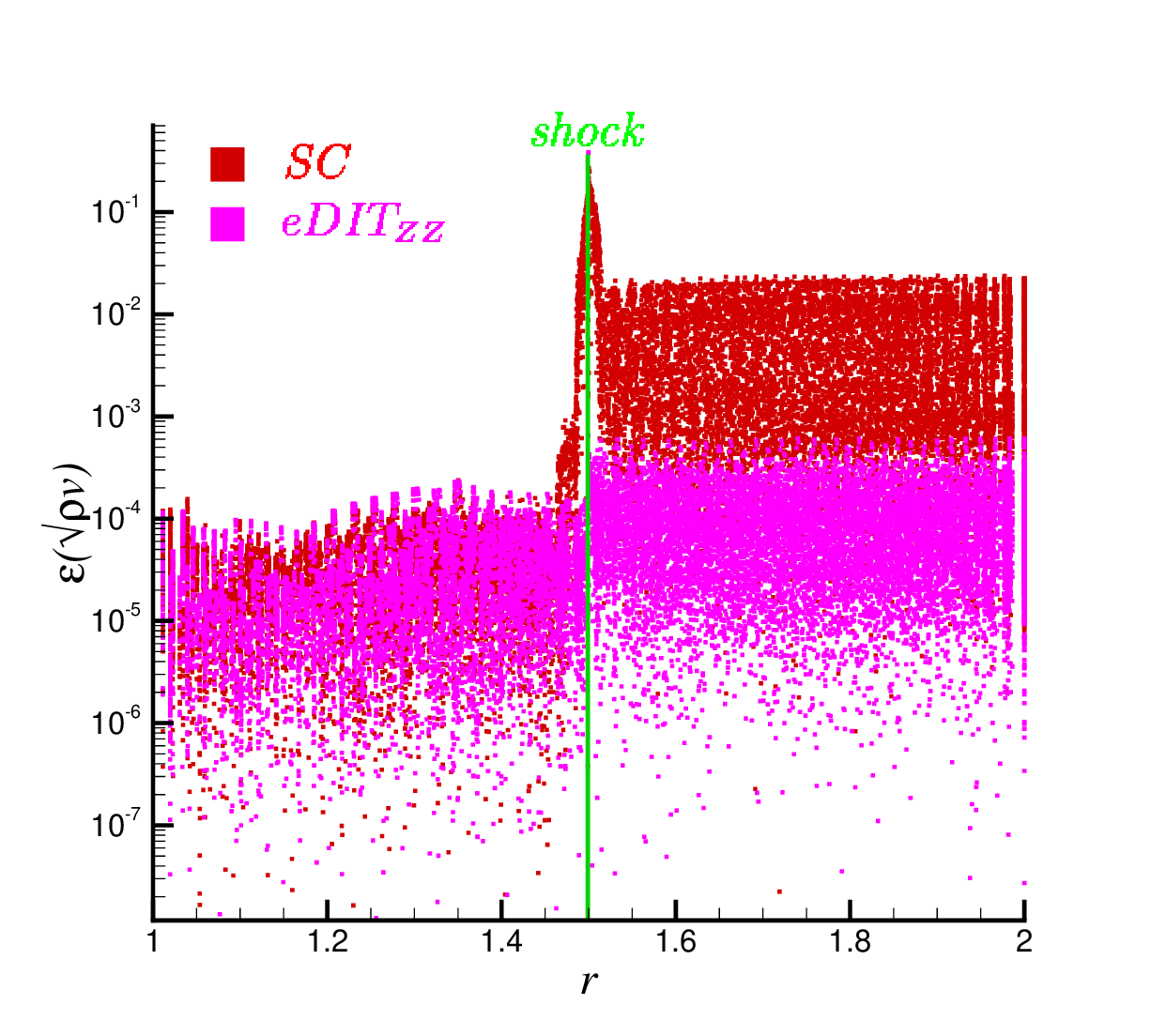}}
\caption{Transonic source flow: local discretization error analysis carried out on \textit{grid level 2} to compare the $SC$ computations to the $eDIT$ ones w.r.t.\ $Z_4=\sqrt{\rho}v$.}\label{Q1D-localerror}
\end{figure*}

\begin{figure}
\centering
\subfloat{%
\includegraphics[width=0.50\textwidth]{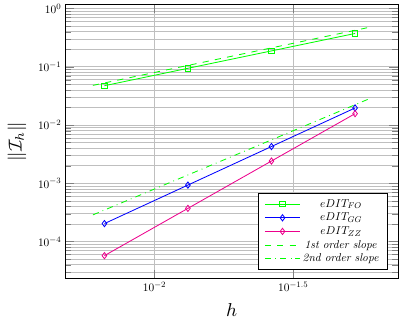}}
\caption{Transonic source flow: grid-convergence analysis on the flux balance between $\tilde\Gamma_U$ and $\tilde\Gamma_D$.}\label{q1d-mass}
\end{figure}

\begin{figure*}
\centering
	\subfloat[$SC$]{\includegraphics[width=0.45\textwidth]{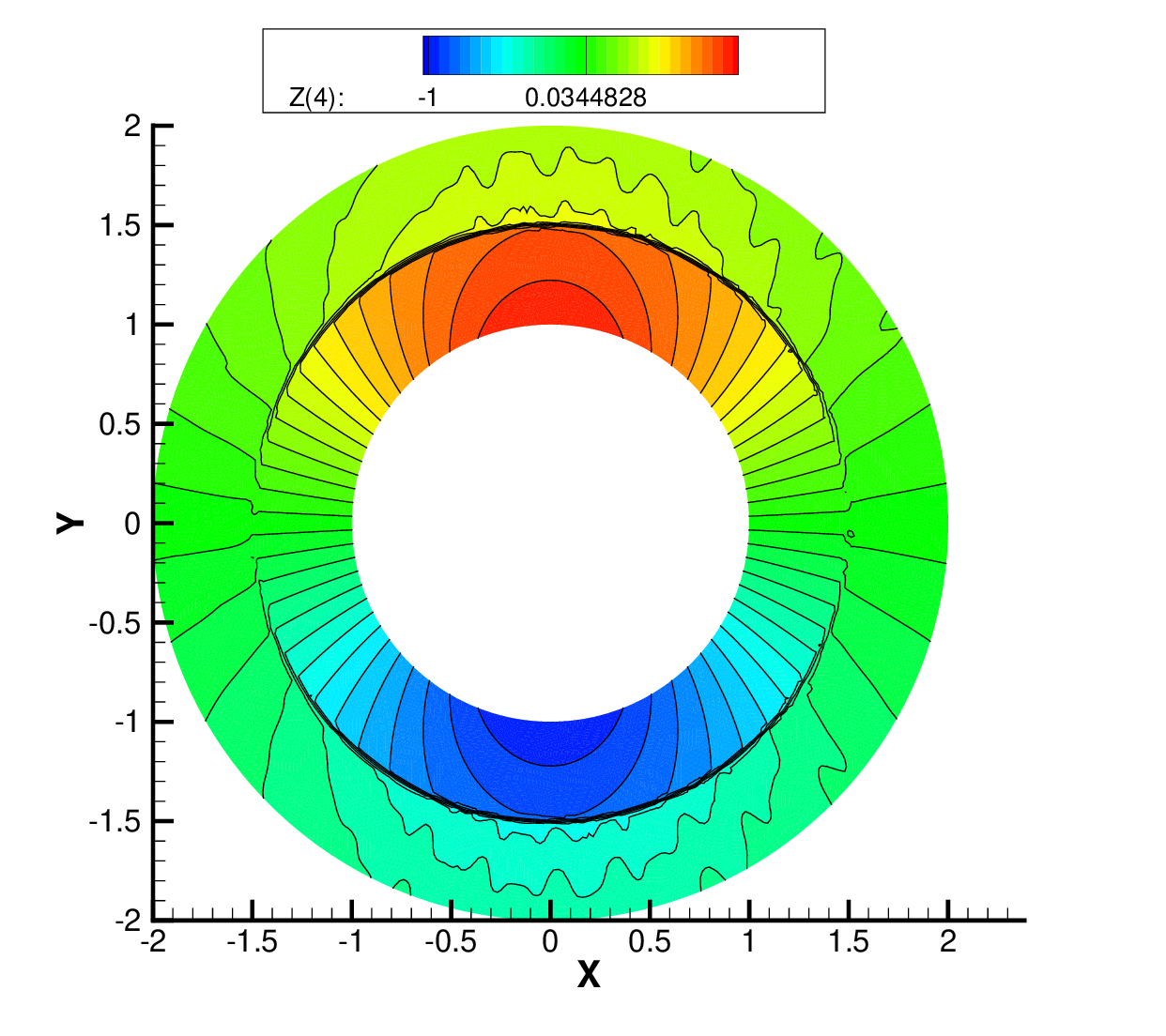}\label{Q1D-SCvsFOvsGGa}}
\subfloat[$eDIT_{GG}$ vs.\ $eDIT_{FO}$] {\includegraphics[width=0.45\textwidth]{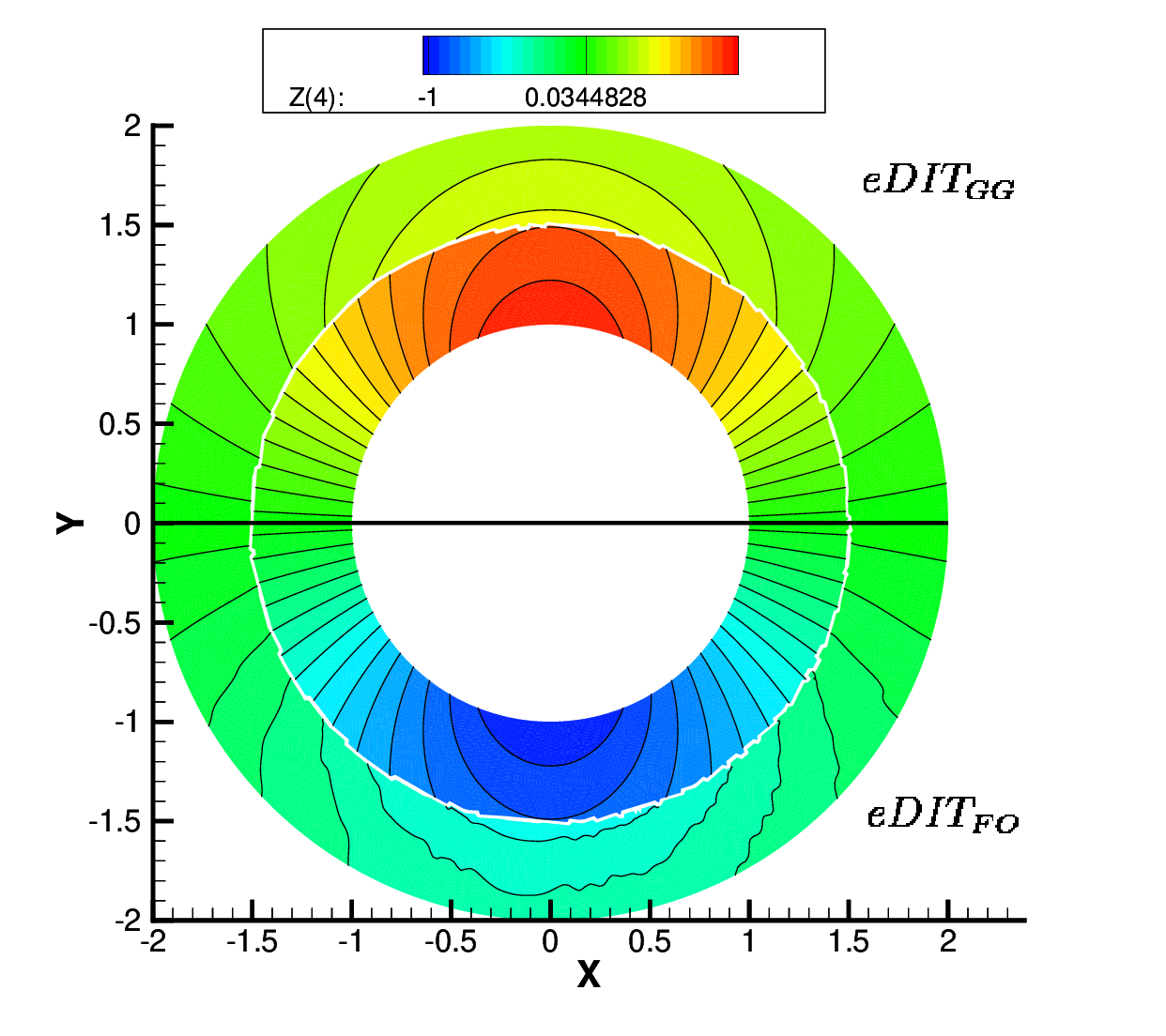}\label{Q1D-SCvsFOvsGGb}}
\caption{Transonic source flow: iso-contours comparison on \textit{grid level 2} w.r.t.\ $Z_4=\sqrt{\rho}v$.}\label{Q1D-SCvsFOvsGG}
\end{figure*}

Firstly, the global errors obtained on the shock-upstream side of the domain have been displayed to corroborate the fact that,
for the supersonic region, the results \revAB{obtained} with all methods \revAB{are} almost identical. 
Indeed, Fig.~\ref{Q1D-supersonic} shows that, on the shock-upstream side, all approaches are able to retain 
the formal order of accuracy of the method.
Instead, when looking at the results shown in Fig.~\ref{Q1D-subsonic}, 
\revAB{which shows the global errors within the subsonic region},
the situation is \revAB{considerably} different.
\revRP{An almost asymptotic convergence is observed}
for all second order versions of the extrapolated tracking method ($eST$, $eDIT_{GG}$ and $eDIT_{ZZ}$).
\revAB{Instead}, for $SC$, the convergence rate drops to one, and even less for the coarsest meshes.
$eST$ and $eDIT_{GG}$ solutions are almost indistinguishable,
with global errors that \revMR{favor} 
more the first approach, 
\revMR{probably due to the use of flow data explicitly accounting for the jump conditions also in the extrapolation. The }  
$eDIT_{ZZ}$ is still comparable having a trend in the middle between the two. 
Although the differences are not remarkable, using a more accurate reconstruction seems \revMR{to have some effect on the magnitude of the error. The improvement is however so small that the }
 $eST$ can be \revMR{safely} replaced by either \revMR{$eDIT_{GG}$ or $eDIT_{ZZ}$. This also shows that our basic idea can be coupled to different extrapolation methods, and others, possibly improved ones, could be suggested in the future.}
As expected, $eDIT_{FO}$ performs as the others in the supersonic side of the domain showing a second-order trend.
However in the subsonic side, due to the first order extrapolations carried out at the shock, 
its trend follows the first order slope curve, because it does not include the \revAB{nodal} gradients
\revAB{in the extrapolation}.\\
\revMC{Figure~\ref{q1d-mass} points out the flux balance computed on the two surrogate boundaries $\tilde\Gamma_U$ and $\tilde\Gamma_D$ showing that convergence
trends strictly depends on the order of the extrapolation. 
It is indeed shown that first (second) order behaviour is displayed for first (second) order extrapolations, confirming the arguments  of  \S~\ref{remarksconservation}. 
}  
%
Finally, Fig.~\ref{Q1D-localerror}, \revAB{which plots} the local discretization error \revAB{against
the radial distance $r$}, clearly reveals
\revRP{the huge reduction of the numerical errors}
that the various $eDIT$ approaches \revMR{provide} w.r.t.\ $SC$. 
The green line in Fig.~\ref{Q1D-localerror} represents the actual position where the shock occurs ($r=1.5$).
As \revMR{already verified in} Fig.~\ref{Q1D-supersonic}, 
the error computed with \revAB{either} $SC$ \revAB{or any of the} $eDIT$ versions \revAB{within} the shock-upstream 
region is almost \revAB{the same}.
\revAB{Iso-contour lines of the fourth component of $\mathbf{Z}$
computed on \textit{grid level 2} are displayed in Fig.~\ref{Q1D-SCvsFOvsGG}}:
\revAB{the $SC$ solution is shown in Fig.~\ref{Q1D-SCvsFOvsGGa}, whereas Fig.~\ref{Q1D-SCvsFOvsGGb}
shows the results of two different $eDIT$ calculations. More precisely, the second-order-accurate
$eDIT_{GG}$ result is displayed in the upper half of Fig.~\ref{Q1D-SCvsFOvsGGb} and
$eDIT_{FO}$ is shown in the lower half of the same frame}.
\revAB{When mutually comparing the three different calculations, it can be seen that
both the $SC$ and $eDIT_{FO}$ calculations} are affected by~\revAB{severe} oscillations
downstream of the shock, which completely disappear \revAB{when using} $eDIT_{GG}$.
\revAB{Even though the $eDIT_{FO}$ solution looks slightly better
than the $SC$ one,
Fig.~\ref{Q1D-subsonic} confirms that $eDIT_{FO}$}
\revMR{cannot be} better than first-order accurate \revAB{behind the shock}.
%
%

        %
        
\subsection{Hypersonic flow past a blunt body}\label{blunt body}
\begin{figure}
\begin{center}
\subfloat[Sketch of the computational domain.]{\includegraphics[width=0.35\textwidth,trim={0cm 0cm 10cm 2cm},clip]{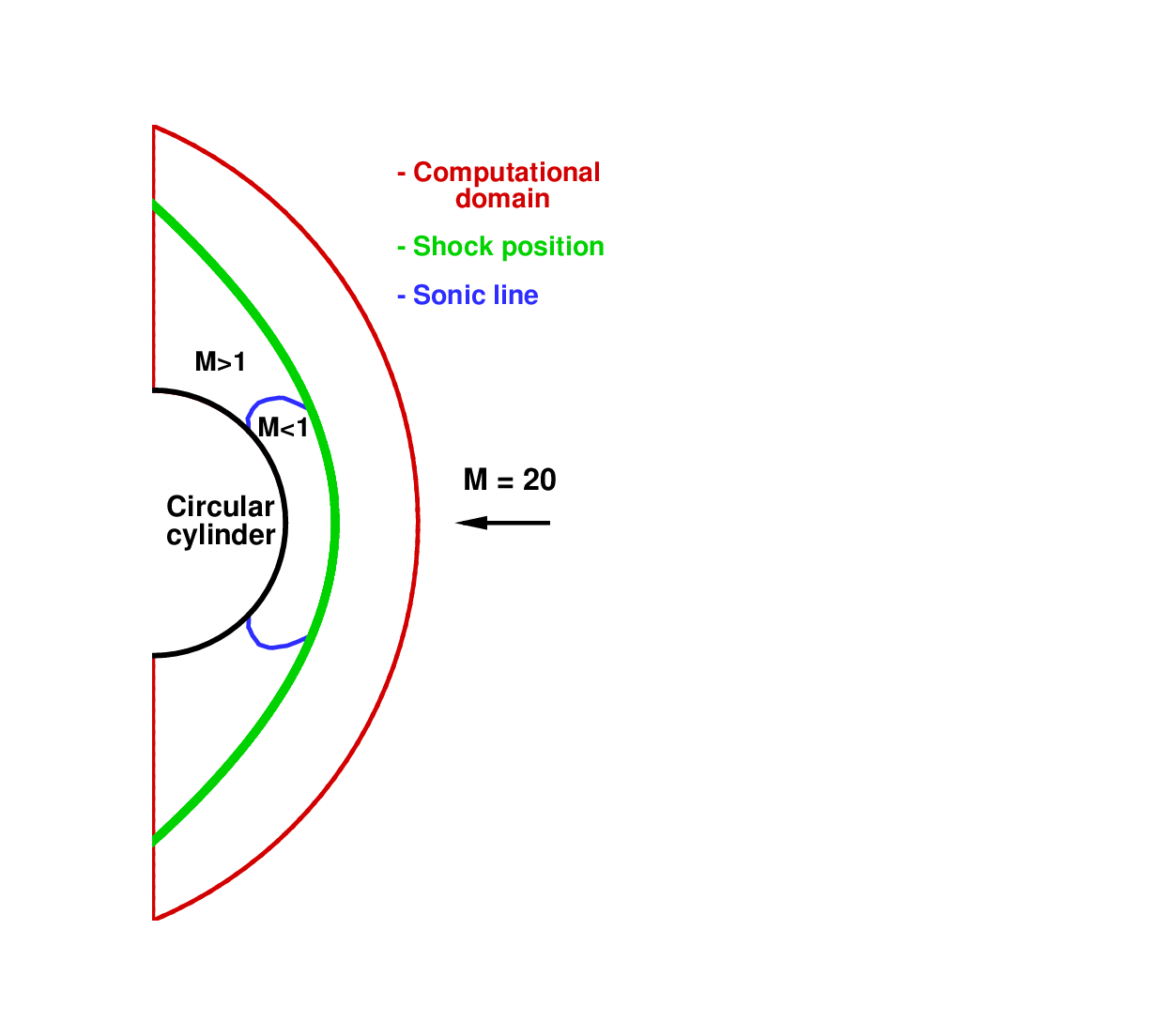}\label{CCsketch}}
\qquad
\subfloat[Detail of the level 0 unstructured mesh.]{\includegraphics[width=0.35\textwidth,trim={0cm 0cm 10cm 2cm},clip]{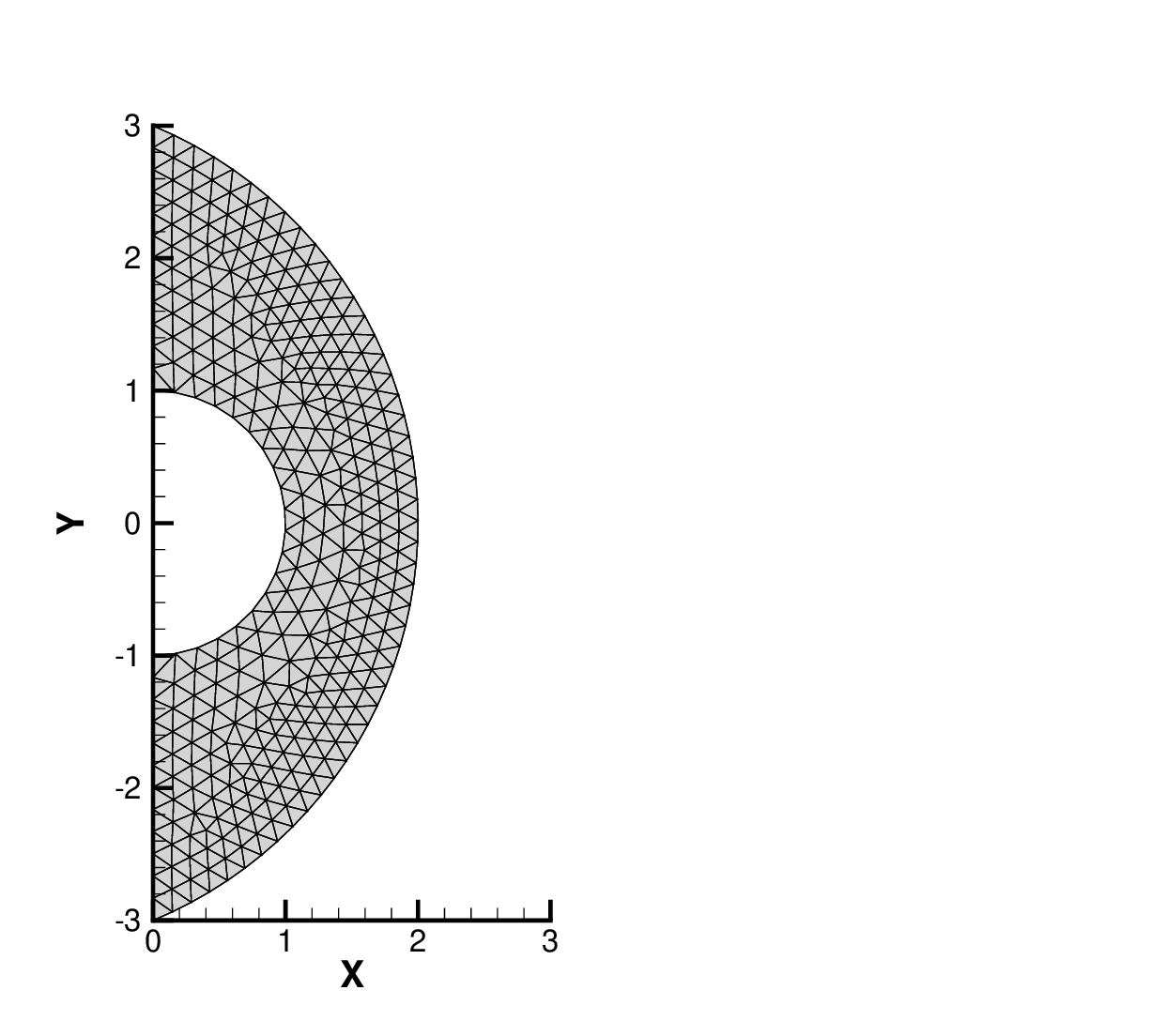}\label{CCmesh}}
\end{center}
\caption{Hypersonic flow past a blunt body.}\label{pic23}
\end{figure}

%
\revMR{We consider now a}
hypersonic ($M_\infty=$ 20) 
flow past the forebody of a circular cylinder.
Compared to the previous test-case, this flow is a comprehensive test-bed for the algorithm, 
because the entire shock-polar is swept whilst moving along the bow shock which \revMR{forms} 
ahead of the blunt body.
Moreover, \revAB{as shown in Fig.~\ref{CCsketch}}, the flow-field is 
much more complex \revAB{than that examined in \S~\ref{sourceflow}}, 
due to the presence \revMR{of a stagnation region within a }
subsonic pocket, \revMR{as well as a smooth re-acceleration of the flow to supersonic conditions along the body.} 
Even though no analytical solution is available for \revAB{the entire flow-field}, 
\revMR{we know that total enthalpy, $H$, is preserved along streamlines at steady state.
For a constant profile of $H=H_\infty$ \revAB{ahead of the shock}, this leads to an exact solution \revAB{featuring} 
a homogeneous total enthalpy field,}
\revMR{thus the condition $H=H_{\infty}$ }
can be used to study the convergence behaviour \revAB{of the various combinations of numerical schemes 
and shock-modeling options}. 
Starting from a coarse mesh \revAB{(shown in Fig.~\ref{CCmesh}}) obtained with the \texttt{\ttfamily delaundo} frontal/Delaunay mesh generator~\cite{Muller1,Muller2}, \revMR{we created}
three \revMR{levels of} nested  \revMR{refined meshes with characteristics summarized} 
in Table~\ref{tab8} \revMR{for both the background and computational grids.} 
%
\begin{table}[htb]
\scriptsize
\centering
\caption{Hypersonic flow past a blunt body: characteristics of the background and computational meshes used to perform the grid-convergence tests.}\label{tab8}
\begin{tabular}{cccccc} \hline\hline
        &\multicolumn{2}{c}{\textit{Background grid}} &\textit{Both grids}  &\multicolumn{2}{c}{\textit{Computational grid}} \\[0.5mm] \cline{2-6}
        Grid level & Grid-points & Triangles & $h$    & Grid-points & Triangles \\[0.5mm] \hline
        0          & 351         & 610       & 0.16   & 351         & 535       \\
        1          & 1,311       & 2,440     & 0.08   & 1,311       & 2,292     \\
        2          & 5,061       & 9,760     & 0.04   & 5,061       & 9,460     \\
        3          & 19,881      & 39,040    & 0.02   & 19,881      & 38,439    \\
\hline\hline
\end{tabular}
\end{table}

\begin{figure}
\centering
\includegraphics[width=0.55\textwidth]{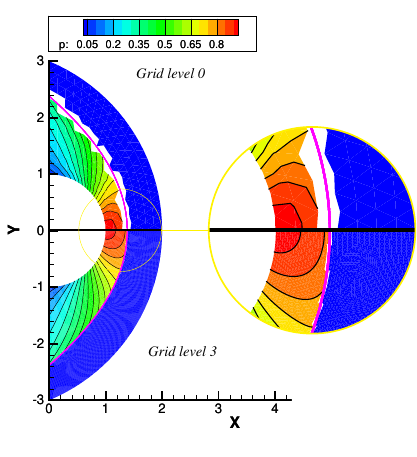}
\caption{Hypersonic flow past a blunt body: comparison between the solutions obtained with $eDIT_{GG}$ on \textit{Grid level 0} and \textit{Grid level 3} in terms of pressure iso-lines.}\label{CCcomp}
\end{figure}

\revMR{The different extrapolated}
tracking methods \revMR{discussed in the previous sections} \revAB{are expected to} 
provide orders of convergence similar to those observed for the source flow \revAB{of \S~\ref{sourceflow}}, 
\revMR{as well as smooth and clean results with very low numerical perturbations}, \revAB{even on the coarsest grids};
\revAB{this latter aspect is illustrated in Fig.~\ref{CCcomp}, where the pressure iso-contour lines 
computed using $eDIT_{GG}$ on the
level 0 (coarsest) and level 3 (finest) grids are mutually compared}.
%
\begin{figure}
\centering
\subfloat[Total enthalpy error]{%
\includegraphics[width=0.48\textwidth]{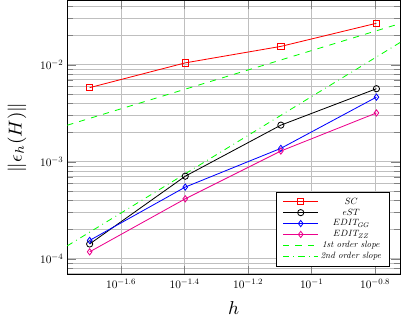}}
\subfloat[Flux balance]{%
\includegraphics[width=0.48\textwidth]{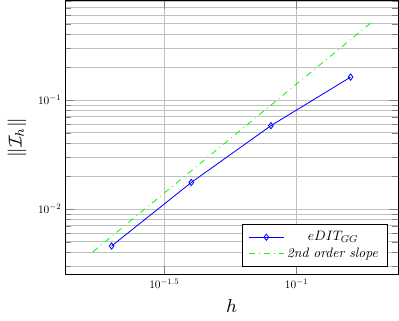}}
\caption{Hypersonic flow past a blunt body: grid-convergence analysis.}\label{pic28}
\end{figure}
Figure~\ref{pic28}a shows the grid-convergence analysis 
pointing out again \revMR{that the different extrapolation techniques provide  error levels  relatively close, and a second order slope, 
while the captured result converges with a less than first order rate, and errors of one or two orders of magnitude larger.} 
\revMR{The patch super-convergent gradient recovery }
($ZZ$) \revMR{provides a better result for this case}.
\revMC{Figure~\ref{pic28}b shows the effect of a second order extrapolation, $eDIT_{GG}$, on the flux balance within the cavity for this test-case. 
As discussed  in \S~\ref{remarksconservation} we  recover the error associated to the second order extrapolation.
Similar results for other second order extrapolations are straightforward to obtain and therefore not included in Fig.~\ref{pic28}b.}\\

\begin{figure*}
\centering
\subfloat[$SC$, \textit{Grid level 0}]{%
\includegraphics[width=0.25\textwidth,trim={1cm 0cm 8cm 0cm},clip]{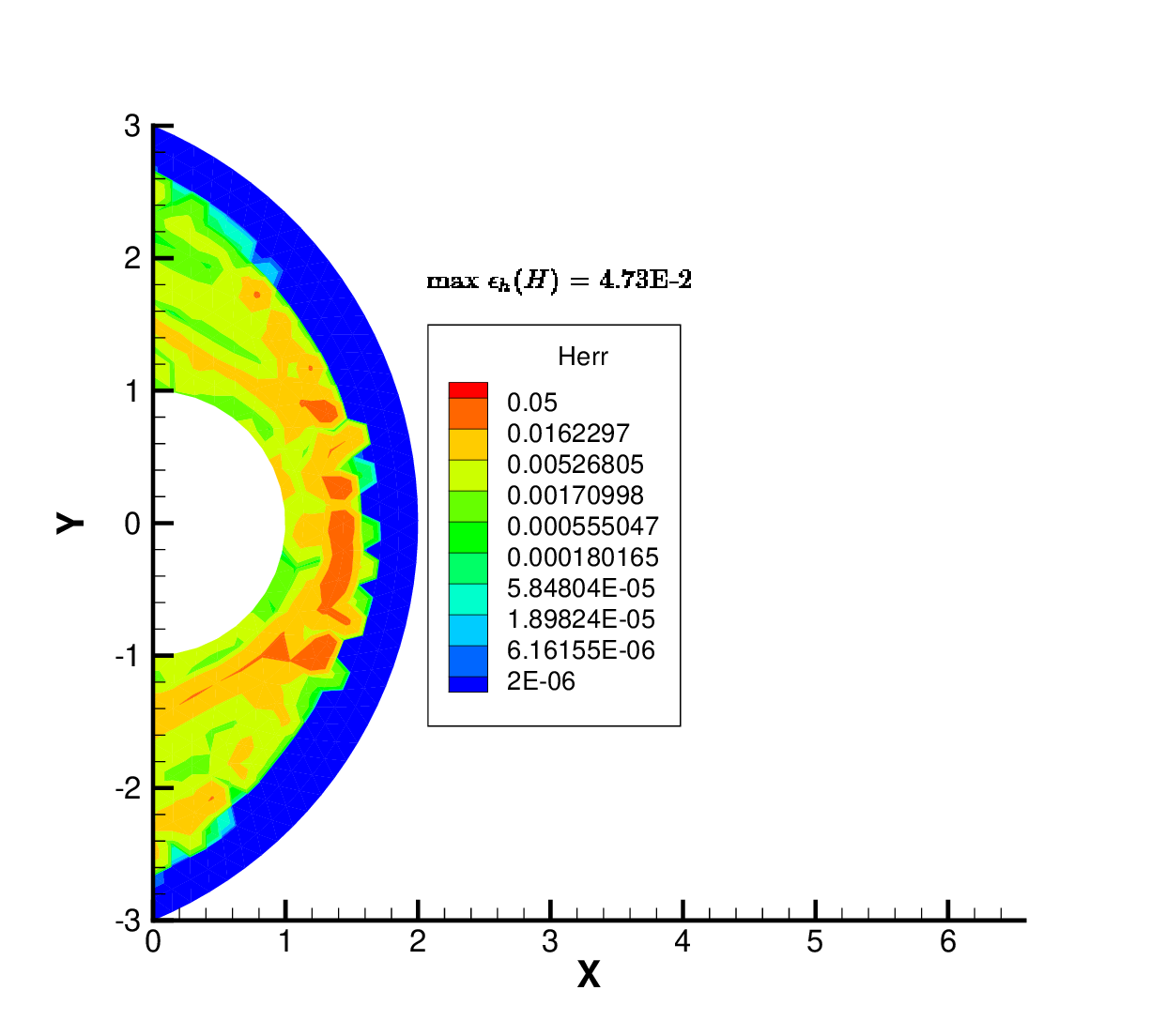}\label{HlocalerrorSCvsGGa}}
\subfloat[$eDIT_{GG}$, \textit{Grid level 0}]{%
\includegraphics[width=0.25\textwidth,trim={1cm 0cm 8cm 0cm},clip]{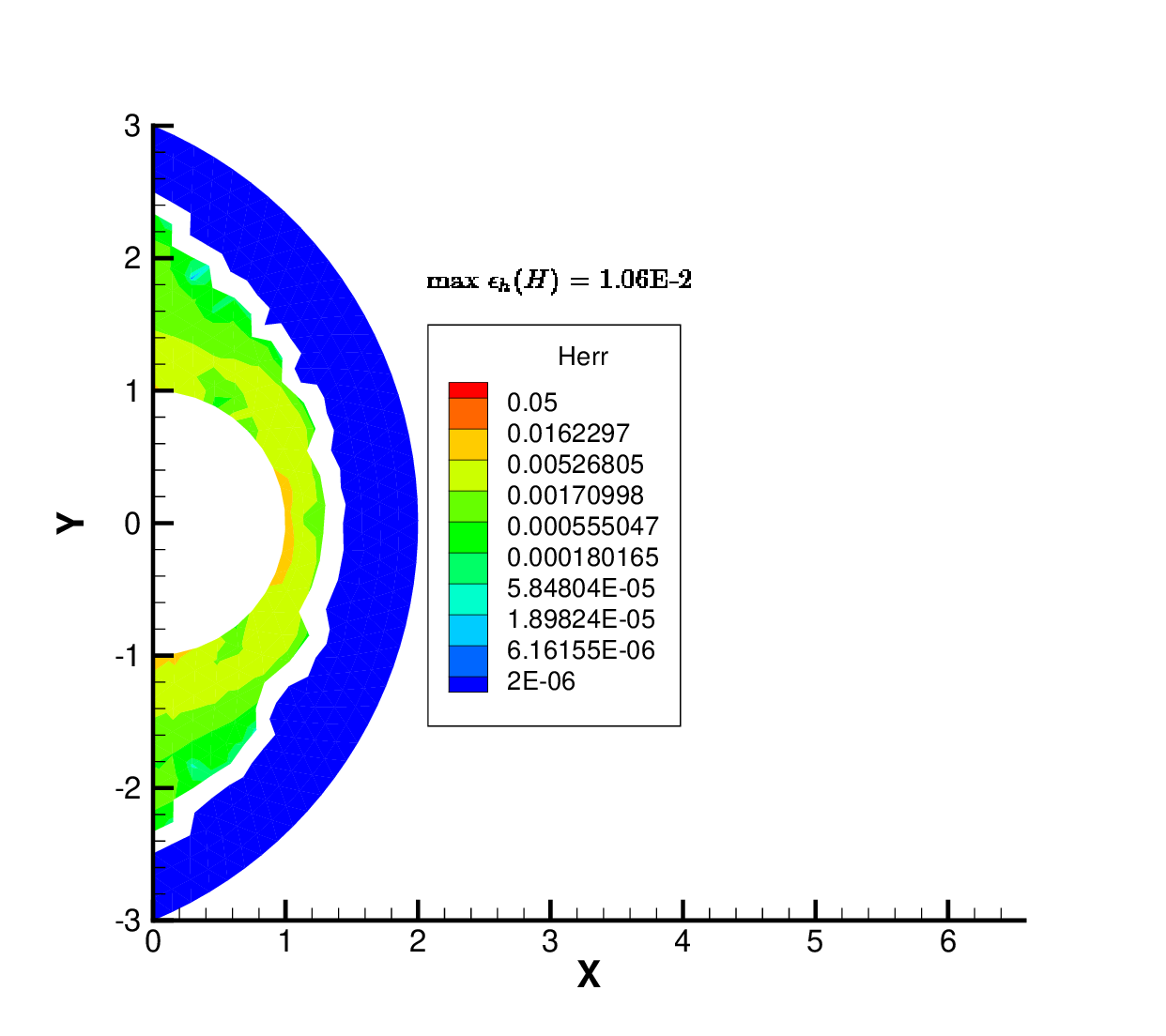}\label{HlocalerrorSCvsGGb}}
\subfloat[$SC$, \textit{Grid level 3}]{%
\includegraphics[width=0.25\textwidth,trim={1cm 0cm 8cm 0cm},clip]{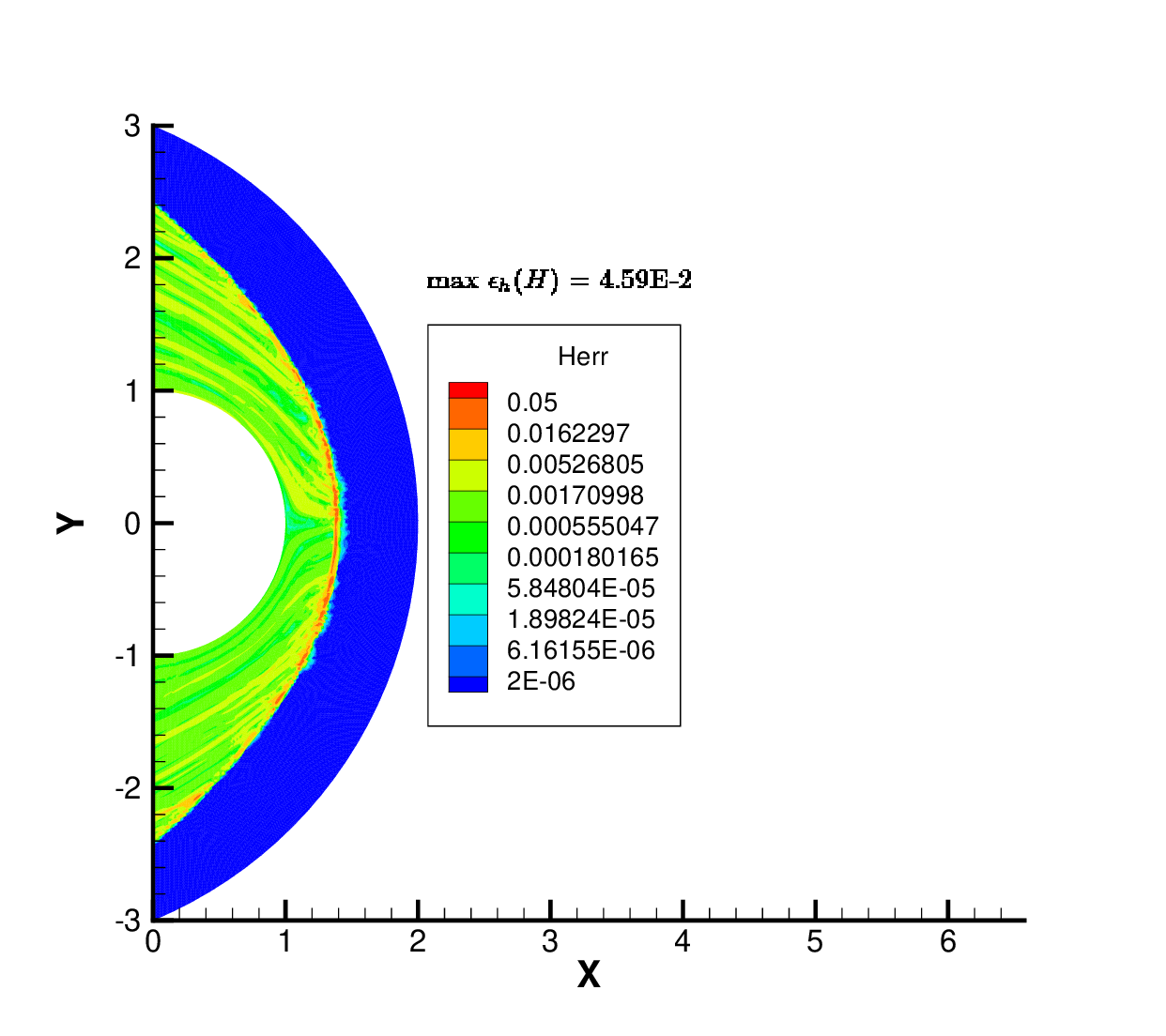}\label{HlocalerrorSCvsGGc}}
\subfloat[$eDIT_{GG}$, \textit{Grid level 3}]{%
\includegraphics[width=0.25\textwidth,trim={1cm 0cm 8cm 0cm},clip]{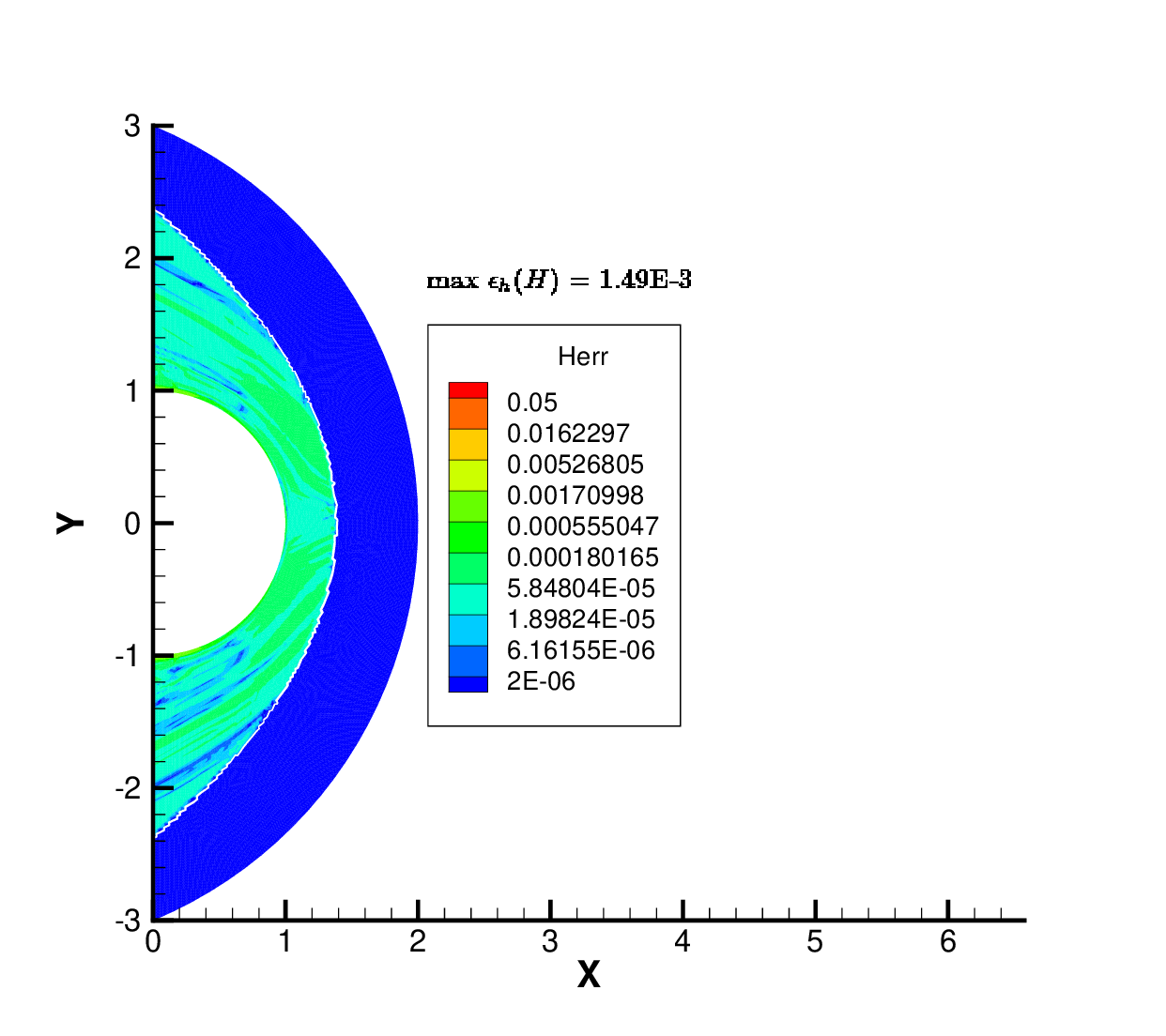}\label{HlocalerrorSCvsGGd}}
\caption{Hypersonic flow past a blunt body: discretization error of the total enthalpy, $\epsilon_h(H)$, obtained with $SC$ and $eDIT_{GG}$ over the entire computational domain for \textit{Grid level 0} and \textit{Grid level 3}.\label{HlocalerrorSCvsGG} }
\end{figure*}

Figure~\ref{HlocalerrorSCvsGG} \revAB{shows the total-enthalpy discretization error 
over the entire computational domain}
for \revAB{the two} solutions obtained with $SC$ and $eDIT_{GG}$ and \revAB{two different grid-levels: the
coarsest and the finest}. 
\revMR{The maximum value of the error is also explicitly given in \revAB{each of the four frames}.
While confirming the
lower error levels on the coarse mesh, and the rapid error convergence obtained with second-order shock-tracking, 
the \revAB{four frames of Fig.~\ref{HlocalerrorSCvsGG}} allow to visually
highlight the substantial difference in error generation.
Indeed, in both captured solutions (Fig.~\ref{HlocalerrorSCvsGGa} and~\ref{HlocalerrorSCvsGGc}), 
the largest error is generated \revAB{along} the shock. 
Conversely, for the two $eDIT_{GG}$ \revAB{calculations} the numerical error is essentially \revAB{connected with}
the wall boundary condition, and undoubtedly 
related to the entropy generated at the stagnation point and advected downstream.
The $eDIT_{ZZ}$ solution is virtually identical and not discussed for brevity.}

\begin{figure}
\centering
\subfloat[$eST$]{%
\includegraphics[width=0.33\textwidth,trim={0cm 0cm 0cm 2cm},clip]{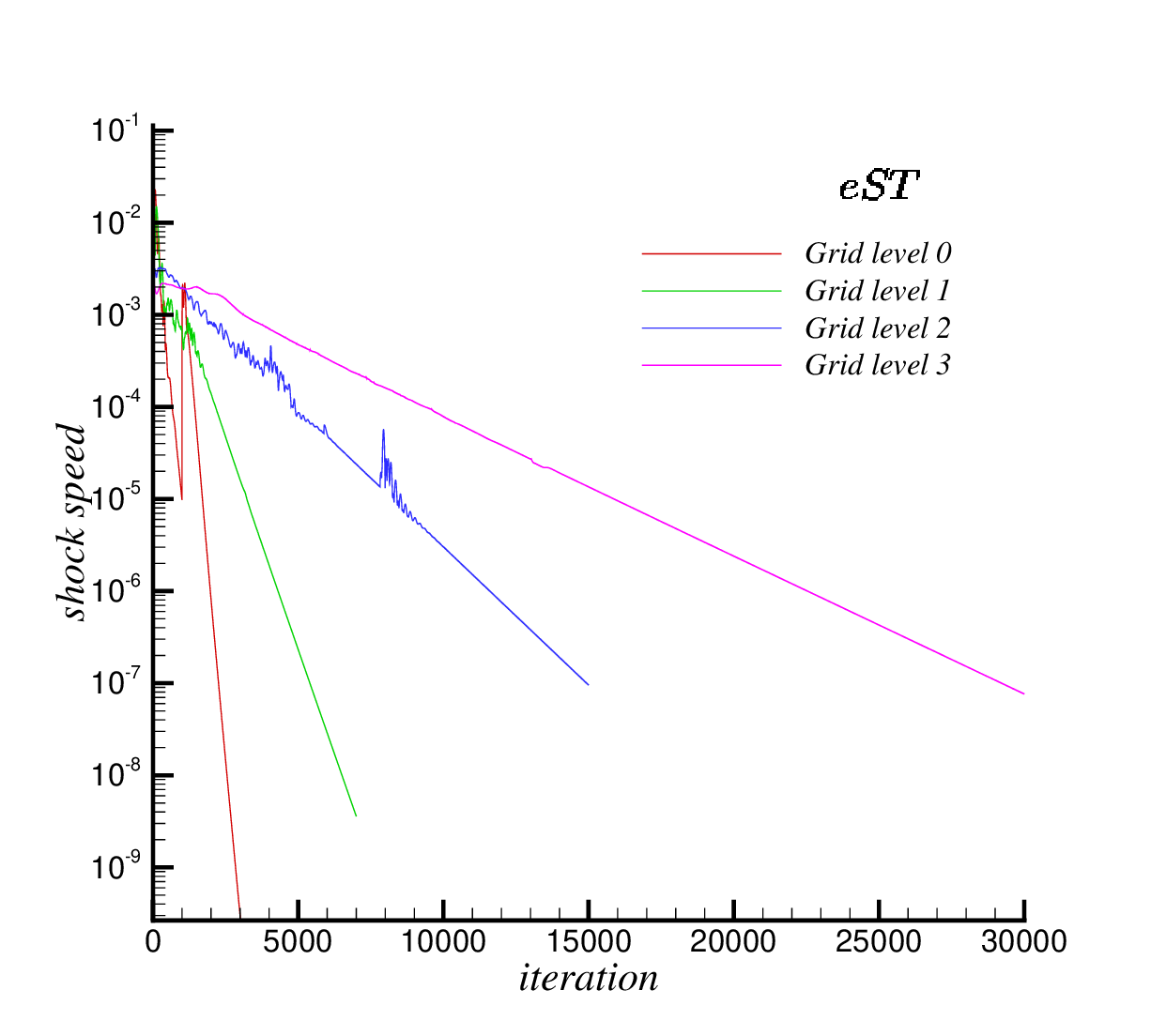}}
\subfloat[$eDIT_{GG}$]{%
\includegraphics[width=0.33\textwidth,trim={0cm 0cm 0cm 2cm},clip]{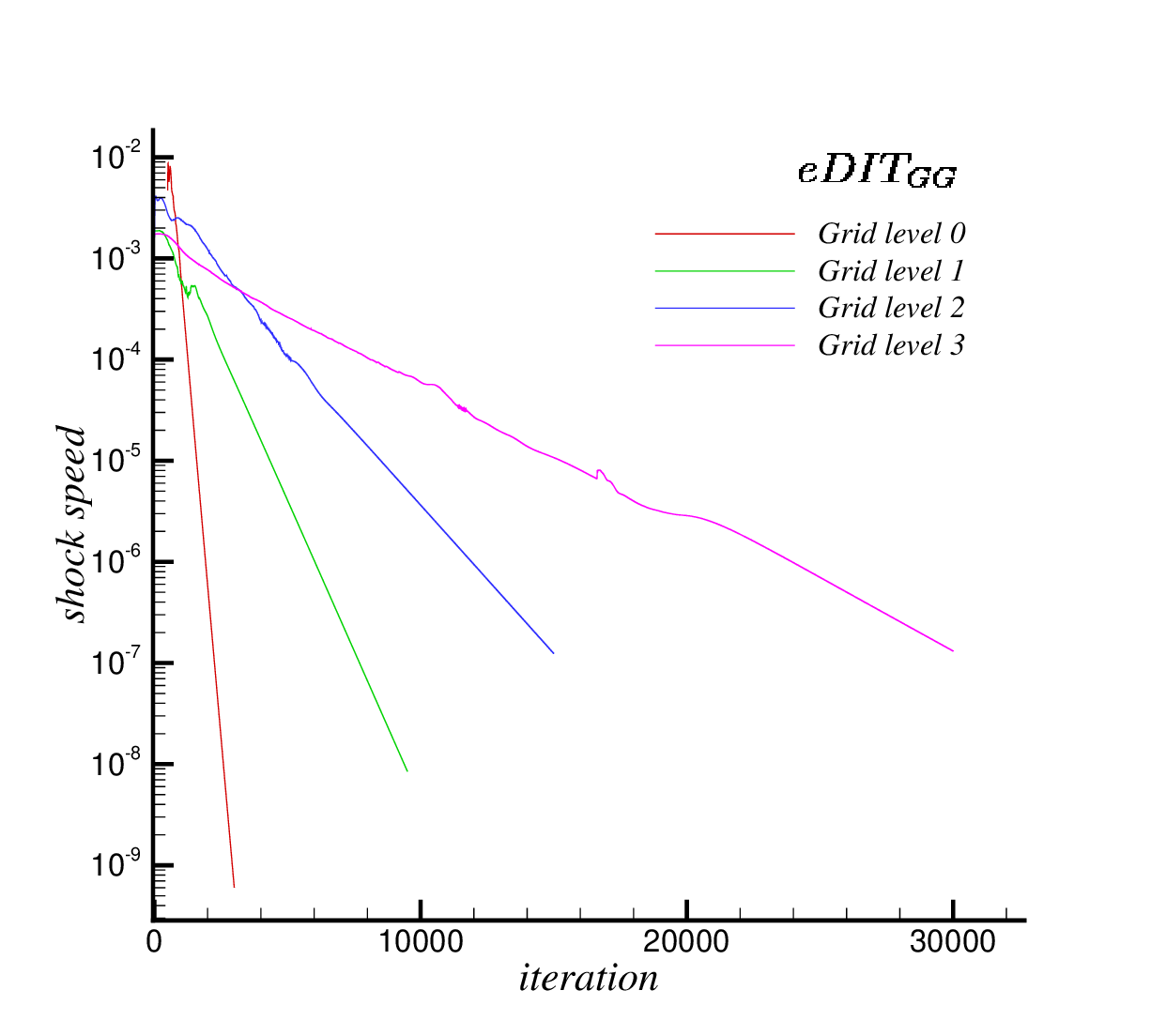}}
\subfloat[$eDIT_{ZZ}$]{%
\includegraphics[width=0.33\textwidth,trim={0cm 0cm 0cm 2cm},clip]{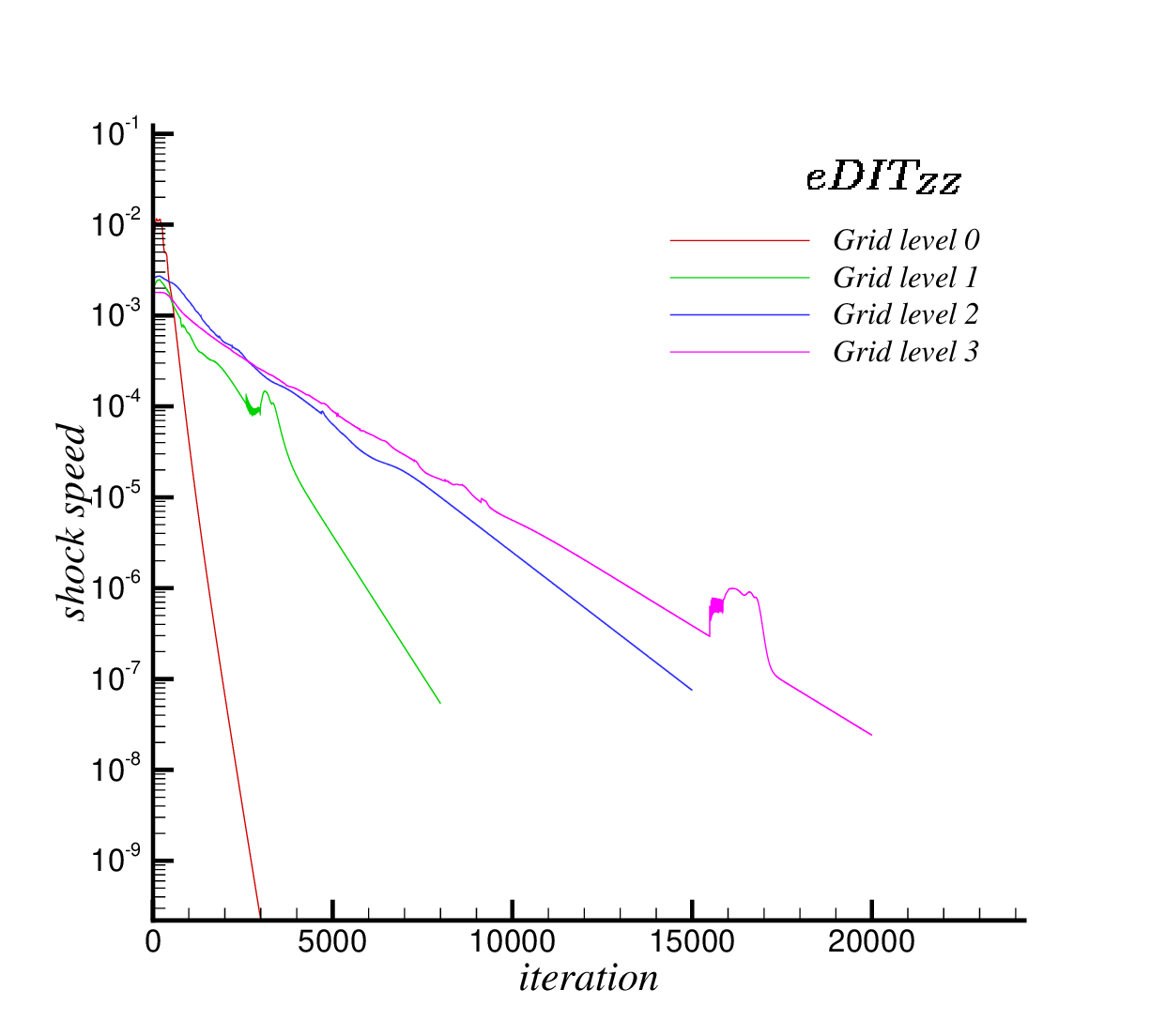}}
\caption{Hypersonic flow past a blunt body: Pseudo-time shock convergence.}\label{pic26}
\end{figure}

\revMR{As a final note, we remark that total} enthalpy is a conserved variable if and only if the time 
derivative \revAB{in Eq.~(\ref{eq:euler0})}
vanishes. This means that \revMR{the convergence to steady state plays a major role in allowing to 
correctly measure \revAB{the decay rate of the discretization error}.
\revAB{Convergence to steady-state} is an important aspect of this type of methods, which is why we also report on 
Fig.~\ref{pic26} the iterative convergence \revAB{of the $eDIT$ algorithm}
when starting from a captured solution. In all our 
\revAB{calculations} we have set as a stopping criterion a threshold on the norm of }
the shock speed of $\sim 10^{-7}$. 
\revMR{As Fig.~\ref{pic26} shows, this level is reached quite monotonically in all the computations. This despite the shock crossing some nodes during the iterations, which produces a change in topology
of the surrogate discontinuities and of the computational domain which can be seen in some of the peaks appearing locally in the  iterative convergence plots. 
Similar convergence curves are often hard to obtain with non-linear shock capturing methods.} 
%
%
%
         
\subsection{Interaction between two shocks of the same family}\label{SS2}
\revAB{The present and subsequent test-cases address the interaction among different
discontinuities, a kind of flow-topology which had not been addressed in
the first journal appearance~\cite{EST2020} of the \revRP {$eST$} algorithm and thus represent one
of the key novelties of this study}.

\begin{figure*}
\centering
\includegraphics[width=0.6\textwidth]{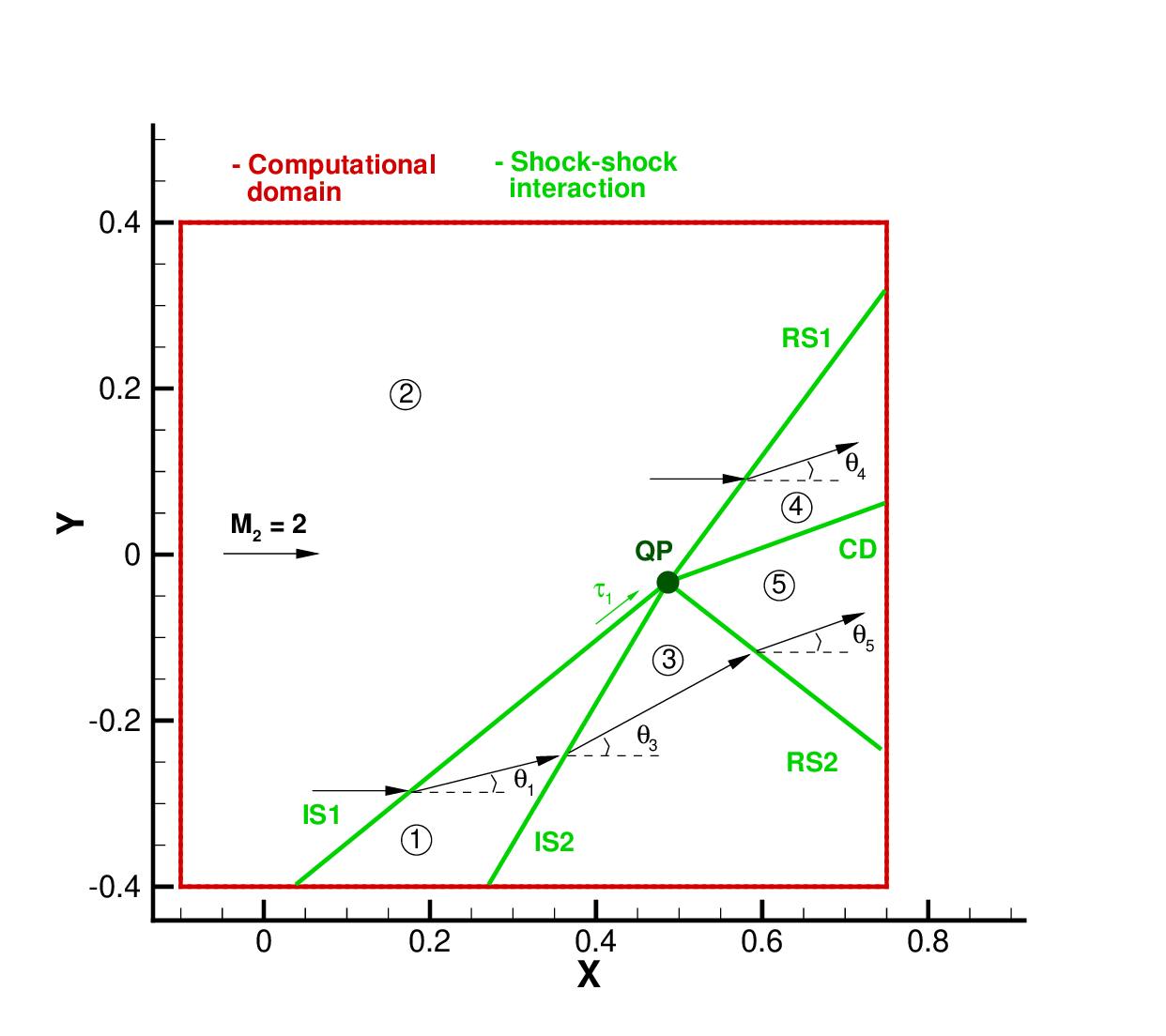}
\caption{Interaction between two shocks of the same family: sketch of the flow.}\label{pic33}
\end{figure*}
\revMR{We consider here the interaction between two oblique shocks}
colliding and giving rise to a five-waves interaction. \revMR{As shown in the sketch of Fig.~\ref{pic33}, }
the  two  incident shocks (IS1 and IS2) \revAB{of the same family}
interact in a quadruple point, referred to as QP in \revMR{the figure, }
generating  \revMR{three reflected discontinuities: a strong reflected} shock 
(RS1), \revMR{with jumps of magnitude larger}
than the two incident ones; a contact discontinuity (CD); 
a \revMR{weak} reflected shock (RS2), \revMR{with jumps so small that it could}
be considered as a \revAB{Mach} wave. 
\revAB{More in general}, depending on the free-stream Mach number and the flow deflection angles
caused by the incident shocks, RS2 \revAB{can} be either a weak expansion or compression wave.
\revMR{With reference to Fig.~\ref{pic33}, in this case} the free-stream Mach number and the flow deflection 
angles are: $M_2 = 2$, $\theta_1 = 10^\circ$ and 
$\theta_2 = 20^\circ$ \revMR{which can be shown to lead to }
states 4 and 5 characterized by $p_4 = p_5 = 2.822$, 
$M_4=1.218$, $M_5=1.28$ and $\theta_4 = \theta_5 = 19.87^\circ$. 
\revMR{For this choice of parameters,  RS2 is so weak that 
we decided not to track this wave since its resolution has negligible overall effect on the flow.}\\

\revMR{Simulations of the interaction have been performed using a Delaunay triangulation 
containing 7,921 nodes and 15,514 triangles.}
\revMR{As before,} the \revMR{same} mesh 
is used to perform the $SC$ simulation and as a background grid for the $eDIT_{GG}$ computations. 
\revMR{Results with the $eDIT_{ZZ}$ are virtually identical to the latter, and not discussed.}
%
\begin{figure*}
\centering
\subfloat[$SC$ solution]              {\includegraphics[width=0.45\textwidth]{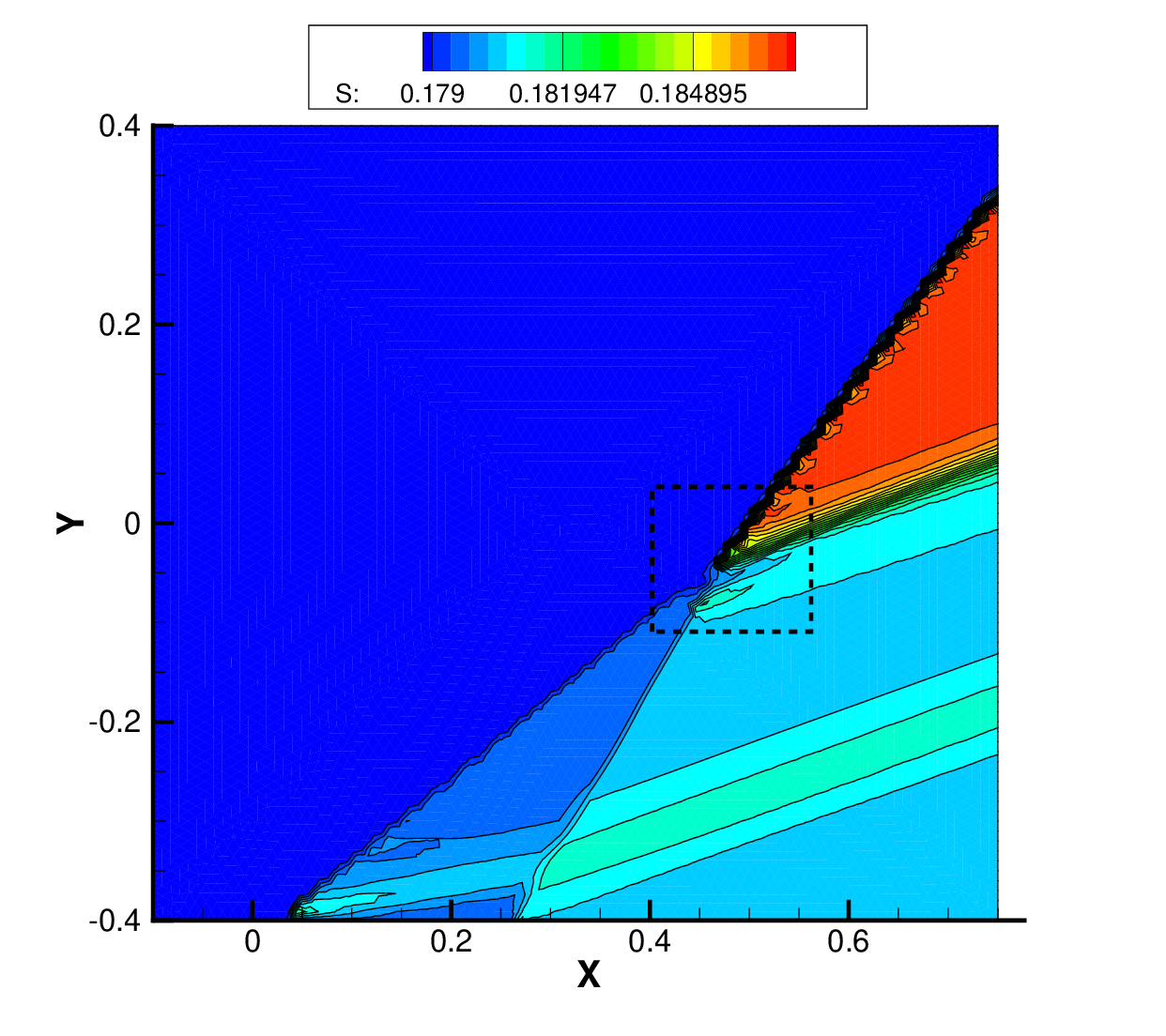}\label{pic34a}}
\subfloat[close-up around the QP]     {\includegraphics[width=0.4\textwidth]{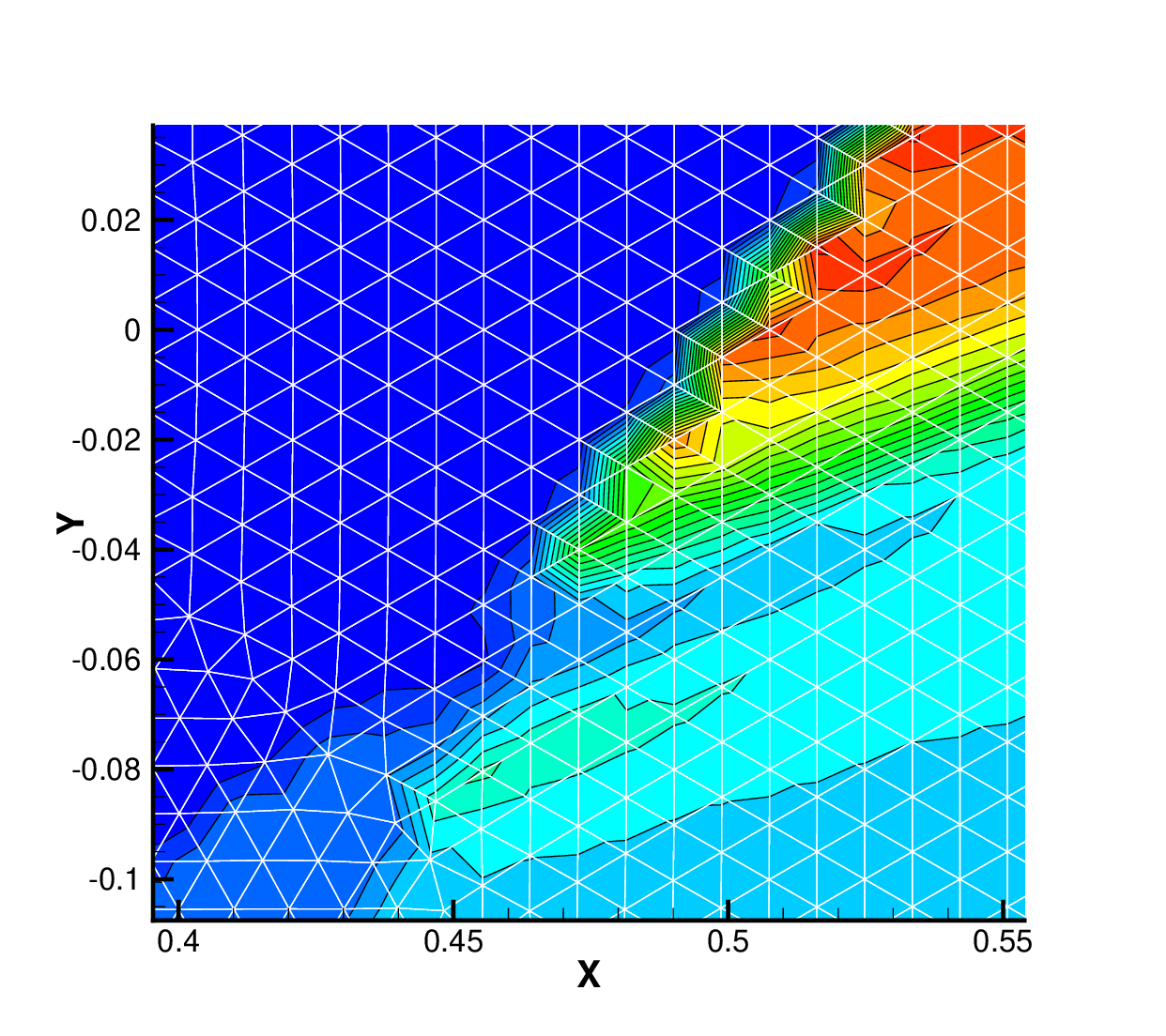}\label{pic34b}}\qquad
	\subfloat[Hybrid $eDIT_{GG}$: \revAB{the CD is captured}.] {\includegraphics[width=0.45\textwidth]{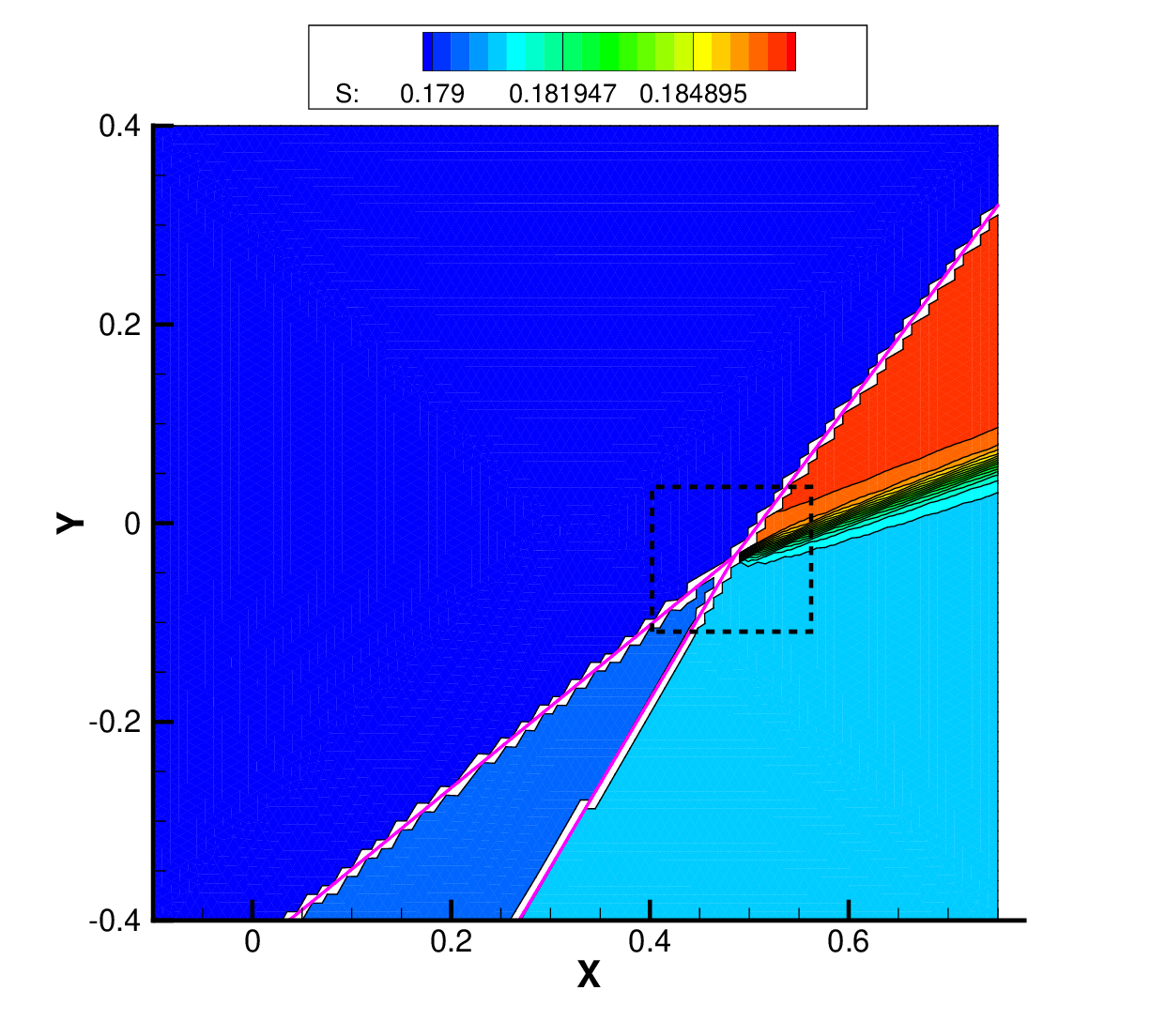}\label{pic34c}}
\subfloat[close-up around the QP]     {\includegraphics[width=0.4\textwidth]{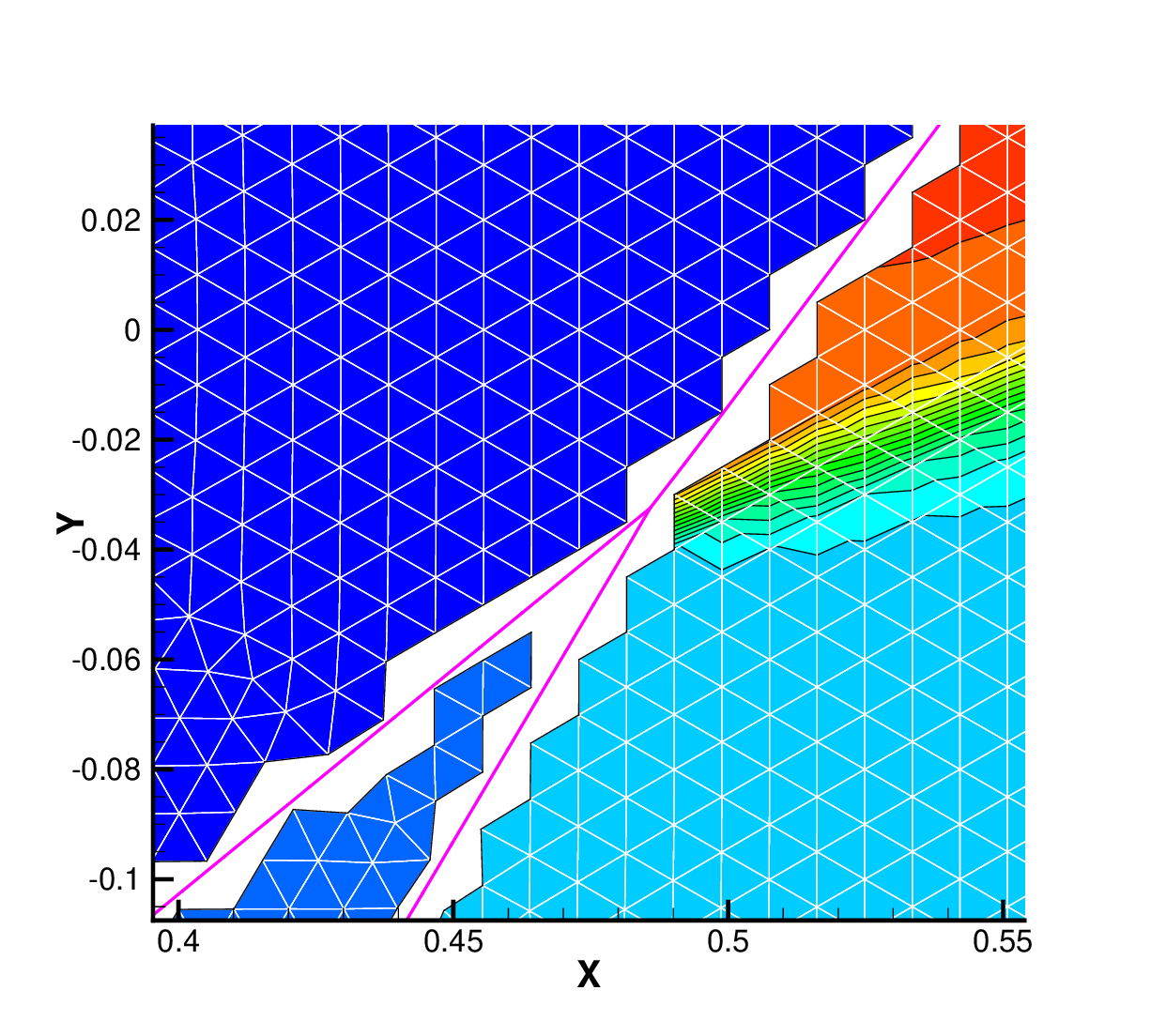}\label{pic34d}}\qquad
	\subfloat[Hybrid $eDIT_{GG}$: \revAB{the CD is tracked}.] {\includegraphics[width=0.45\textwidth]{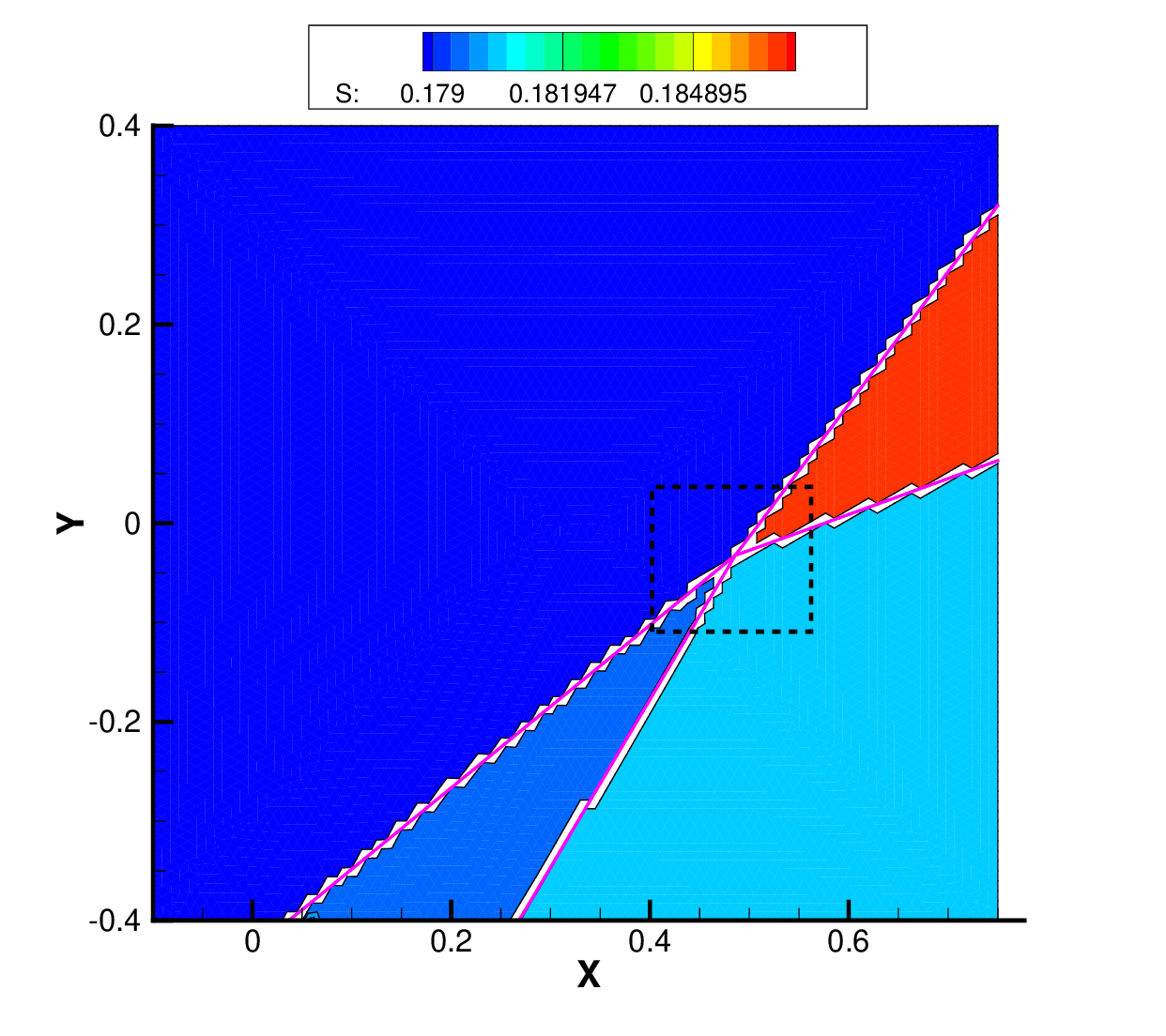}\label{pic34e}}
\subfloat[close-up around the QP]     {\includegraphics[width=0.4\textwidth]{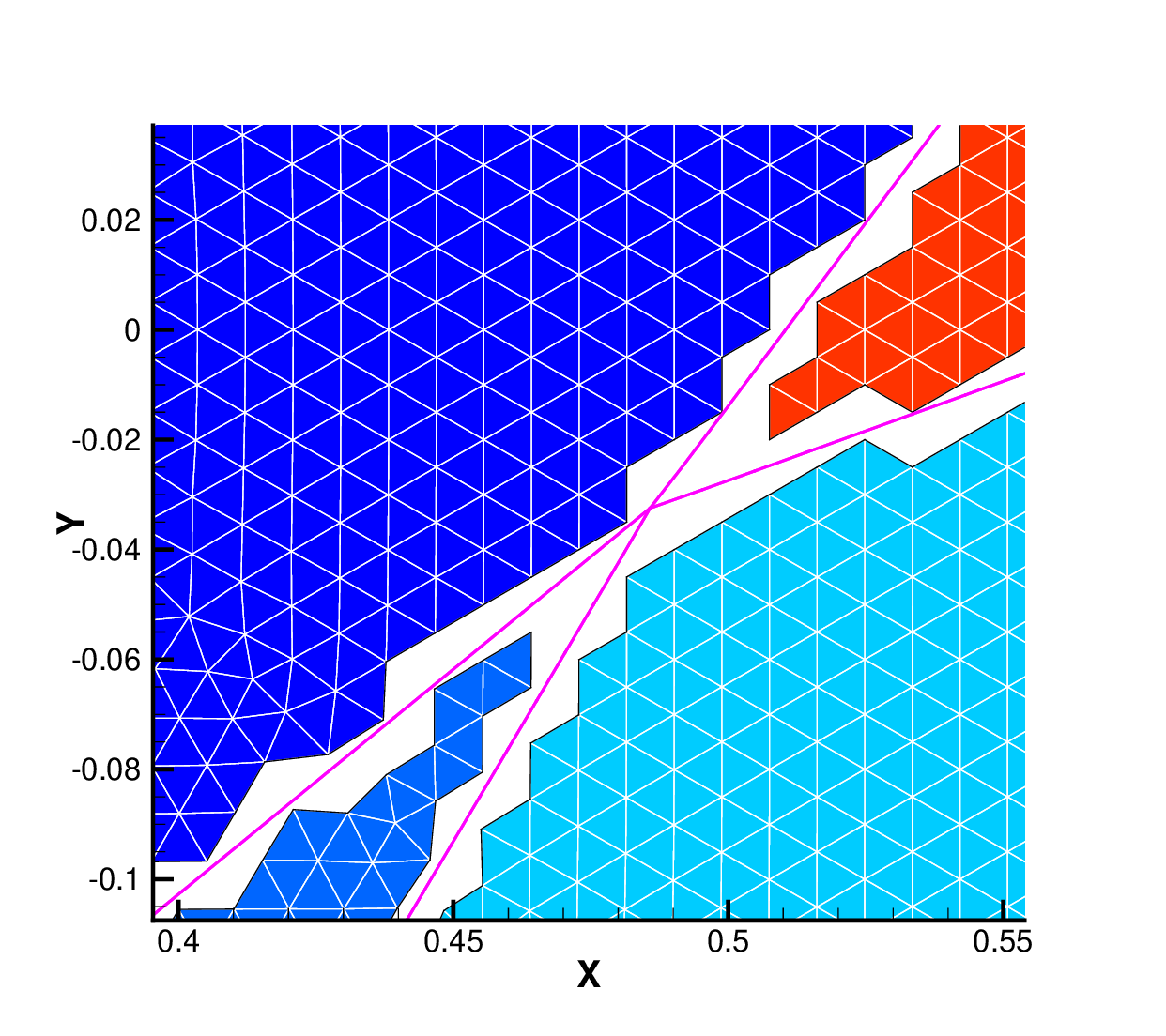}\label{pic34f}}
\caption{Interaction between two shocks of the same family: numerical solutions (comparisons in terms of $S=p\rho^{-\gamma}$).}\label{pic34}
\end{figure*}

\revMR{As mentioned in \S~\ref{shocks calc}, \revAB{the} interaction points need to be modelled 
explicitly on a case by case basis. 
If only some of the discontinuities meeting at the interaction point are tracked, the only possibility 
is to  solve independent
algebraic problems arising from the jump conditions for each discontinuity, and use some approximate formula for 
QP. 
For the type of interaction considered here,} 
\revAB{we have used the following relation to compute the velocity of the interaction point QP}:
\begin{equation}\label{eq:QP_speed}
\boldsymbol{\omega}_{QP}\,=\,\boldsymbol{\omega}_{IS1}\,+\,\left(\boldsymbol{\omega}_{IS2}
\cdot \boldsymbol{\tau_1}\right) \, \boldsymbol{\tau_1}
\end{equation}
where $\boldsymbol{\tau_1}$ is the vector tangential to IS1 (see Fig.\ref{pic33}a). 
\revAB{Equation~(\ref{eq:QP_speed}) is similar, but not identical, to the
formula proposed in~\cite{chang2019adaptive} for the same purpose}.
%

\revMR{We compare in Fig.~\ref{pic34} three sets of results: \revAB{{\it i})} a fully captured solution
\revAB{(top row of Fig.~\ref{pic34})};
\revAB{{\it ii})} a hybrid solution in which CD is captured, 
and all the shocks (IS1, IS2 and RS1)  are tracked
\revAB{(central row of Fig.~\ref{pic34})};
\revAB{{\it iii})} a fully tracked result
\revAB{(bottom row of Fig.~\ref{pic34})} \revRP{where also CD is fitted}. First of all, this result}
shows that it is possible to track several discontinuities and capture others, 
\revMR{all involved in the same  interaction. Compared to the fully captured solution, this hybrid simulation
with captured CD provides a much cleaner entropy field, as seen comparing the top and central
rows of frames in Fig.~\ref{pic34}. 
However,  we can also see a small anomaly due to the capturing of the CD,
which is still   observable in part of the flow downstream of the reflected shock RS1:
see the close up of Fig~\ref{pic34d}.}
These spurious disturbances
\revMR{are removed  by }
adding \revAB{the contact discontinuity} to the set of discontinuities to 
be tracked by the algorithm, \revMR{as visible on the bottom row of results of Fig.~\ref{pic34}.} 
%
\revMR{Note that despite its apparent simplicity, the results obtained are surprisingly good as this is a difficult simulation. Indeed, on a mesh as  coarse 
as the one used here, the elements crossed by the discontinuity generate a relatively large cavity around the interaction point. This is clearly visible for example 
in} Fig~\ref{pic34f}.  
\revMR{The extrapolation strategy proposed allows to handle this delicate geometrical situation.}
%

%
        
\subsection{Shock-wall interaction: Mach reflection}\label{MR}
\begin{figure*}
\centering
 \includegraphics[width=0.6\textwidth]{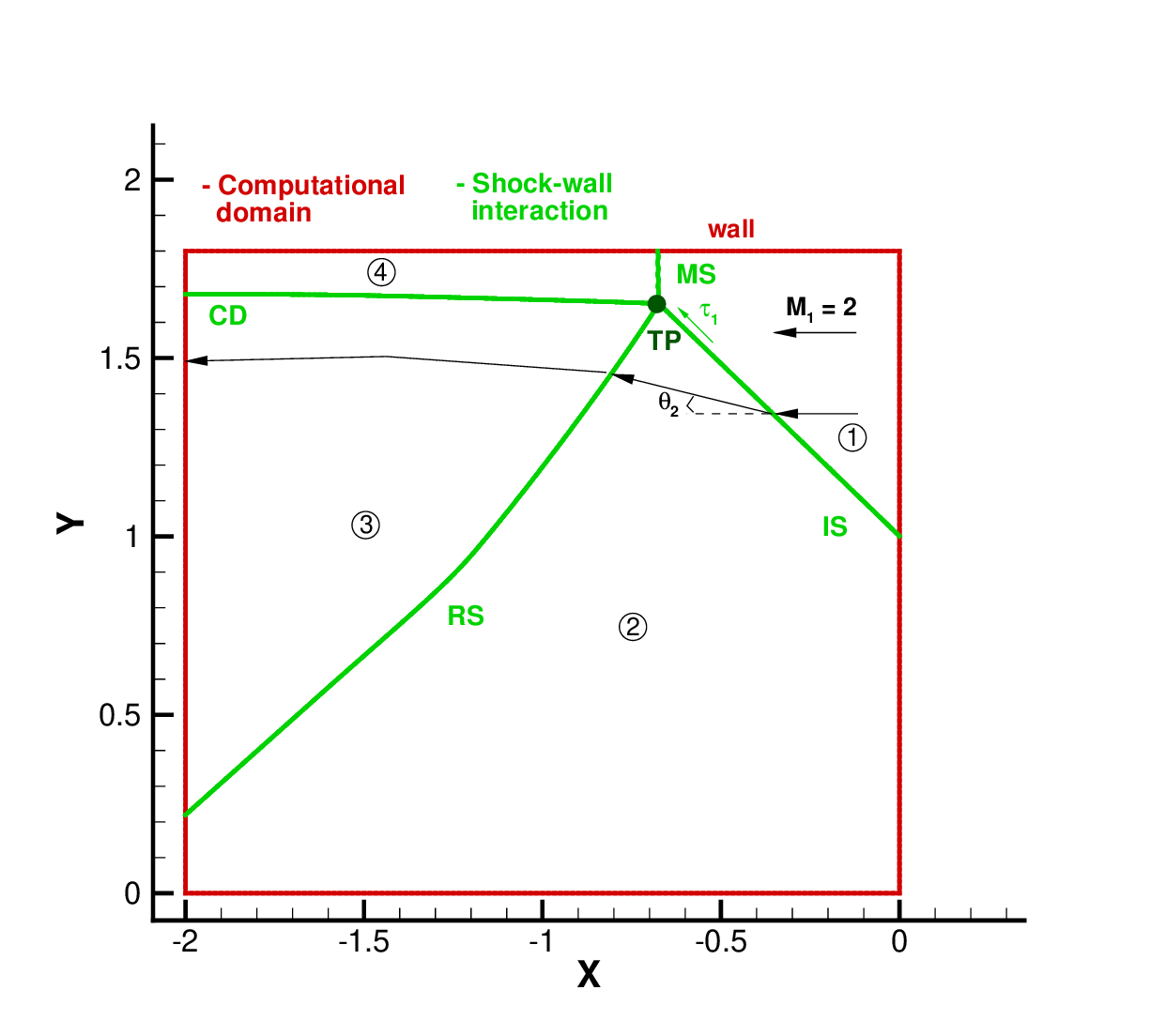}
\caption{Mach reflection: sketch of the flow.}\label{pic35}
\end{figure*}
\revAB{An oblique shock that impinges on a straight wall gives rise to either a regular or a Mach
reflection, depending on the combination of free-stream Mach number, $M_1$, and flow deflection 
angle, $\theta_2$, that the free-stream flow undergoes 
while passing through the oblique shock (hereafter also referred to as the incident shock).
Whenever $\theta_2$ is
larger than the maximum deflection that the supersonic stream behind the 
incident shock can sustain, a Mach reflection takes place.
As sketched in Fig.~\ref{pic35},
a Mach reflection consists in a fairly complex three-shocks system
(the incident shock, IS, the reflected shock, RS, and the
Mach stem, MS)} interacting in the triple point, TP, from which a contact discontinuity, CD, also arises. 
\revMR{Note that differently from the \revAB{test-case addressed in \S~\ref{SS2}}, 
here not all the regions \revAB{surrounding TP} involve \revRP{uniform} flow. 
Moreover, \revAB{CD} presents a slight angle, which makes it in general not mesh-aligned}.\\ 

\begin{figure*}
\centering
\subfloat[$SC$ solution]	      {\includegraphics[width=0.45\textwidth]{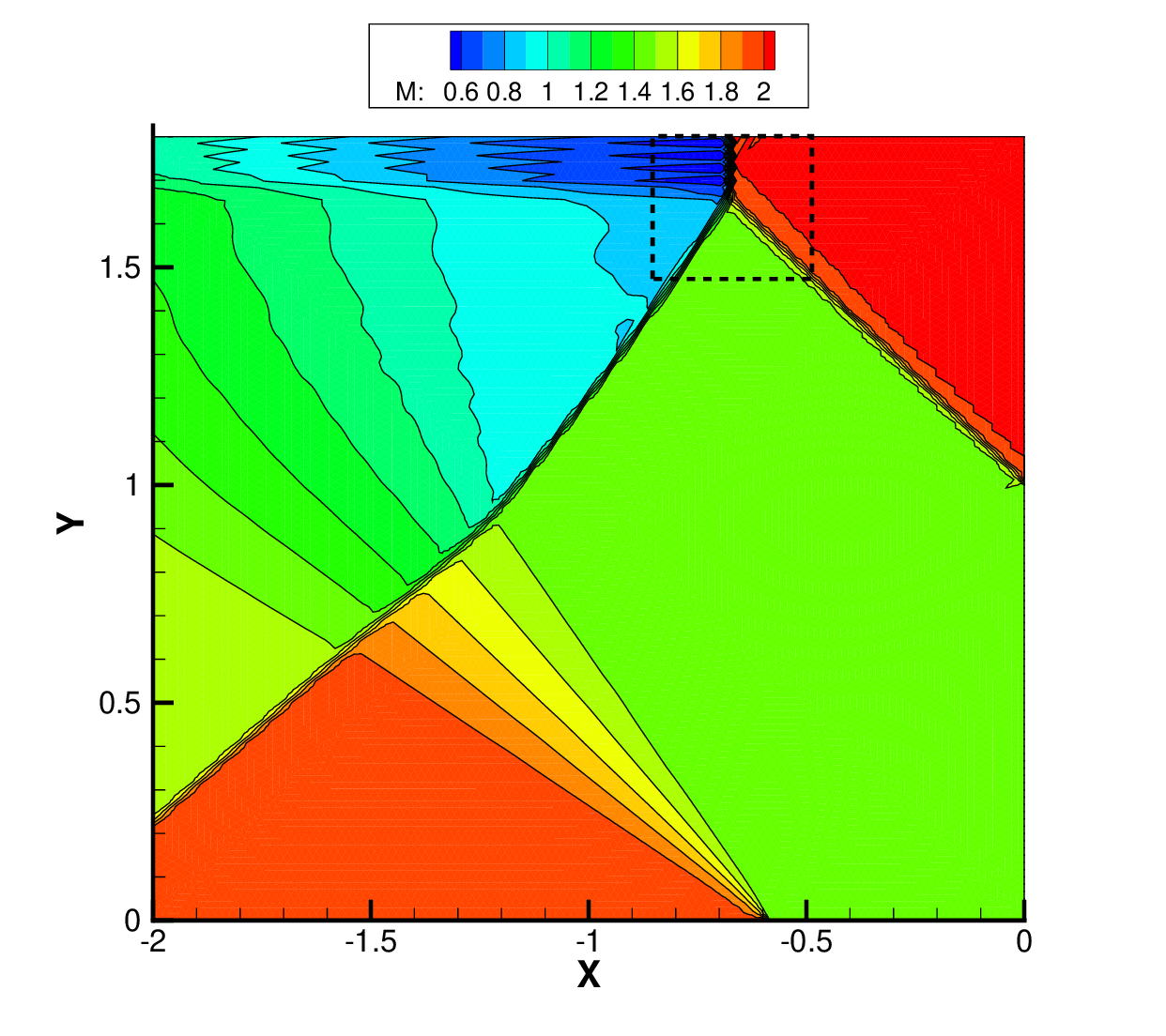}}
\subfloat[close-up around the TP]     {\includegraphics[width=0.40\textwidth]{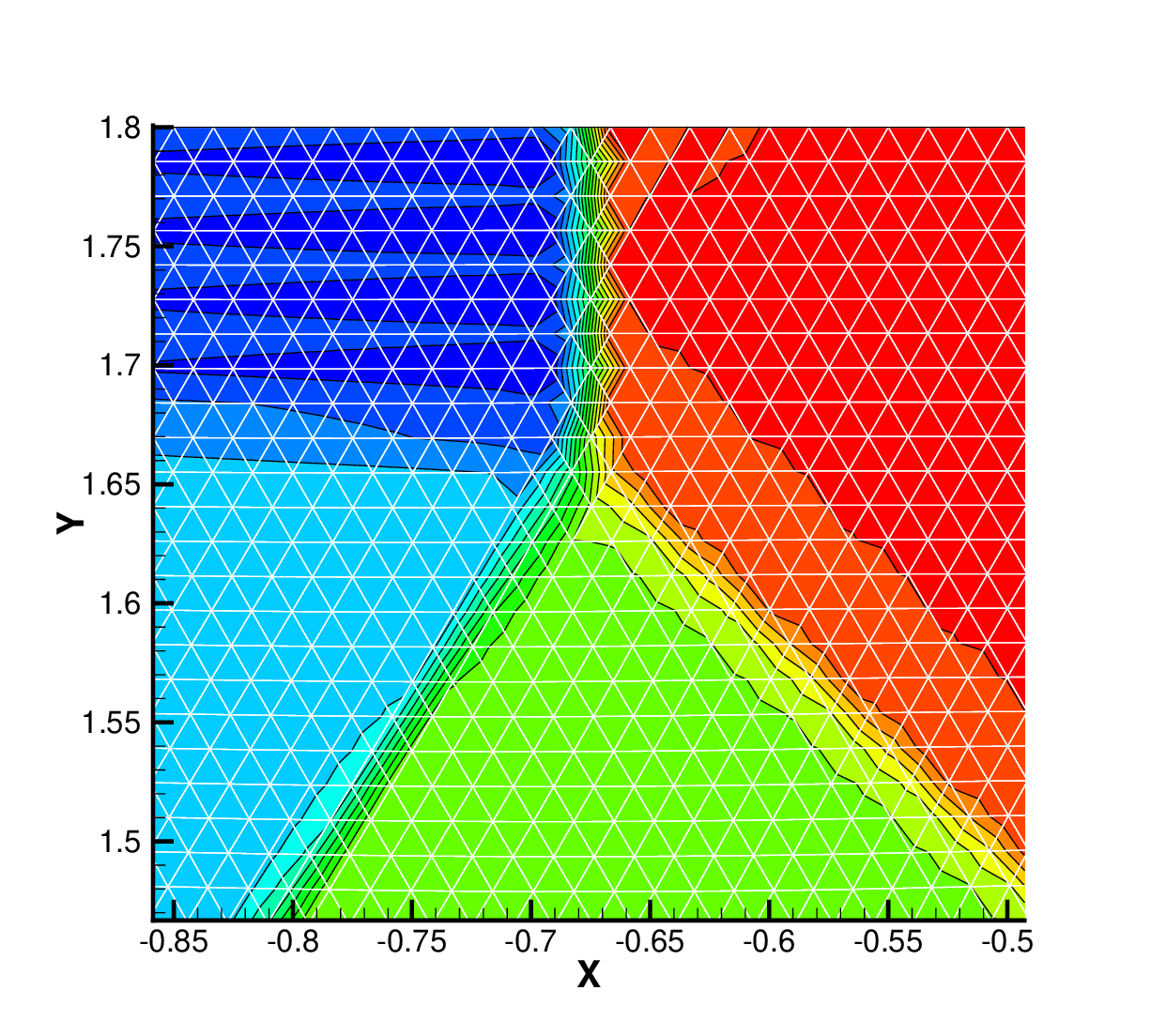}}\qquad
	\subfloat[Hybrid $eDIT_{GG}$: \revAB{configuration 2}] {\includegraphics[width=0.45\textwidth]{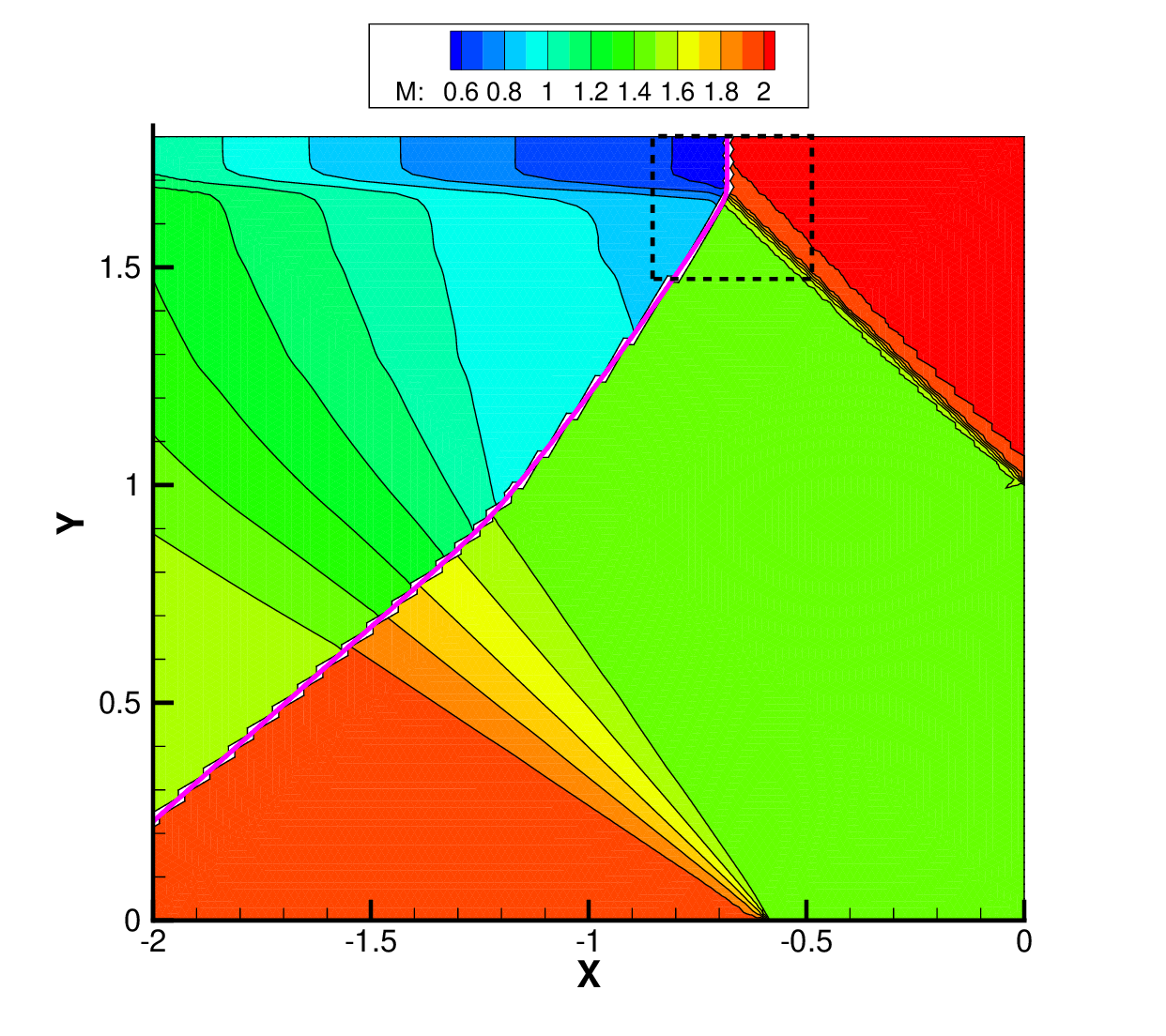}}
\subfloat[close-up around the TP]     {\includegraphics[width=0.40\textwidth]{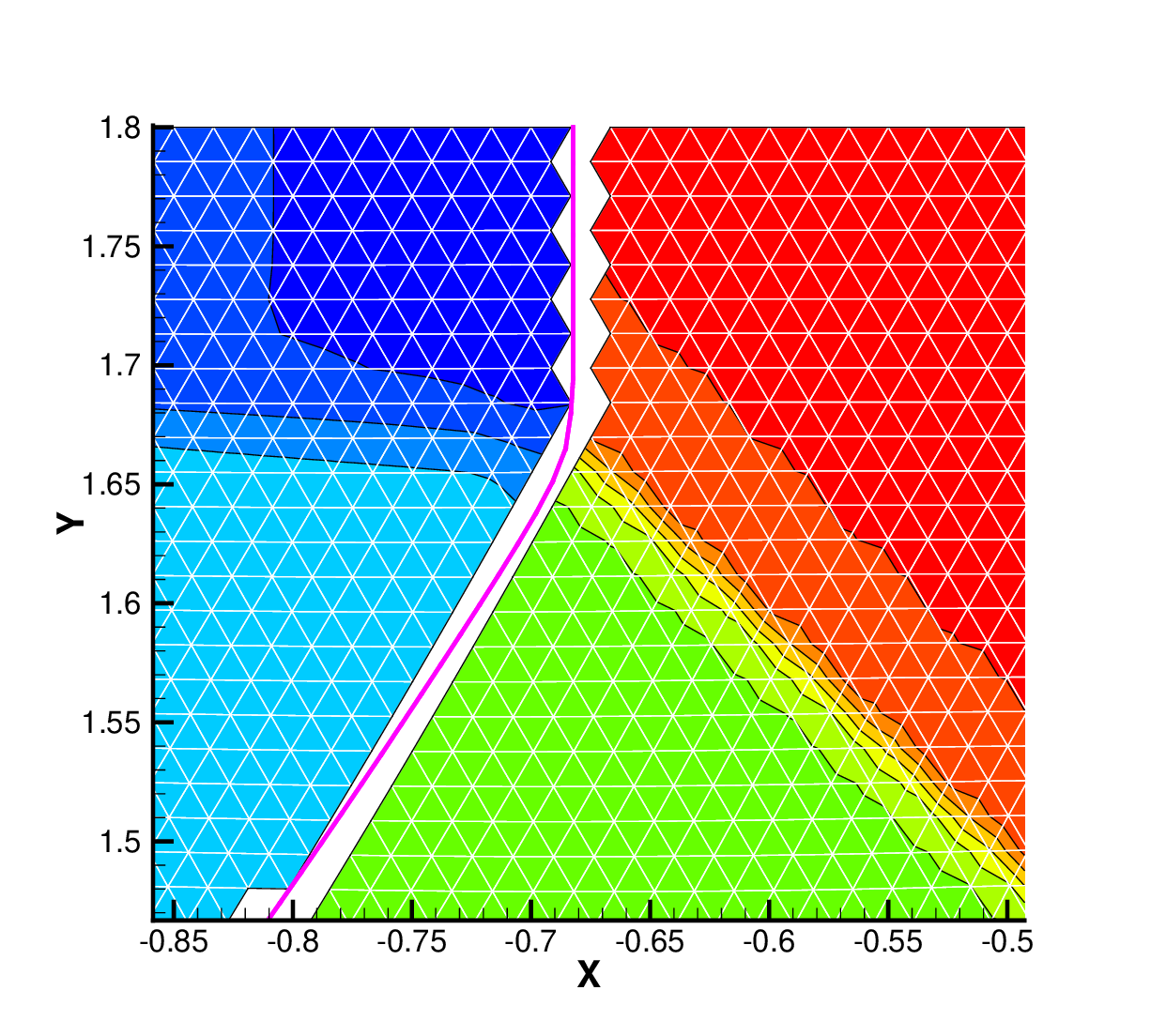}}\qquad
	\subfloat[Hybrid $eDIT_{GG}$: \revAB{configuration 3}] {\includegraphics[width=0.45\textwidth]{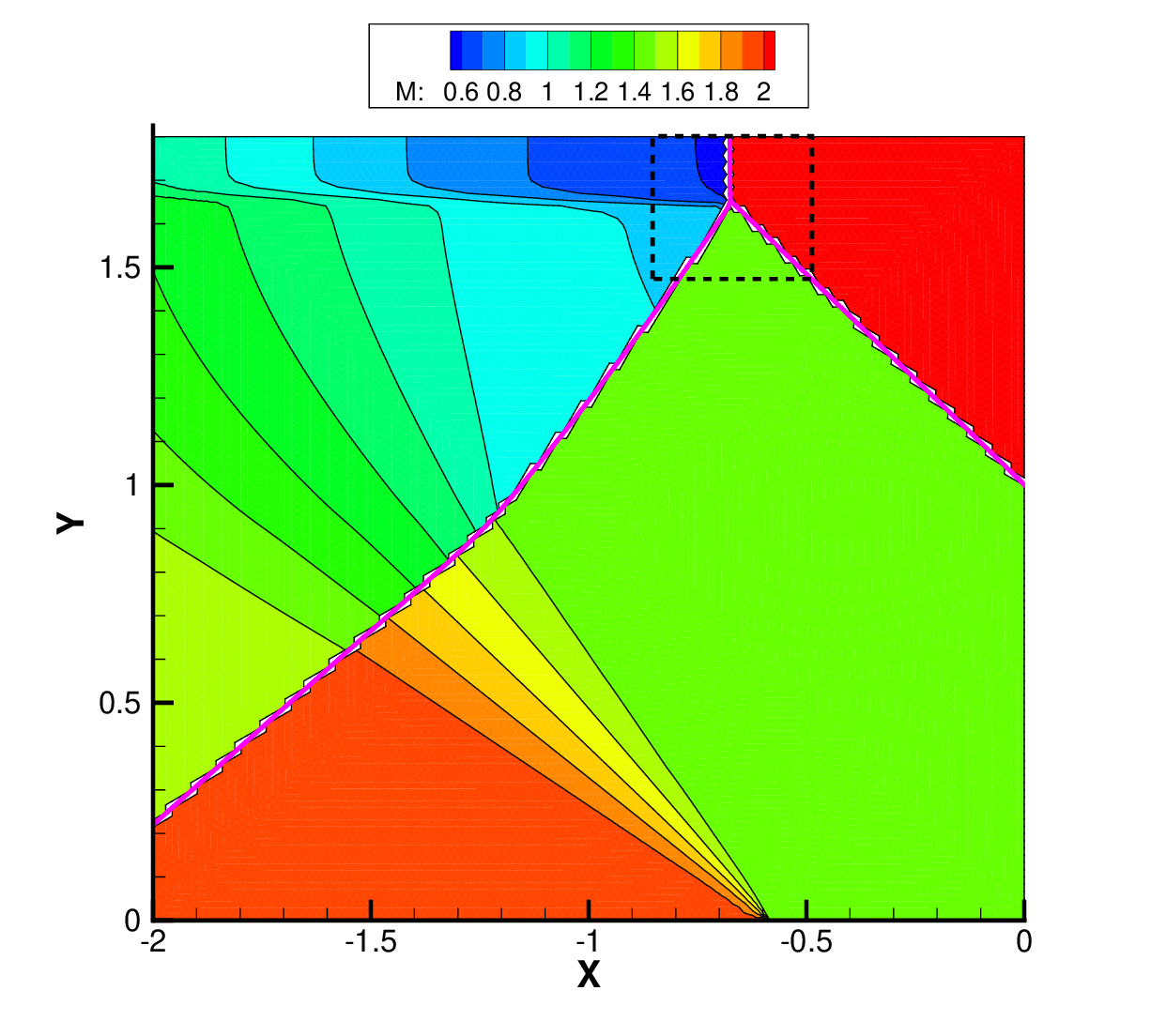}}
\subfloat[close-up around the TP]     {\includegraphics[width=0.40\textwidth]{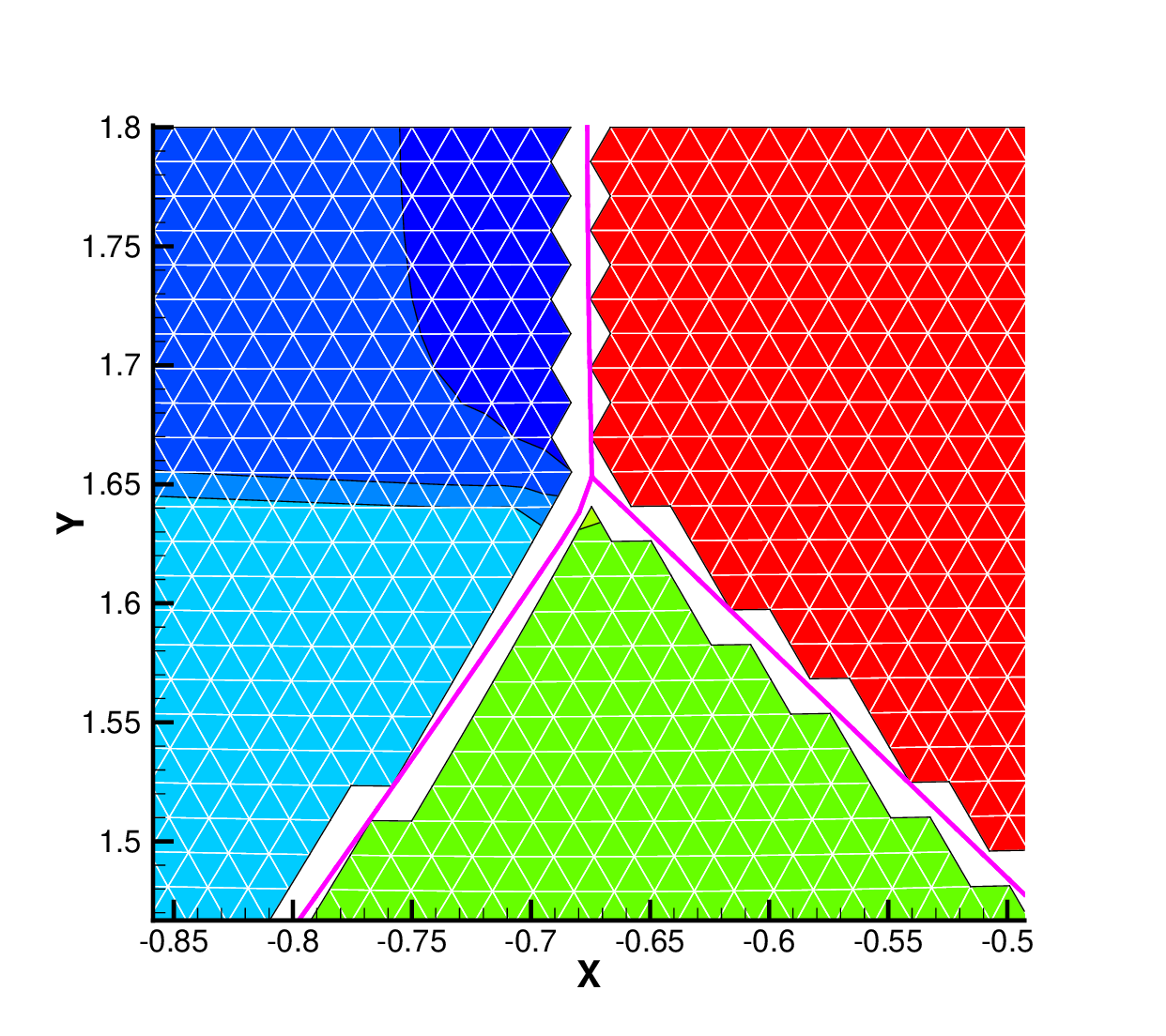}}
\caption{Mach reflection: \revAB{the Mach number iso-lines} of the numerical solutions.}\label{pic36}
\end{figure*}
%

\revAB{For the chosen setting, $M_1 = 2$ and $\theta_2 = 14^\circ$},
\revMR{we consider four different types of simulations: 
\begin{enumerate}
\item fully  captured;
\item  an extrapolated shock-tracking setting in which the Mach stem (MS) and the reflected shock (RS) are tracked as a single discontinuity;
\item a hybrid setting in which  only the shocks are tracked, but CD is captured; 
\item full-fledged discontinuity\revAB{-tracking} for all shocks plus the CD.
\end{enumerate}
Note that no modelling of the triple point TP is necessary in configuration 2 which involves a unique shock-mesh. 
In case 3 we need to provide an explicit model for the triple-point.
Not having enough equations to compute \revAB{the} TP, \revAB{because the CD is captured, rather than being tracked}, 
we need to resort to a heuristic  approach. As in the case \revAB{of \S~\ref{SS2}}, 
a nonlinear algebraic problem is solved independently for each discontinuity,
and for \revAB{the} TP \revRP{velocity}  we use the following simplified relation}:
\begin{equation}\label{eq:TP_speed}
\boldsymbol{\omega}_{TP}\,=\,\boldsymbol{\omega}_{IS}\,+\,\left(
\boldsymbol{\omega}_{MS} \cdot \boldsymbol{\tau_1}\right) \, \boldsymbol{\tau_1}
\end{equation}
where $\boldsymbol{\tau_1}$ is the unit vector tangential to the IS (see Fig.\ref{pic35}a).\\ 
\revAB{Equation~(\ref{eq:TP_speed}) is similar, but not identical, to
the formula used in~\cite{chang2019adaptive} for the same purpose}.

\begin{figure*}
\centering
\subfloat[Initial solution]    {\includegraphics[width=0.35\textwidth,trim={1cm 0.5cm 2.4cm 0cm},clip]{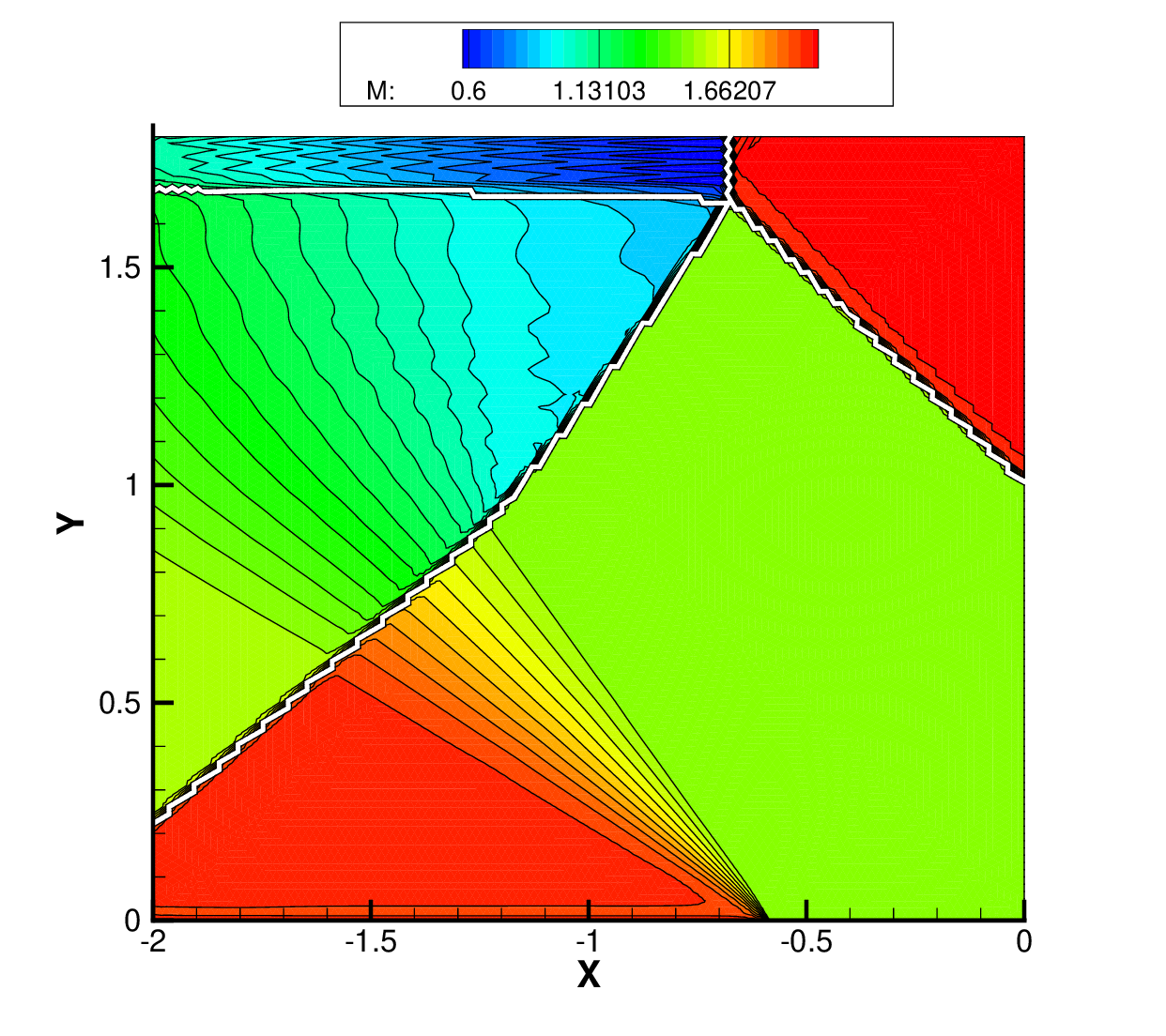}}
\subfloat[After 200 iterations]{\includegraphics[width=0.35\textwidth,trim={1cm 0.5cm 2.4cm 0cm},clip]{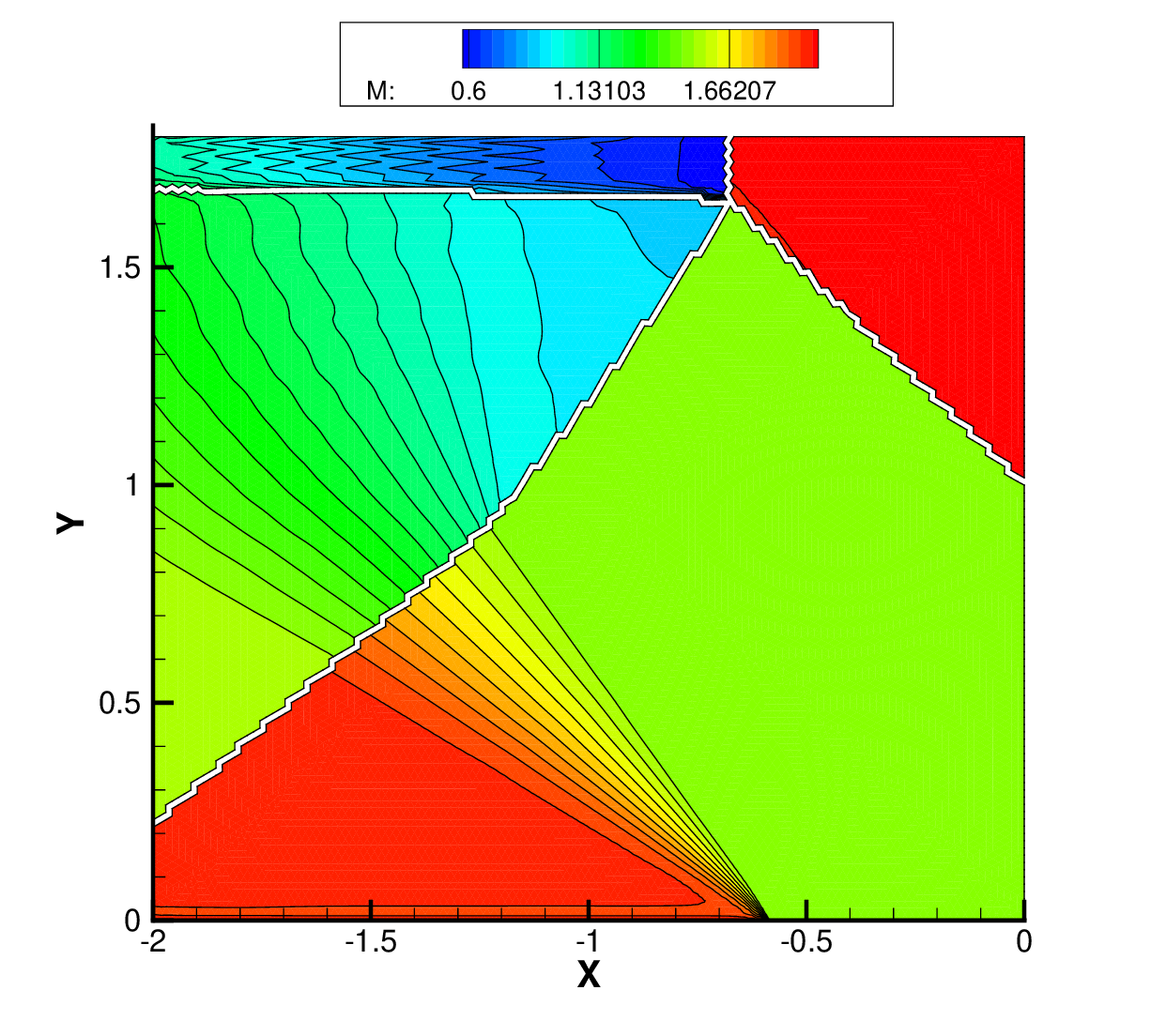}}
\qquad
\subfloat[After 400 iterations]{\includegraphics[width=0.35\textwidth,trim={1cm 0.5cm 2.4cm 0cm},clip]{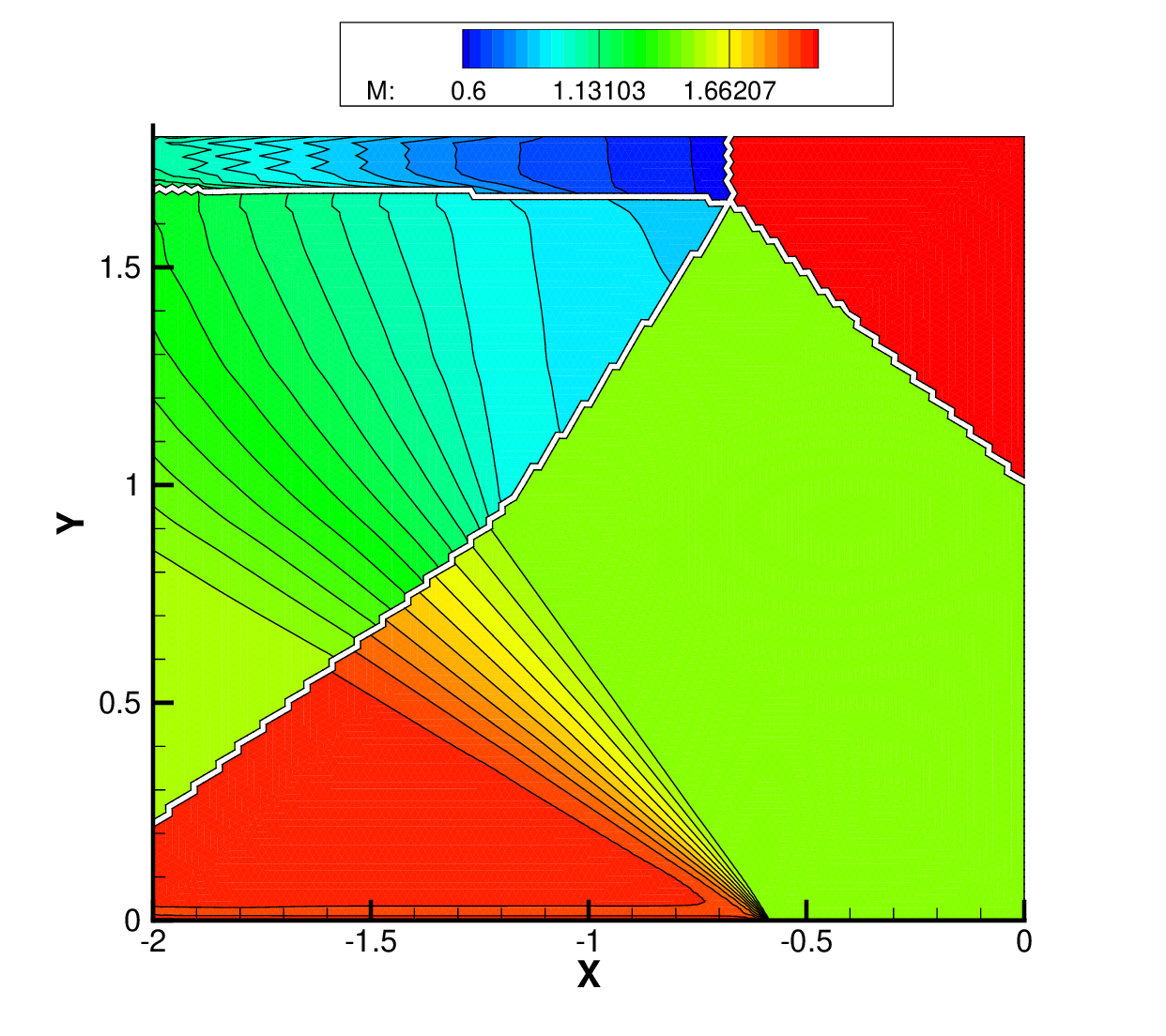}}
\subfloat[Converged solution]  {\includegraphics[width=0.35\textwidth,trim={1cm 0.5cm 2.4cm 0cm},clip]{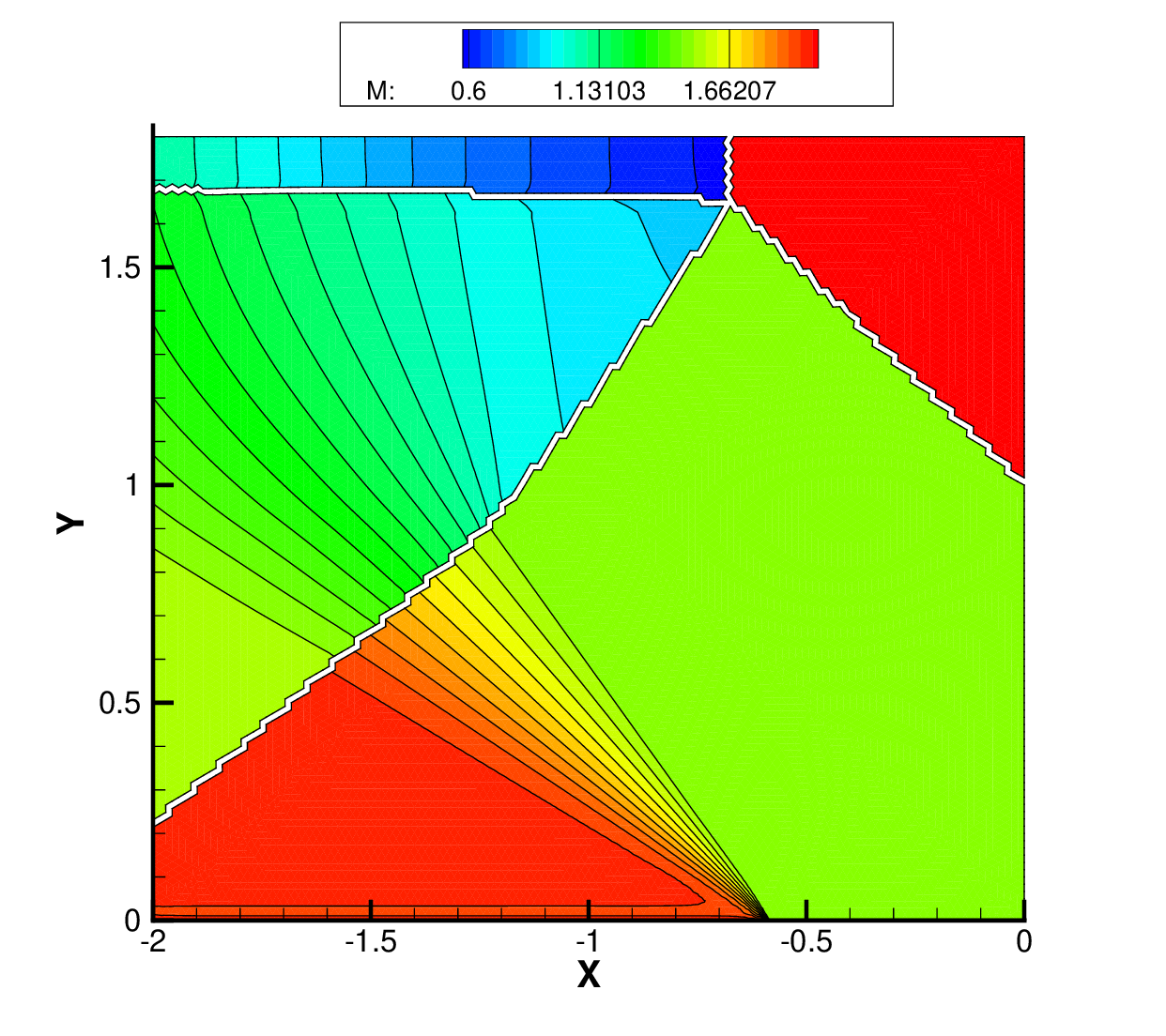}}
\caption{Mach reflection: pseudo-time convergence of the full-fledged $eDIT_{GG}$ simulation (Mach number iso-contour lines)}\label{pic38}
\end{figure*}

Finally, in configuration  4,
\revAB{all four states surrounding the TP
and its unkown velocity, $\boldsymbol{\omega}_{TP}$,
can be coupled into a single non-linear system of algebraic equations,
whose solution provides updated values downstream of the RS and MS
as well as $\boldsymbol{\omega}_{TP}$.
Full details are given in~\cite{Paciorri5}}.\\
%
%
\begin{figure*}
\centering
\subfloat[Full-fledged $eDIT_{GG}$ solution] 	{\includegraphics[width=0.45\textwidth]{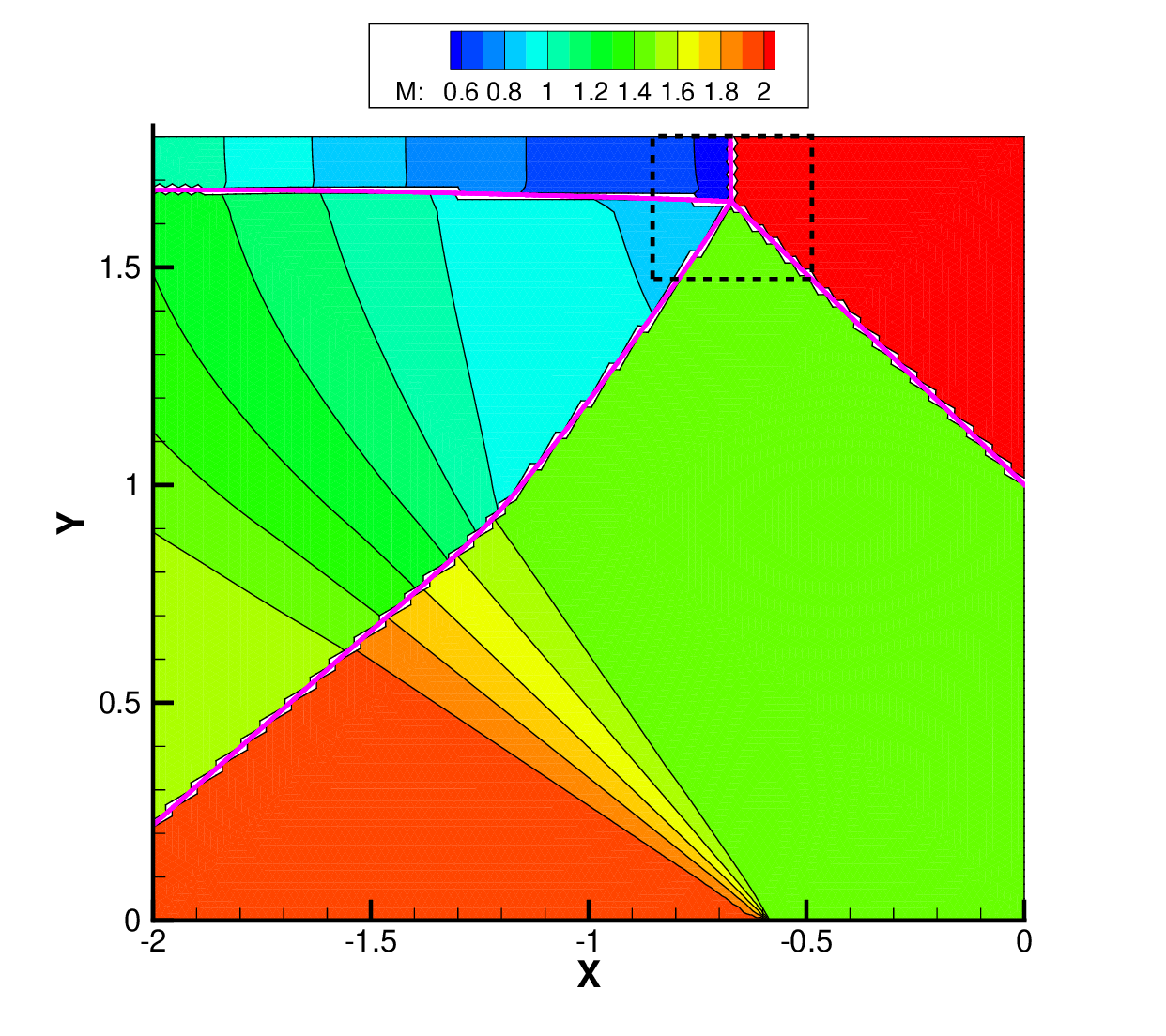}}
	\subfloat[close-up around \revAB{the} TP]  			{\includegraphics[width=0.40\textwidth]{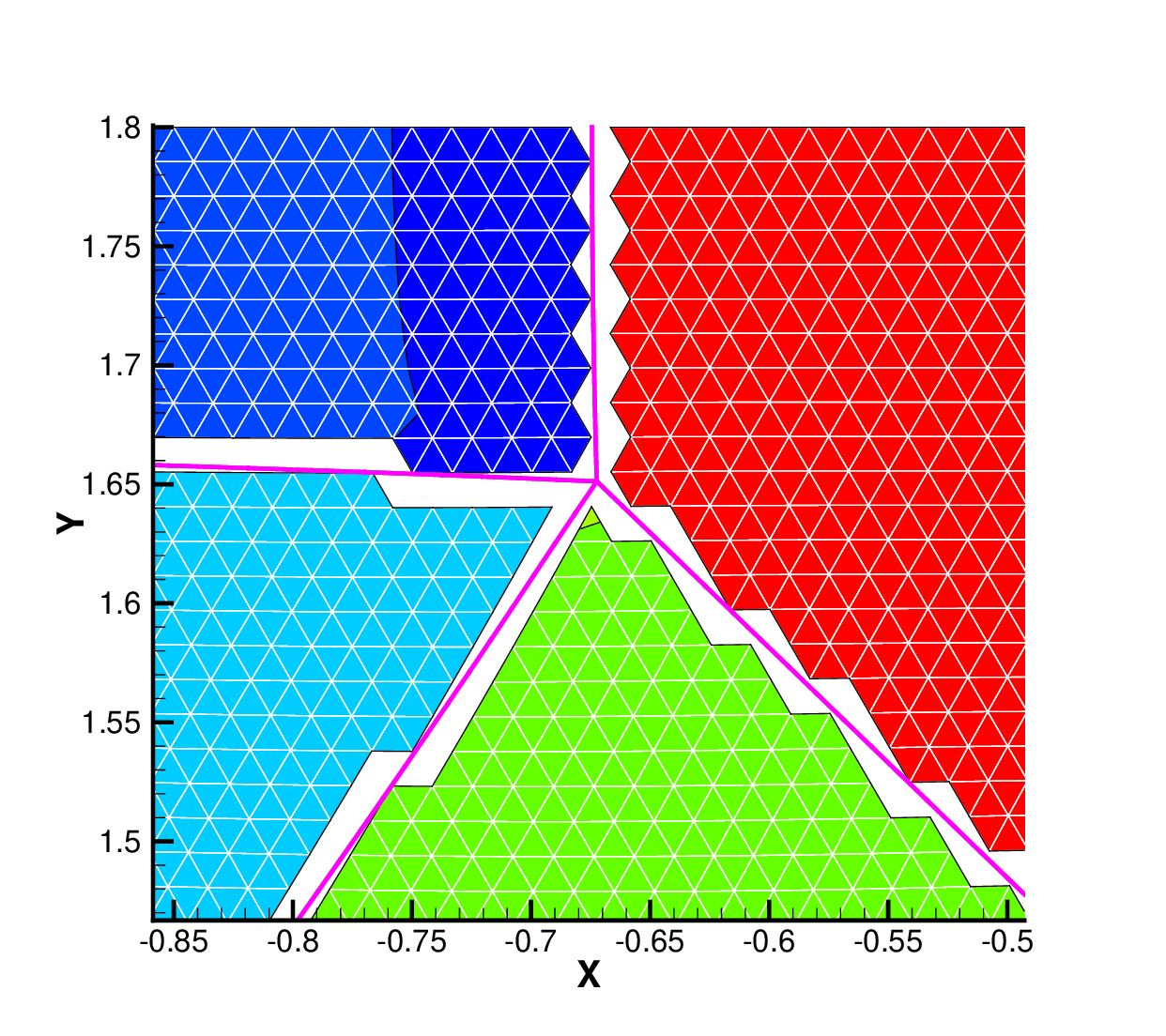}}
	\caption{Mach reflection: \revAB{the Mach number iso-lines} of the full-fledged (configuration 4) numerical solutions.}\label{pic37}
\end{figure*}
\revMR{The results are arranged in four rows in Figures~\ref{pic36} and~\ref{pic37},
showing on the left an overview of the Mach contours, and on the right a zoom of the \revAB{TP}.}
\revAB{Compared to the SC calculation (top row in Fig.~\ref{pic36}),}
\revMR{tracking MS and RS (second row in Fig.~\ref{pic36}) already  provides an enormous improvement in the quality of the solution. 
The \revAB{CD} remains however poorly captured on this coarse mesh.
Adding the \revAB{IS} to the tracking set (third row in Fig.~\ref{pic36}) leads to a cleaner flow, 
which is visible in the nicer straight contour lines downstream of RS. 
However, \revAB{CD} is still poorly resolved. Finally, in full-fledged $eDIT$ mode of Fig.~\ref{pic37} we are able to 
resolve the \revAB{CD} in a single row of elements. 
A small kink in the contour lines in the non-constant region around CD is still visible. We assume these to be due to 
perturbations arising in correspondence of the steps in the surrogate \revAB{CD} and propagating downstream. 
A qualitative view of the error cleaning and pseudo-time convergence of the full-fledged simulation is displayed in Fig.~\ref{pic38}.}  

%
%

\subsection{Supersonic channel flow}\label{channel}
\begin{figure*}
\centering
	\subfloat[Sketch of the computational domain \revAB{and flow-pattern}.]{\includegraphics[width=0.65\textwidth,trim={2cm 1cm 2cm 4cm},clip]{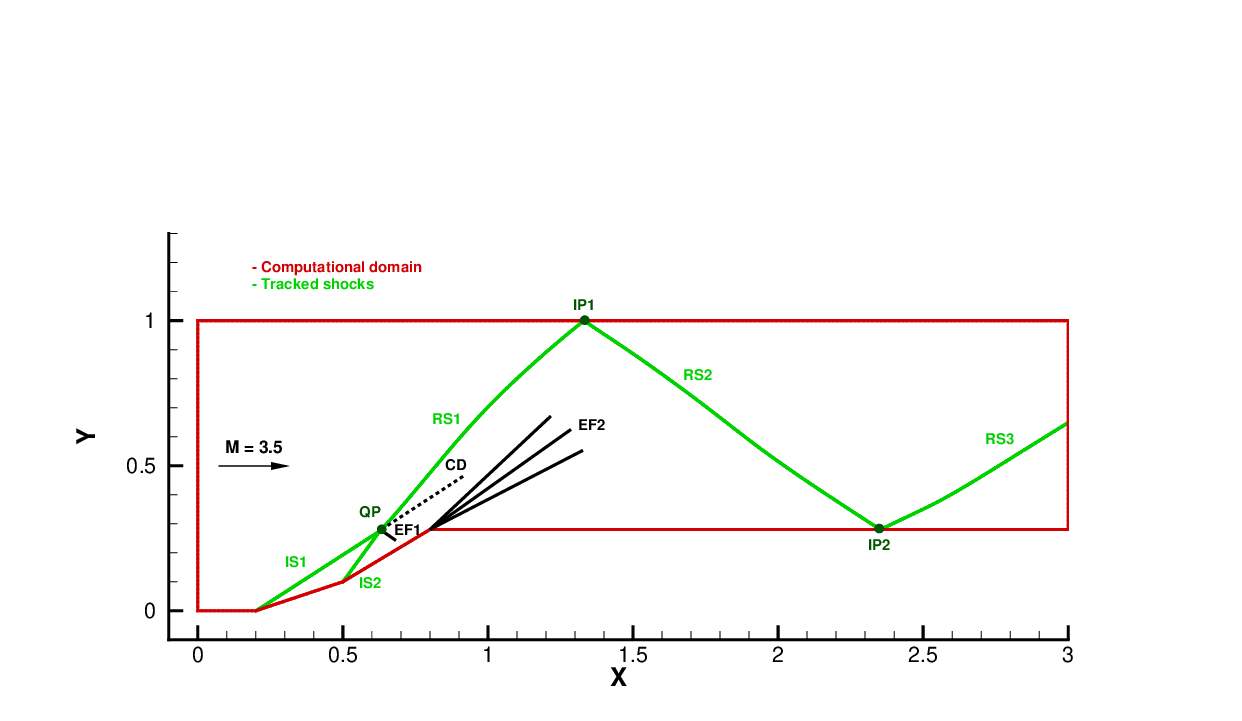}\label{pic40a}}
\subfloat[\revAB{Detail of the background triangulation.}]{\includegraphics[width=0.35\textwidth,trim={0cm 0cm 0cm 0cm},clip]{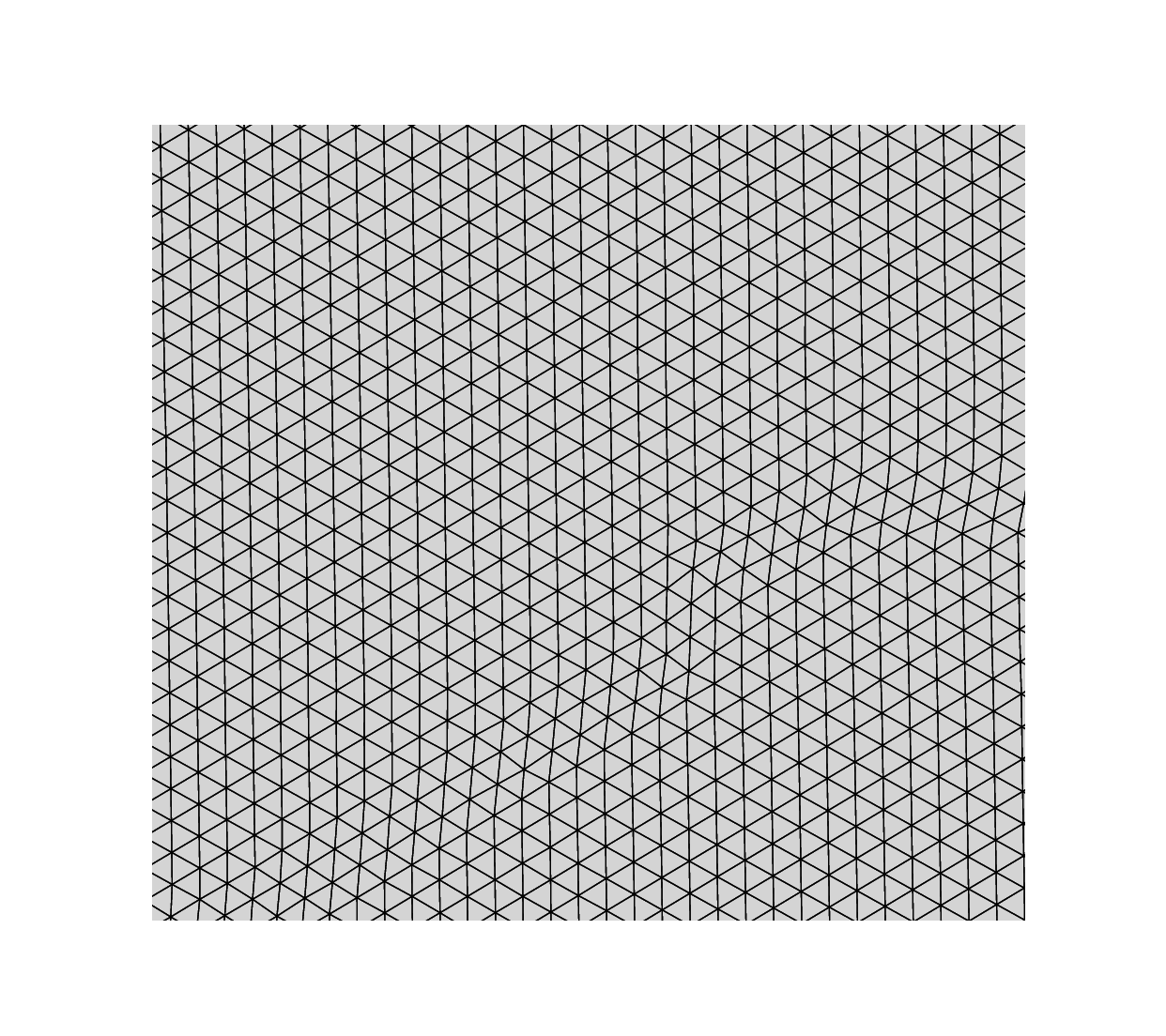}\label{pic40b}}
\caption{Supersonic channel flow.}\label{pic40}
\end{figure*}
The last test-case considered consists in \revAB{the} supersonic, $M_{\infty} = 3.5$, flow in a planar channel, 
whose variable-area geometry is shown in Fig.~\ref{pic40a} and reported in~\cite{zou2021}: 
the two constant-area \revAB{portions of the duct} are joined through a double ramp. 
\revAB{The flow pattern is as follows}: two oblique, straight shocks of the same family, labeled
IS1 and IS2, originate at the two \revAB{convex} corners of the \revAB{ramp and} their interaction 
gives rise to a wave configuration 
already explored in~\S~\ref{SS2} where IS1, IS2, a new shock RS1, a \revAB{slip-stream} CD and 
an expansion fan EF1 meet \revAB{at the interaction point QP}.
A second, stronger expansion fan EF2 takes place at the \revAB{concave} corner of the lower wall and interacts 
with RS1, which bends
before being reflected from the upper wall. The reflected shock RS2 is again reflected by the lower wall and leaves the duct as shock RS3.   

\begin{figure}
\centering
\includegraphics[width=0.35\textwidth,trim={5cm 3cm 5cm 5cm},clip]{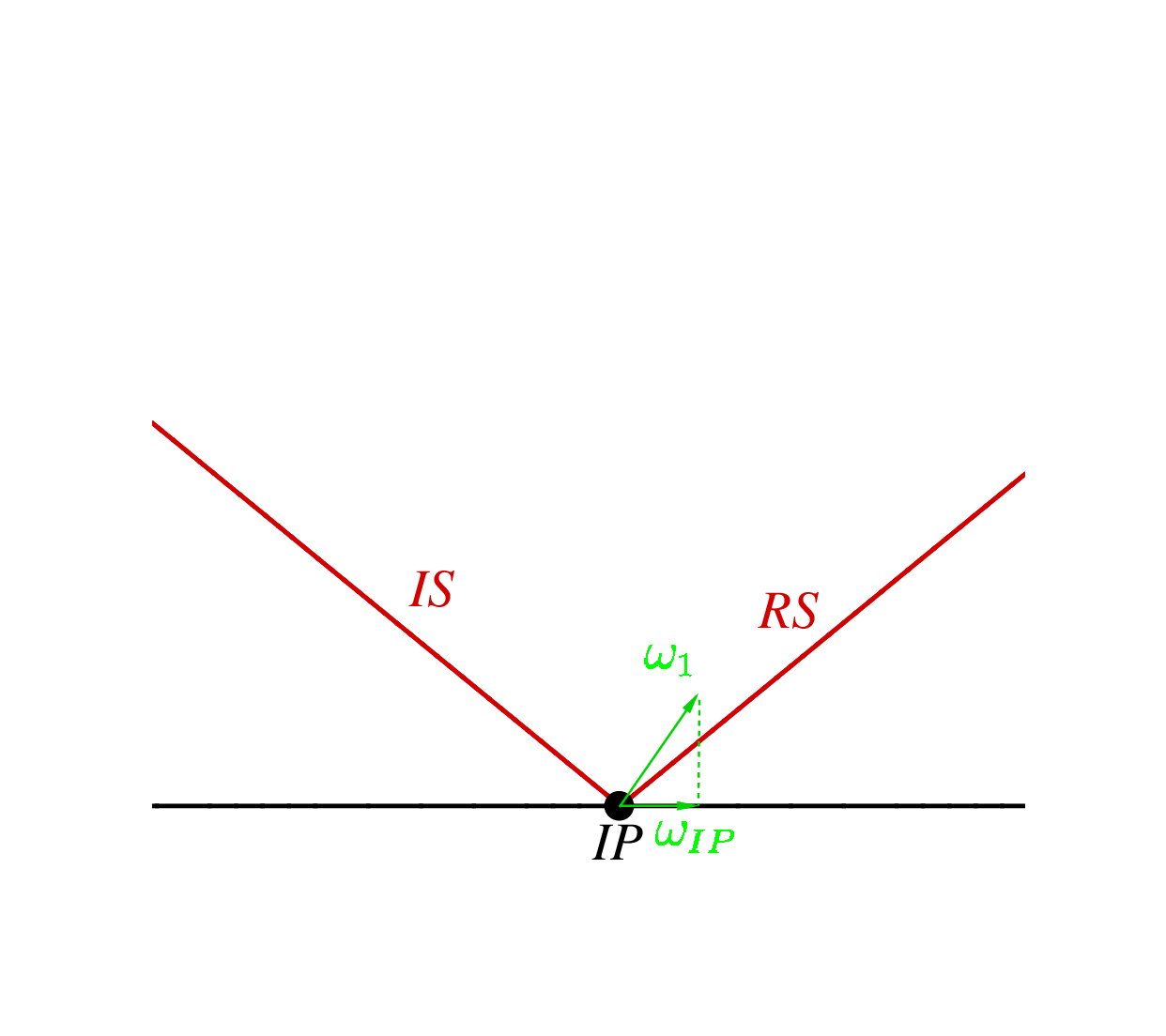}
	\caption{Supersonic channel flow: \revAB{regular reflection at a solid wall and calculation of
	of the displacement velocity of the point, IP, where the IS and RS meet}.}\label{RRspeed}
\end{figure}

\revAB{The goal of the present test-case consists in assessing that the}
different features that have been 
implemented within the $eDIT$ algorithm, \revAB{already described in the previous sections},
work well also when \revAB{combinedly used to simulate
a fairly complex shock-pattern, such as the one illustrated in Fig.~\ref{pic40a}}.
\revAB{To improve simulation fidelity, we track as many discontinuities as we can, i.e.\ the}
three shock waves IS1, IS2 and RS1, which meet at QP, and the two regular reflections
made up by RS1, RS2 and RS3 taking place on the upper and lower walls. \revAB{Only the CD has been captured}.\\\
When simulating this \revAB{testcase} with the $eDIT$ algorithm, two \revAB{different} models 
have \revAB{been} employed: 
\revAB{the one already described in~\S~\ref{SS2} for}
QP, \revAB{where the interaction between shocks of the same family takes place,}
and a different one for points IP1 and IP2, \revAB{where a regular reflection takes place on the upper, 
resp.\ lower wall}.
\revAB{As shown in Fig.~\ref{RRspeed}, the velocity of points IP1 and IP2 is set}
equal to the component tangential to the wall of the corresponding incident shock velocity.

\revAB{Simulations have been performed using both $SC$ and $eDIT_{GG}$}.
The unstructured grid used in the $SC$ calculation and as
the background triangulation in the $eDIT_{GG}$ calculation,
which
is made up of 12,417 grid-points and 24,310 \revAB{nearly equilateral} triangles,
has been generated using the frontal mesh generator of the 
\texttt{\ttfamily gmsh} software~\cite{geuzaine2009gmsh}.
\revAB{A detail of the mesh is shown in Fig.~\ref{pic40b}.}

\begin{figure*}
\centering
\subfloat[$SC$ solution]  {\includegraphics[width=0.60\textwidth,trim={1cm 0cm 2cm 4.5cm},clip]{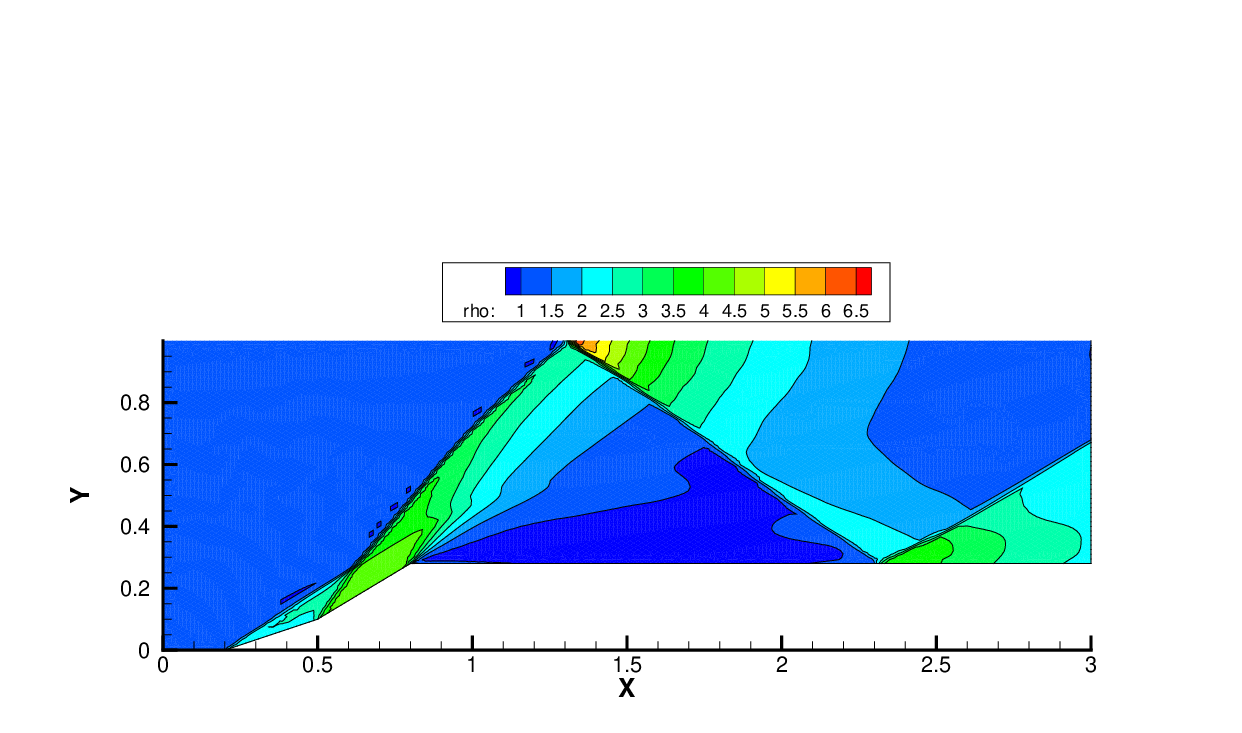}\label{pic41a}}
\qquad
\subfloat[$eDIT_{GG}$ solution]{\includegraphics[width=0.60\textwidth,trim={1cm 0cm 2cm 4.5cm},clip]{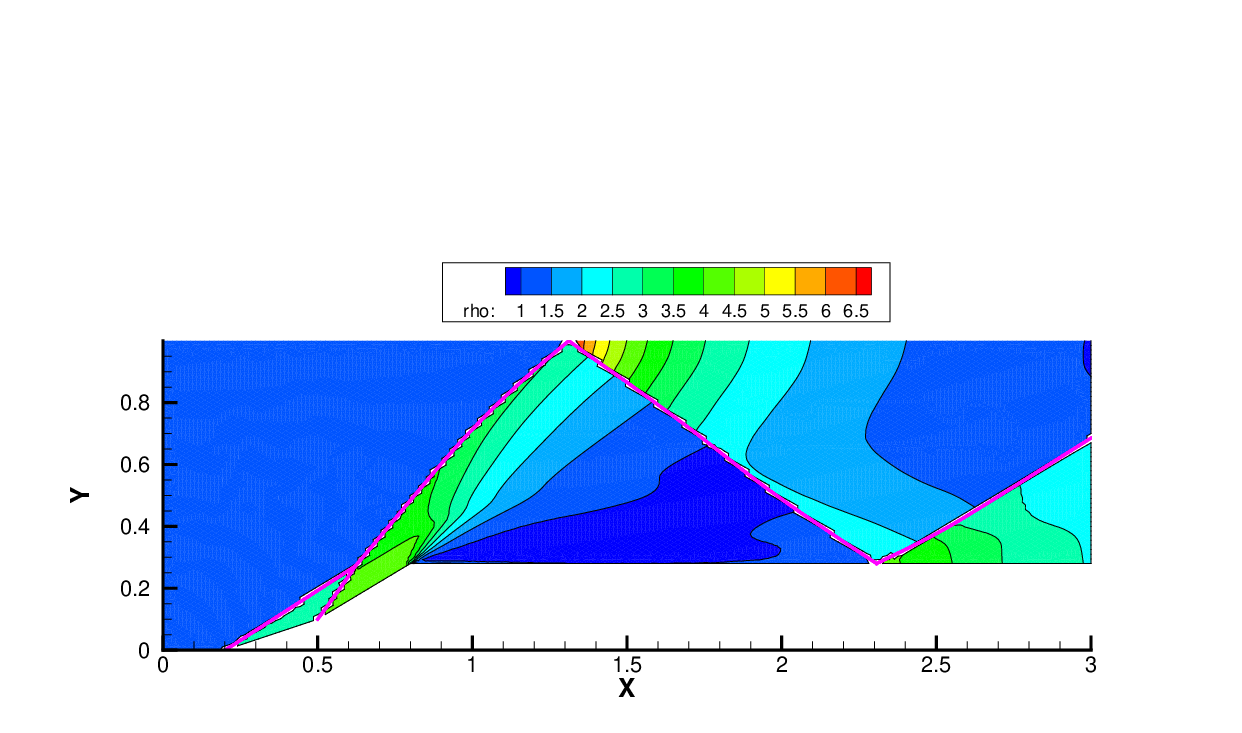}\label{pic41b}}
\caption{Supersonic channel flow: \revAB{dimensionless density $\rho/\rho_{\infty}$} iso-contours.}\label{pic41}
\end{figure*}

\begin{figure*}
\centering
\subfloat[$SC$ solution]  {\includegraphics[width=0.60\textwidth,trim={1cm 0cm 2cm 4.5cm},clip]{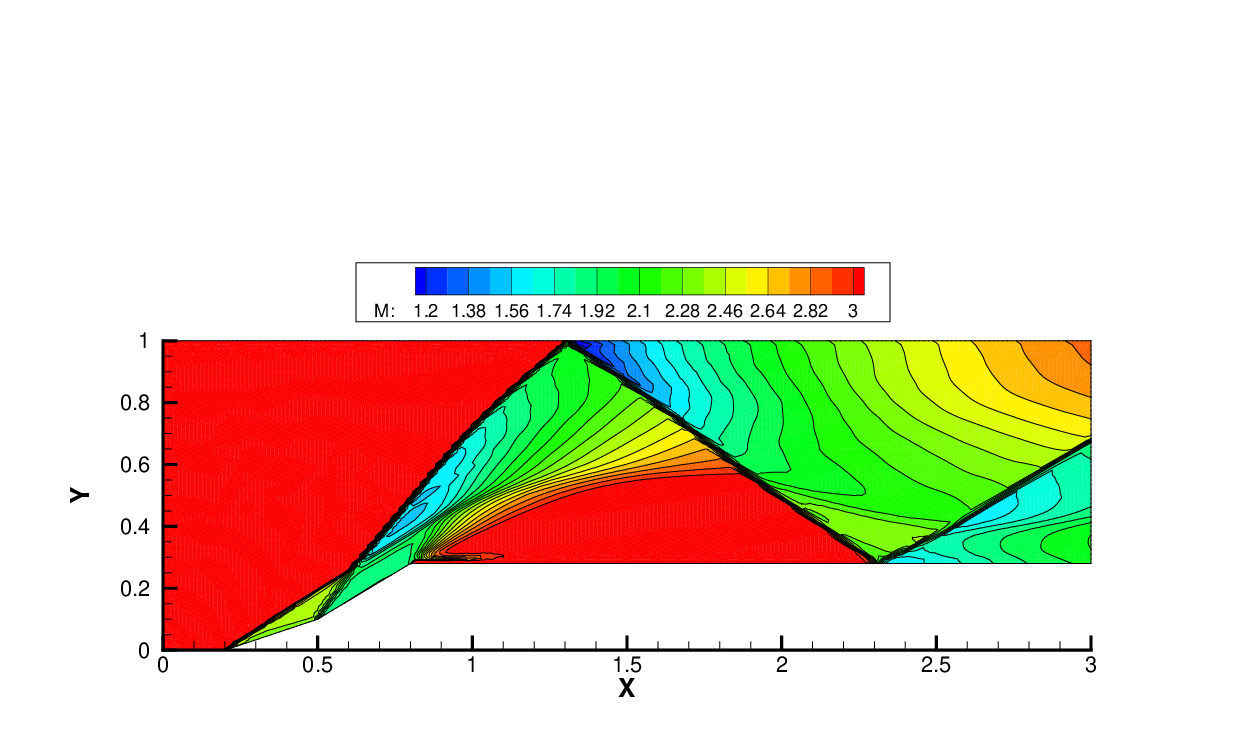}\label{pic42a}}
\qquad
\subfloat[$eDIT_{GG}$ solution]{\includegraphics[width=0.60\textwidth,trim={1cm 0cm 2cm 4.5cm},clip]{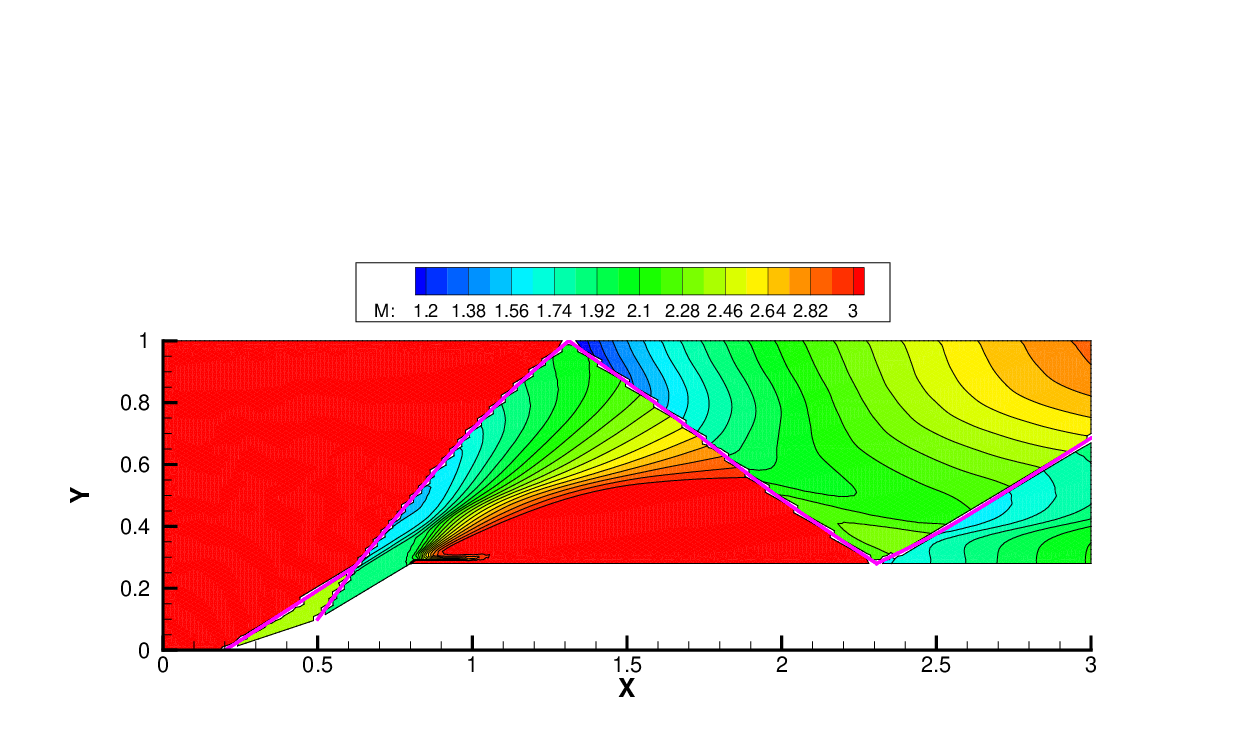}\label{pic42b}}
\caption{Supersonic channel flow: \revAB{Mach iso-contours}.}\label{pic42}
\end{figure*}

\begin{figure*}
\centering
\subfloat[Upper wall.]{\includegraphics[width=0.40\textwidth,trim={1cm 0cm 1cm 0cm},clip]{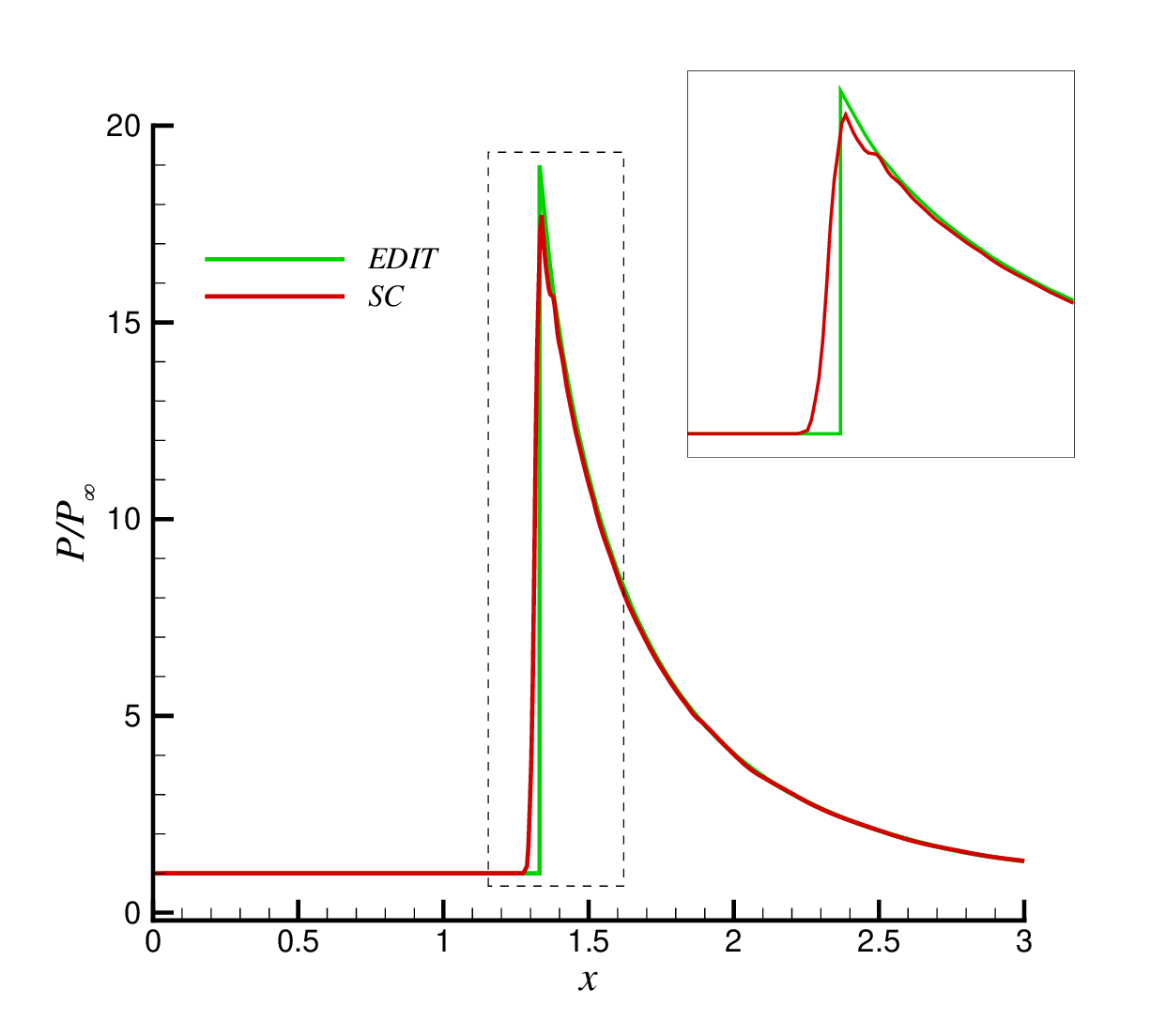}}
\subfloat[Lower wall.]{\includegraphics[width=0.40\textwidth,trim={1cm 0cm 1cm 0cm},clip]{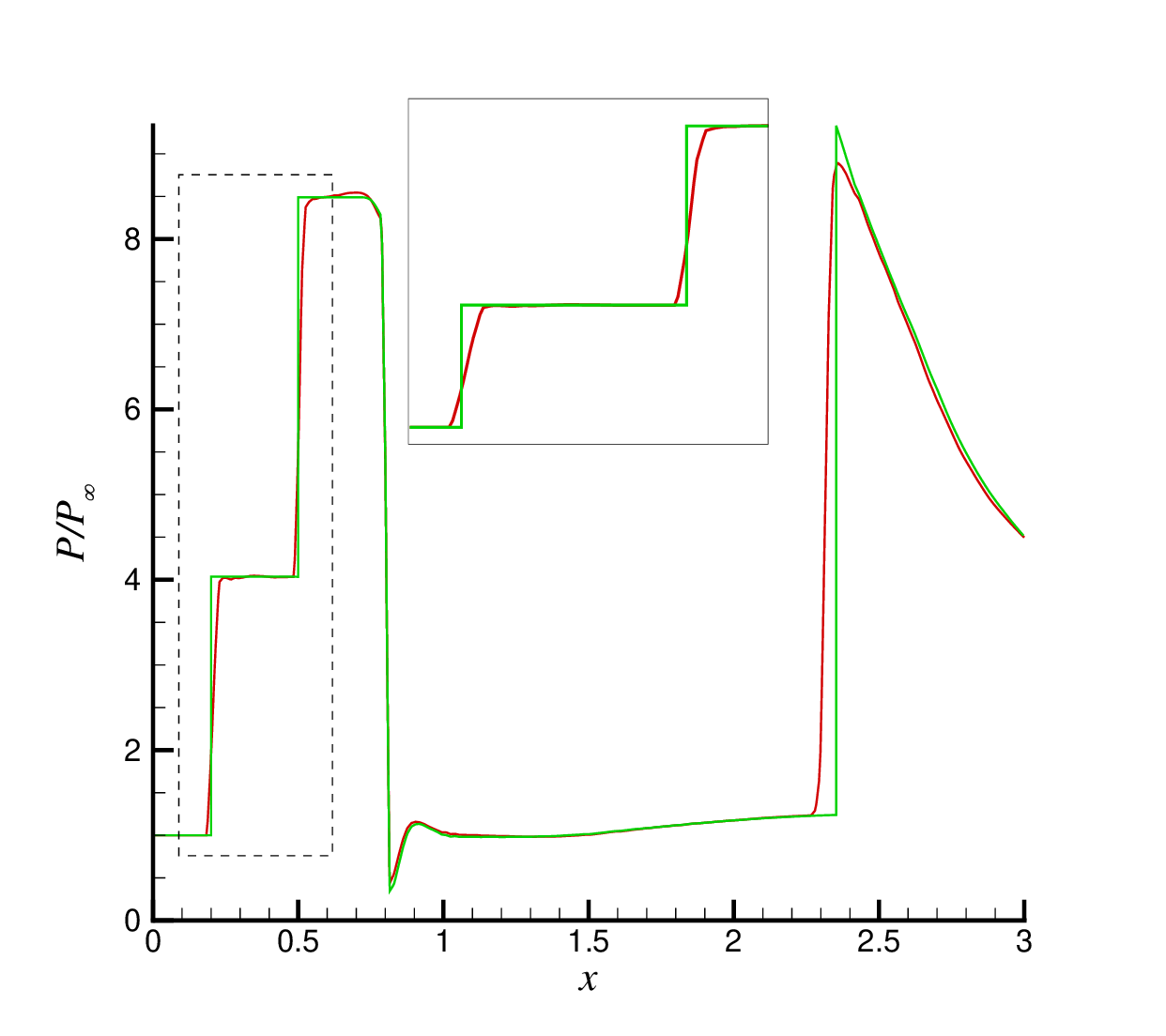}}
\caption{Supersonic channel flow: dimensionless pressure distribution along the walls of the channel computed by means of $SC$ and $eDIT$.}\label{pic43}
\end{figure*}

\revAB{The comparison between the two sets of calculations is reported in Figs.~\ref{pic41} and~\ref{pic42},
were density and Mach number iso-contours are respectively displayed}.
\revAB{It can be seen that}
the capture of the discontinuities, \revAB{see Figs.~\ref{pic41a} and~\ref{pic42a},}
\revAB{gives rise to} spurious disturbances (\revAB{in particular downstream of RS2})
which are not present in the $eDIT_{GG}$ calculation (Figs.~\ref{pic41b} and~\ref{pic42b})
pointing out a globally \revAB{cleaner} flow-field \revAB{featuring} smoother iso-contours. 
The oscillations that occur downstream of RS2 are more clearly visible in Fig.~\ref{pic43},
\revAB{which shows
the dimensionless pressure profiles along the upper and lower walls.
It is evident that by resorting to a shock-tracking approach, not only we are able to 
get rid of the major oscillation introduced
by $SC$ and exactly represent the discontinuities as having zero-thickness,
but we also get a 
better prediction of wall pressure peaks, which are smoothed out in the $SC$ calculation}.

        \section{Conclusion}\label{conclusion}

A novel technique recently proposed \revRP{in~\cite{EST2020}} to simulate flows with shock waves 
has been further improved to \revAB{make it capable of dealing}
with different kinds of discontinuities (\revAB{both} shock waves and contact discontinuities) 
\revRP{as well as shock-shock and shock-wall interactions, thus
opening the possibility to compute  complex flows. Moreover, since  the algorithm 
can \revAB{be} run in hybrid mode, it is  also possible to study complicated flows  by tracking some 
of the discontinuities and leaving others to be captured.}
The \revRP{proposed} technique provides genuinely second-order-accurate results even
for flows featuring very strong shocks, 
without the complexity of the re-meshing/adaptation phase of previous fitting approaches. 
\revRP{For this reason,  the present algorithm is substantially independent from the data structure 
of the gasdynamic solver and, therefore, 
it can be applied with small modifications to  both cell-centered and vertex-centered solvers 
on both unstructured and structured grids.
Although many improvements of the fitting/tracking technique have been made in the present paper, 
further developments are in progress, \revAB{including its}
coupling with a structured-grid, finite volume solver~\cite{Assonitis2021}.}

        \appendix
\section{Nodal gradient reconstruction}\label{appendixA}
\revAB{As explained in \S~\ref{setting}, in the 
present version of the $eDIT$ algorithm, the 
computational domain is tesselated into triangles and
bi-linear shape functions with unknowns stored
at the vertices of the triangulation that provide
the functional representation of the 
dependent variable, which is Roe's parameter vector $\mathbf{Z} = \sqrt{\rho} \left(1,H,u,v\right)^t$.
Using the aforementioned setting,
the gradient of $\mathbf{Z}$ is constant within each
triangular cell $T$ and can be readily computed as follows:
\begin{equation}\label{eq:cell_gradient}
\nabla  \mathbf{Z}_T\,=\,\frac{\sum_{k=1}^3(\mathbf{Z}_k\,\mathbf{n}_k)}{2\,A_T}
\end{equation}
where $\mathbf{n}_k$ denotes the normal to the face opposite vertex $k$, scaled by its length
and pointing inside triangle $T$ and $A_T$ denotes the area of triangle $T$.

As explained in~\S~\ref{transfer1}, however, in order to perform the data trasfer
between the surrogate boundaries and the discontinuity-front,
it is also necessary to compute the gradient of $\mathbf{Z}$ within every grid-point of the surrogate boundaries;
this is accomplished using either the Green-Gauss (GG) or Zienkiewicz-Zhu (ZZ) approaches to be detailed hereafter}.
%
%
\paragraph{Green-Gauss reconstruction}\label{GG}
\revAB{The GG reconstruction consists in camputing the nodal gradient in grid-point $i$
as the area-weighted average of the cell-wise constant gradients of all triangles
that surround grid-point $i$:
\begin{equation}\label{eq:nabla}
\nabla \mathbf{Z}_i\,=\,\sum_T \, \theta_T \nabla\mathbf{Z}_T
\end{equation}
}
where: 
\begin{equation}
	\label{eq:thetat}
\theta_T \,=\, \frac{A_T}{\sum_{T} A_T} \;\;\; \text{such that} \;\;\; \sum_{T} \theta_T \,=\, 1
\end{equation}
\revAB{The summations in Eqs.~(\ref{eq:nabla}) and~(\ref{eq:thetat}) range over all 
triangles surrounding the grid-point $i$}.
%
%
\paragraph{Zienkiewicz-Zhu reconstruction}\label{ZZ}
\revAB{Similarly to the GG reconstruction, all steps involved in the ZZ gradient-reconstruction 
have to be repeated within all grid-points of
the surrogate boundaries. In order to alleviate the notation, however, we shall hereafter drop the index that refers
to the grid-point and use subscript $k$ to refer to one of the four components
of the parameter vector and subscript $j$ 
to refer the cartesian coordinates, i.e.\
$\left(x_1,x_2\right) =\left(x,y\right)$}.
The starting point in the ZZ reconstruction consists in using a linear polynomial expansion
of each component of the gradient, i.e.
%
\begin{equation}
	\left(\frac{\partial Z_k}{\partial x_j}\right)\,=\,\mathbf{P}\cdot\mathbf{a}_{k,j}
\end{equation}
where:
\begin{equation}\label{polynom1}
	\mathbf{P}\,=\,[\,1\,,\,x_1\,,\,x_2\,] \;\;\; \text{and} \;\;\; \mathbf{a}_{k,j}\,=\,[\,a_0\,,\,a_1\,,\,a_2\,]_{k,j}
\end{equation}
Following the Zienkiewicz-Zhu (ZZ) patch recovery procedure, 
the unknown parameters $\mathbf{a}_{k,j}$ \revAB{are computed via}
a least-squares fit \revAB{which amounts to}
minimizing the following test function:
\begin{equation}\label{LSfunction}
	F(\mathbf{a}_{k,j})\,=\,\sum_{T} \theta_T \left[\left(\frac{\partial Z_k}{\partial x_j}\right)_{T}\,-\,\mathbf{P}(x_1^T,x_2^T)\cdot \mathbf{a}_{k,j}\right]^2
\end{equation}
\revAB{The summation in Eq.~(\ref{LSfunction}) ranges over all triangles surrounding
the given grid-point, the gradient in Eq.~(\ref{LSfunction}) is computed 
element-wise according to Eq.~(\ref{eq:cell_gradient}) 
and vector $\mathbf{P}$ is computed in the barycentric coordinates of triangle $T: \left(x_1^T,x_2^T\right)$}.
Given that three unknowns, $\mathbf{a}_{k,j}$, must be computed for each component of $\mathbf{Z}$
and each cartesian coordinate, a stencil of at least three triangles is required.
The minimization problem defined by Eq.~\eqref{LSfunction} can be solved in matrix form:
\begin{equation}\label{eq:ZZ}
\mathbf{a}_{k,j}\,=\,\underline{\underline A}_{k,j}^{-1}\cdot \mathbf{b}_{k,j}
\end{equation}
where:
\begin{equation}\label{LSmatrix}
\underline{\underline A}_{k,j}\,=\,\sum_T \, \theta_T \,\mathbf{P}^t(x^T_1,x^T_2) \, \mathbf{P}(x^T_1,x^T_2) \;\;\; \text{and} \;\;\; \mathbf{b}_{k,j}\,=\,\sum_T\, \theta_T \, \mathbf{P}^t(x^T_1,x^T_2) \, \left( \frac{\partial Z_k}{\partial x_j} \right)_{T}
\end{equation}
\revAB{The $GG$ reconstruction, Eq.~(\ref{eq:nabla}), can be recovered 
from the $ZZ$ reconstruction~(\ref{eq:ZZ}) by choosing}
a constant polynomial expansion: $\mathbf{P}\,=\,[\,1\,,\,0\,,\,0\,]$.\\ 
        \section*{References}
	\bibliographystyle{elsarticle-num} 
	\bibliography{biblio}
	\end{document}